\documentclass[12pt]{article}
\usepackage{amsmath,amssymb,amsthm}
\usepackage[english]{babel}
\usepackage[cp1251]{inputenc}
\def\leq {\leqslant}
\def\le {\leqslant}
\def\ge {\geqslant}
\def\geq {\geqslant}
\def\@bibitem[#1]#2{\item\@biblabel{#1}.\if@filesw
{\def\protect##1{\string##1\space}\immediate\write
\@auxout{\string\bibcite{#2}{#1}}}\fi\ignorespaces\@showtag{#2}}

\theoremstyle{plain}
\newtheorem{theorem}{Theorem}[section]
\newtheorem{rem}{Remark}[section]
\newtheorem{lemma}{Lemma}[section]

\renewcommand{\theequation}%
{\arabic{section}.\arabic{equation}}
\begin{document}
\pagestyle{myheadings}
\thispagestyle{empty}


\begin{center}
{\sc G. Akishev}
\end{center}

\begin{center}
{\bf Lipschitz space with mixed logarithmic smoothness and embedding theorems}
\end{center}

\vspace*{0.2 cm}

\begin{abstract}
\begin{quote}
\noindent{\bf Abstract.}
This article considers the Lipschitz space with mixed logarithmic smoothness 
 of $2\pi$ periodic functions of several variables. 
We obtain equivalent descriptions of the norm of the Lipschitz space and prove
embedding theorems between Besov and Lipschitz spaces.
\end{quote}
\end{abstract}
\vspace*{0.2 cm}


\textbf{MCS}: 26A16, 26B35, 42B35. 

\vspace*{0.2 cm}

\section*{Introduction} 
\subsection{Lorentz space}\label{subsec0.1}
 
 Let $\mathbb{R}^{m}$ be an $m$--dimensional Euclidean space of points $\mathbf{x}=(x_{1},\ldots,x_{m})$ with real coordinates. Let $\mathbb{T}^{m}=\{\mathbf{x}=(x_{1},\ldots,x_{m})\in \mathbb{R}^{m};\  0\leq x_{j}< 2\pi;\  j=1,\ldots,m\}$ be an $m$--dimensional cube, and  $\mathbb{I}^{m}=[0, 1)^{m}$.
 \smallskip

We denote by $L_{p,\tau}(\mathbb{T}^{m})$,  $1<
p<\infty, 1\leqslant \tau <\infty$, 
the Lorentz space of all real-valued Lebesgue measurable functions $f$,  $2\pi$--periodic in each variable, and for which the quantity
\begin{equation*}
\|f\|_{p,\tau} = \left\{\int\limits_{0}^{1}\bigl(f^{*}(t)
\bigr)^{\tau}t^{\frac{\tau}{p}-1}dt
\right\}^{\frac{1}{\tau}} 
\end{equation*}
is finite. Here $f^{*}(t)$ is a non-increasing rearrangement of the function $|f(2\pi\mathbf{x})|$, $\mathbf{x} \in \mathbb{I}^{m }$ (see \cite[Chapter~1, Sec. 3, P.~213--216]{34}.
If $\tau=p$, then $\|f\|_{p}=\|f\|_{p,p}$.

We denote by ${\mathring L}_{p, \tau}(\mathbb{T}^{m})$  the set of all functions $f\in
L_{p, \tau}(\mathbb{T}^{m})$ such that
\begin{equation*}
\int\limits_{0}^{2\pi }f(\mathbf{x}) dx_{j}  =0,\; \;
j=1,...,m .
\end{equation*}
Below $a_{\mathbf{n}}(f)$ are the Fourier coefficients of a function $f\in {\mathring L}_{1}\left(\mathbb{T}^{m} \right)$ with respect  the system $\{e^{i\langle\mathbf{n}, 2\pi\mathbf{x}\rangle}\}_{\mathbf{n}\in\mathbb{Z}^{m}}$ and $\langle\mathbf{y}, \mathbf{x}\rangle=\sum\limits_{j=1}^{m}y_{j} x_{j}$;
 \begin{equation*}
    \delta _{\mathbf{s}}(f, 2\pi\mathbf{x})
=\sum\limits_{\mathbf{n} \in \rho (\mathbf{s})
}a_{\mathbf{n}}(f) e^{i\langle\mathbf{n}, 2\pi\mathbf{x}\rangle},
 \end{equation*}
 \begin{equation*}
\rho(\mathbf{s})=\left\{\mathbf{k} =(k_{1},...,k_{m}) \in \Bbb{Z}^{m}: \, \,  2^{s_{j} -1} \leq \left| k_{j} \right|<2^{s_{j} } ,j=1,...,m\right\},
\end{equation*}
where $s_{j} = 1,2,... .$
Further, $\mathbb{Z}_{+}^{m}$ is the set of points with non-negative integer coordinates.

The value
\begin{equation*}
Y_{l_{1},\ldots,l_{m}}(f)_{p, \tau} = \inf_{T_{l_{j}}} \|f-\sum_{j=1}^{m}T_{l_{j}}\|_{p, \tau}^{*} \;\;,\;\; l_{j} = 0,1,2,...
\end{equation*}
is called the best approximation by ''angle'' of a function $f\in {\mathring L}_{p, \tau}(\mathbb{T}^{m})$ by trigonometric polynomials, where $T_{l_{j}} \in {\mathring L}_{p, \tau}(\mathbb{T}^{m})$ is a trigonometric polynomial of order $l_{j}$ with respect to the variable $x_{j}, \,\, j = 1,\ldots,m$ (in the case $\tau=p$, see  \cite{29}--\cite{32}).

By $C(p,q,r,y)$ we denote positive constants depending on
the parameters indicated in parentheses, generally speaking, different in different
formulas. For positive values $A(y), B(y)$, $A(y) \asymp
B(y)$ means that there are positive numbers $C_{1},\,C_{2}$ such that $C_{1}A(y) \leq B(y) \leq C_{2}A(y)$. For the sake of brevity, in the case $B\ge C_{1}A$ or $B\le C_{2}A$, we often write $B \gg A$ or $B \ll A$, respectively.

\subsection{Smoothness module of a function}\label{subsec0.2}

We consider the standard basis $\{\mathbf{e}_{k}\}_{k=1}^{m}$ in $\mathbb{R}^{m}$.

{\bf Definition} (see, for example, \cite{23},  \cite{32} and the bibliography therein).
Let $\alpha\in (0, \infty)$. For a function $f \in L(\mathbb{T}^{m})$, the positive order difference $\alpha$ with respect to the variable $x_{k}$ with step $h \in \mathbb{R}$ is given  by the formula
\begin{equation*}
\Delta_{h}^{\alpha}f(\mathbf{x})=\sum_{\nu=0}^{\infty}(-1)^{\nu}
\left(\begin{smallmatrix}\alpha \\ \nu \end{smallmatrix}\right) f(\mathbf{x}+(\alpha - \nu)h\mathbf{e}_{k}),
\end{equation*}
where $\left(\begin{smallmatrix}\alpha \\ \nu \end{smallmatrix}\right)=1$ for $\nu=0$, $\left(\begin{smallmatrix}\alpha \\ \nu \end{smallmatrix}\right)=\alpha$ for $\nu=\alpha$ and $\left(\begin{smallmatrix}\alpha \\ \nu \end{smallmatrix}\right)=\frac{\alpha(\alpha -1)\cdot \ldots \cdot (\alpha-\nu+1)}{\nu !}$ for $\nu\geqslant\alpha$.

The mixed difference of positive orders $\alpha_{k}$ with respect to the variable $x_{k}$ of a function $f \in L(\mathbb{T}^{m})$ is defined by induction as $\Delta_{h_{m}}^{\alpha_{m}}(\Delta_{h_{m-1}}^{\alpha_{m-1}}(...(\Delta_{h_{1}}^{\alpha_{1}}f(\mathbf{x}))...))$
and denoted by the symbol $\Delta_{\mathbf{h}}^{\boldsymbol{\alpha}}f(\mathbf{x})$.

{\bf Definition}. The mixed modulus of smoothness of order $\overline\alpha$ of a function $f\in L_{p, \tau}(\mathbb{T}^{m})$ is defined by the formula
\begin{equation*}
\omega_{\boldsymbol{\alpha}}(f, \mathbf{t})_{p, \tau}:=\omega_{\alpha_{1},...\alpha_{m}}(f, t_{1},..., t_{m})_{p, \tau}=\sup_{|h_{1}|\leq t_{1},...,|h_{m}|\leq t_{m}}\|\Delta_{\mathbf{h}}^{\boldsymbol{\alpha}} (f)\|_{p, \tau}.
\end{equation*}

\begin{rem}
In the case of $\tau=p$, the mixed smoothness modulus of a function $f\in L_{p}(\mathbb{T}^{m})$ is defined and its properties are studied in  \cite{23}, \cite{32},  as well as in \cite[Ch. 1, Sec. 11]{7} (and see the bibliography therein).
\end{rem}

{\bf Definition}. The full smoothness modulus of order $\alpha>0$ of a function $f \in L_{p, \tau}(\mathbb{T}^{m})$ is defined as  
\begin{equation*}
\omega_{\alpha}(f, t)_{p, \tau} = \sup_{\|\mathbf{h}\|\leq t}\|\Delta_{\mathbf{h}}^{\alpha}f\|_{p, \tau},
\end{equation*}
where $\mathbf{h}\in \mathbb{R}$ and 
\begin{equation*}
\Delta_{\mathbf{h}}^{\alpha}f(\mathbf{x})=\sum_{\nu=0}^{\infty}(-1)^{\nu}
\left(\begin{smallmatrix}\alpha \\ \nu \end{smallmatrix}\right) f(\mathbf{x}+(\alpha - \nu)\mathbf{h}), \, \, \|\mathbf{h}\|:=\biggl(\sum_{j=1}^{m}h_{j}^{2}\biggr)^{1/2}.
\end{equation*}
In the case of $\tau=p$ and $\alpha\in \mathbb{N}$, the complete smoothness module is defined and its properties are studied   in \cite[Ch. 4, Sec. 4. 2]{27}, \cite[Ch. 3, Sec. 3.2]{36}, \cite{35}, and for  $\alpha>0$ in \cite{23} and the bibliography therein.
  \smallskip

\subsection{Lipschitz spaces with mixed smoothness modulus}\label{subsec0.3}

Define  the Lipschitz  class consisting of all $2\pi$ periodic
functions $f$ satisfying the condition
\begin{equation*}
|f(x)-f(y)|\leqslant M|x-y|^{\alpha}
\end{equation*}
for some positive number M (see, for example, \cite[Ch.3, Sec. 2]{17}, \cite{19},  \cite[Ch. 3,  Sec. 3.2]{36}).  

The first embedding theorems for the class $\text{Lip}(\alpha, p)$ were obtained by E.K. Titchmarsh \cite{37}. In the article \cite{19}, the Lipschitz class $\text{Lip}(\alpha, p)$ is defined in the space $L_{p}(\mathbb{T})$, $1\leqslant p \leqslant \infty$, $\alpha\in [0, 1)$, and embedding theorems for Lipschitz spaces are proved.   This topic was developed by P.L. Ul'yanov \cite{38}, and a more detailed review was given by V.I. Kolyada \cite{25}.  The results of P.L. Ul'yanov were further developed in joint works by Yu.S. Kolomoitsev and S.Yu. Tikhonov \cite{22}, Z. Ditzian and S.Yu. Tikhonov \cite{13} (see also the bibliography in these works). 
 
The article \cite{14}, the logarithmic Lipschitz space $\text{Lip}_{p, q}^{(\alpha, -b)}(\mathbb{ T}^{m})$ is defined as the family of all functions $f\in L_{p}(\mathbb{ T}^{m})$ such that 
 \begin{equation*}
\|f\|_{\text{Lip}_{p, q}^{(\alpha, -b)}}: =\|f\|_{p} + \biggl(\int_{0}^{1}\Bigl( t^{-\alpha}(1-\log t)^{-b}\omega_{\alpha}(f, t)_{p}\Bigr)^{q}\frac{dt}{t} \biggr)^{1/q} <+\infty,
\end{equation*}
for numbers $\alpha> 0$, $1\leqslant p\leqslant \infty$, $0< q \leqslant \infty$, $b\in \mathbb{R}$. It is assumed that $b> 1/q$ ($b\geqslant 0$ if $q=\infty$). Otherwise, $L_{p, q}^{(\alpha, -b)}(\mathbb{ T}^{m})=\{0\}$.
Setting $\alpha=1$, we obtain the spaces studied in \cite{11}, \cite{12}, \cite{21}.

In the articles \cite{11}, \cite{12}, \cite{18} the following generalization of the Besov space is defined:  
$\mathbf{B}_{p, \theta}^{r, b}(\mathbb{T}^{m})$ is the set of all functions $f\in L_{p}(\mathbb{T}^{m})$, $1\leq p < \infty$ for which
\begin{equation*}
\|f\|_{\mathbf{B}_{p, \theta}^{r, b}}: = \|f\|_{p} + \Bigl(\int_{0}^{1}(t^{-r}(1 - \log t)^{b}\omega_{k}(f, t)_{p})^{\theta} \frac{dt}{t}\Bigr)^{1/\theta} < \infty, \, \, 
\end{equation*}
for $0< \theta \leq \infty$, $b> -1/\theta$, $k\in \mathbb{N}$, $k > r>0$.
A more general  Besov space  in symmetric space was defined by M.Z. Berkolaiko \cite{9}.

In the case of $\alpha=1$, the space $\text{Lip}_{p, q}^{(1, -b)}(\mathbb{ T}^{m})$ is defined and studied by D.D. Haroske \cite{21}. The relationships between the spaces   $\text{Lip}_{p, q}^{(1, -b)}(\mathbb{ T}^{m})$ and the Besov spaces $B_{p, q}^{1, b}(\mathbb{ T}^{m})$ were investigated by F. Cobos, O. Dominguez \cite{11}  and F. Cobos, O. Dominguez, and Triebel H. \cite{12}. In the general case $\alpha > 0$ , the problems of embedding one Lipschitz space into another such space and into Besov spaces were studied in detail by O. Dominguez, D.D. Haroske, S. Tikhonov \cite{14}, O. Dominguez, S. Tikhonov \cite{15}, Yu.S. Kolomoitsev, and S.Yu. Tikhonov \cite{22}.

If $p=q=\infty$, $\alpha=1$ and $b=0$ in $\text{Lip}_{p, q}^{(1, -b)}(\mathbb{ T}^{m})$, then we obtain the classical Lipschitz space   $\text {Lip} 1 (\mathbb{ T}^{m})$ (see \cite{14}, \cite{37})

 We consider the following well-known generalized Nikol'skii--Besov space (see, for example, \cite{30}Ц-\cite{32}).

{\bf Definition.} Let $\boldsymbol{\alpha} = (\alpha_{1},...,\alpha_{m})$, $\mathbf{r} = (r_{1},...,r_{m}),$ $\alpha_{j}>r_{j} > 0$, $j =1,...,m,$ $1\leq p< \infty$, $1\leq \tau < \infty$, $0< \theta \leq \infty$,  $b_{j}\in \mathbb{R}$, $j=1,...,m$. By $S_{p, \tau, \theta}^{\mathbf{r}, \mathbf{b}}\mathbf{B}$ we denote the space of all functions $f\in {\mathring L}_{p, \tau} (\mathbb{T}^{m})$ for which
\begin{equation*}
\| f\|_{S_{p, \tau, \theta}^{\mathbf{r}, \mathbf{b}}\mathbf{B}}
 = \|f\|_{p, \tau} + \Biggl[\int_{0}^{1}...\int_{0}^{1}\omega_{\boldsymbol{\alpha}}^{\theta}(f, \mathbf{t})_{p, \tau}\prod_{j=1}^{m} \frac{(1-\log t_{j})^{\theta b_{j}}}{t_{j}^{1+\theta r_{j}}}dt_{1}...dt_{m}\Biggr]^{\frac{1}{\theta}} < +\infty,
\end{equation*}
where $\mathbf{b}=(b_{1},...,b_{m})$, $\mathbf{k}=(k_{1},...,k_{m})$, $k_{j}\in \mathbb{N}$, $j=1,...,m$.

 It is known that (in particular, in the case of $\tau=p$ and $b_{j}=0$, $j =1,...,m,$ see \cite{26})
 \begin{equation}\label{eq0.1}
\| f\|_{S_{p, \tau, \theta}^{\mathbf{r}, \mathbf{b}}B} \asymp
\left\{\sum\limits_{\mathbf{s} \in \Bbb{Z}_{+}^{m}} \prod_{j=1}^{m}2^{s_{j}r_{j}\theta}(s_{j} + 1)^{b_{j}\theta} \|\delta_{\mathbf{s}}(f)\|_{p, \tau}^{\theta}\right\}^{\frac{1}{\theta}},
 \end{equation}
for  $1 < p < +\infty,$ $1\leq \tau< \infty$, $0< \theta \le +\infty$.

In this article, we consider $\text{Lip}_{p, \tau, \theta}^{(\boldsymbol{\alpha}, -\mathbf{b})}(\mathbb{T}^{m})$ a logarithmic Lipschitz space with mixed smoothness in the Lorentz space.

{\bf Definition.} Let $\alpha_{j}> 0$,  $b_{j}\in \mathbb{R}$ for $j=1,\ldots ,m$ and 
$\boldsymbol{\alpha}=(\alpha_{1},\ldots ,\alpha_{m})$, $\mathbf{b}=(b_{1},\ldots ,b_{m})$, $1<p<\infty$, $1<\tau<\infty$, $0<\theta \leqslant \infty$. A logarithmic Lipschitz space with mixed smoothness is the space of all functions $f\in L_{p, \tau}(\mathbb{ T}^{m})$ for which
\begin{equation*}
\|f\|_{\text{Lip}_{p, \tau, \theta}^{(\boldsymbol{\alpha}, -\mathbf{b})}}:=\|f\|_{p, \tau} + \Biggl[\int_{0}^{1}...\int_{0}^{1}\omega_{\boldsymbol{\alpha}}^{\theta}(f, \mathbf{t})_{p, \tau}\Bigl(\prod_{j=1}^{m}t_{j}^{-\alpha_{j}} (1-\log t_{j})^{- b_{j}}\Bigr)^{\theta} \prod_{j=1}^{m}\frac{dt_{j}}{t_{j}}\Biggr]^{\frac{1}{\theta}} < +\infty
\end{equation*}
and is denoted by the symbol $\text{Lip}_{p, \tau, \theta}^{(\boldsymbol{\alpha}, -\mathbf{b})}(\mathbb{T}^{m})$. In this case, it is assumed that $b_{j}> 1/\theta$ ($b_{j}\geqslant 0$, if $\theta=\infty$). Otherwise, $\text{Lip}_{p, \tau, \theta}^{(\boldsymbol{\alpha}, -\mathbf{b})}(\mathbb{T}^{m})=\{0\}$, since if $b_{j_{0}}\leqslant 1/\theta$ for some $j_{0}$, then 
 \begin{equation*}
\int_{0}^{1}\omega_{\boldsymbol{\alpha}}^{\theta}(f, \mathbf{t})_{p, \tau}\Bigl(t_{j_{0}}^{-\alpha_{j_{0}}} (1-\log t_{j_{0}})^{- b_{j_{0}}}\Bigr)^{\theta}\frac{dt_{j_{0}}}{t_{j_{0}}} =+\infty.
\end{equation*}

The main goal of the article is to find equivalent norms of the space  $\text{Lip}_{p, \tau, \theta}^{(\boldsymbol{\alpha}, -\mathbf{b})}(\mathbb{T}^{m})$ and prove embedding theorems between the spaces $\text{Lip}_{p, \tau, \theta}^{(\boldsymbol{\alpha}, -\mathbf{b})}(\mathbb{T}^{m})$ and $S_{p, \tau, \theta}^{\mathbf{r}, \mathbf{b}}B$. 
In the first section, we formulate the main results of the article and compares them with known theorems. In the second section we prove the auxiliary statements necessary for proving the main results. In the third section we provides proofs of the main results.

\setcounter{equation}{0}
\setcounter{lemma}{0}
\setcounter{theorem}{0}

\section{Main results}\label{sec1}

The following theorem gives equivalent norm characterizations in a Lipschitz space.
 \begin{theorem}\label{th1.1} 
Let $\alpha_{j}> 0$,  $b_{j}>1/\theta$ for $j=1,\ldots ,m$ and $1<p<\infty$, $1<\tau<\infty$, $0<\theta \leqslant \infty$. Then  for the function $f\in \text{Lip}_{p, \tau, \theta}^{(\boldsymbol{\alpha}, -\mathbf{b})}(\mathbb{T}^{m})$, the following relation holds
\begin{equation*}
\|f\|_{\text{Lip}_{p, \tau, \theta}^{(\boldsymbol{\alpha}, -\mathbf{b})}}\asymp \|f\|_{p, \tau} + (J_{p, \tau, \theta, \boldsymbol{\alpha}}(f))^{\frac{1}{\theta}}. 
\end{equation*}
Here and further   
\begin{equation*}
J_{p, \tau, \theta, \boldsymbol{\alpha}}(f):=\sum\limits_{l_{m} =0}^{\infty}...\sum\limits_{l_{1} =0}^{\infty}
\prod_{j=1}^{m}2^{l_{j}(\frac{1}{\theta}-b_{j})\theta}
\Biggl\|\Biggl(\sum\limits_{s_{m}=2^{l_{m}}}^{2^{l_{m}+1}-1}...\sum\limits_{s_{1}=2^{l_{1}}}^{2^{l_{1}+1}-1}\prod_{j=1}^{m}2^{s_{j}\alpha_{j}2}|\delta_{\mathbf{s}}(f)|^{2}\Biggr)^{\frac{1}{2}}\Biggr\|_{p, \tau}^{\theta}
\end{equation*}
 \end{theorem}
We set $\mathbf{e}=\sum_{\nu=1}^{m}\mathbf{e}_{\nu}=(1,...,1)$, where $\{\mathbf{e}_{\nu}\}_{\nu=1}^{m}$ is the standard basis in $\mathbb{R}^{m}$. 
 \begin{theorem}\label{th1.2}
Let $\alpha_{j}> 0$,  $b_{j}>1/\theta$ for $j=1,\ldots ,m$ and $\boldsymbol{\alpha}=(\alpha_{1},\ldots ,\alpha_{m})$, $\mathbf{b}=(b_{1},\ldots ,b_{m})$, $1<p<\infty$, $1<\tau<\infty$, $0<\theta \leqslant \infty$. Then 
 \begin{equation}\label{eq1.1} 
S_{p, \tau, \theta}^{\boldsymbol{\alpha}, -\mathbf{b}+\frac{1}{\min\{2, \tau, \theta\}}\mathbf{e}}B(\mathbb{T}^{m}) \hookrightarrow \text{Lip}_{p, \tau, \theta}^{(\boldsymbol{\alpha}, -\mathbf{b})}(\mathbb{T}^{m})\hookrightarrow S_{p, \tau, \theta}^{\boldsymbol{\alpha}, -\mathbf{b}+\frac{1}{\max\{2, \tau, \theta\}}\mathbf{e}}B(\mathbb{T}^{m}), 
\end{equation}
\begin{equation}\label{eq1.2} 
S_{p, \tau, \min\{2, \tau, \theta\}}^{\boldsymbol{\alpha}, -\mathbf{b}+\frac{1}{\theta}\mathbf{e}}B(\mathbb{T}^{m}) \hookrightarrow \text{Lip}_{p, \tau, \theta}^{(\boldsymbol{\alpha}, -\mathbf{b})}(\mathbb{T}^{m})\hookrightarrow S_{p, \tau, \max\{2, \tau, \theta\}}^{\boldsymbol{\alpha}, -\mathbf{b}+\frac{1}{\theta}\mathbf{e}}B(\mathbb{T}^{m}). 
\end{equation}
 \end{theorem}
 \begin{rem}\label{rem1.1}
  Theorem 1.2 is an analogue of Theorem 4.1 in \cite{14}.
 \end{rem}

 \begin{theorem}\label{th1.3} 
Let $\alpha_{j}> 0$,  $b_{j}>1/\theta$ for $j=1,\ldots ,m$ and
$\boldsymbol{\alpha}=(\alpha_{1},\ldots ,\alpha_{m})$, $\mathbf{b}=(b_{1},\ldots ,b_{m})$, $1<p_{0}<p<p_{1}<\infty$, $1<\tau_{0}, \tau_{1}, \tau<\infty$, $0<\theta \leqslant \infty$.
Then  
 \begin{equation}\label{eq1.3} 
S_{p_{0}, \tau_{0}, \theta}^{\boldsymbol{\alpha}+(\frac{1}{p_{0}}-\frac{1}{p})\mathbf{e}, -\mathbf{b}+\frac{1}{\min\{\tau, \theta\}}\mathbf{e}}B(\mathbb{T}^{m}) \hookrightarrow \text{Lip}_{p, \tau, \theta}^{(\boldsymbol{\alpha}, -\mathbf{b})}(\mathbb{T}^{m})\hookrightarrow S_{p_{1}, \tau_{1}, \theta}^{\boldsymbol{\alpha}+(\frac{1}{p_{1}}-\frac{1}{p})\mathbf{e}, -\mathbf{b}+\frac{1}{\max\{\tau, \theta\}}\mathbf{e}}B(\mathbb{T}^{m}), 
\end{equation}
\begin{equation}\label{eq1.4} 
S_{p_{0}, \tau_{0}, \min\{\tau, \theta\}}^{\boldsymbol{\alpha}+(\frac{1}{p_{0}}-\frac{1}{p})\mathbf{e}, -\mathbf{b}+\frac{1}{\theta}\mathbf{e}}B(\mathbb{T}^{m}) \hookrightarrow \text{Lip}_{p, \tau, \theta}^{(\boldsymbol{\alpha}, -\mathbf{b})}(\mathbb{T}^{m})\hookrightarrow S_{p_{1}, \tau_{1}, \max\{\tau, \theta\}}^{\boldsymbol{\alpha}, -\mathbf{b}+\frac{1}{\theta}\mathbf{e}}B(\mathbb{T}^{m}). 
\end{equation}
 \end{theorem}
\begin{rem}\label{rem1.2} Theorem 1.3 is analogous to Theorem 4.4 in \cite{14}.
  \end{rem}
  Next, we consider Lipschitz spaces $\text{Lip}_{p, \tau, \theta}^{(\boldsymbol{\alpha}, -\mathbf{b})}(\mathbb{T}^{m})$ for fixed $p, \tau$.

 \begin{theorem}\label{th1.4} 
Let $\alpha_{j}^{(i)}>0$, $b_{j}^{(i)}>\frac{1}{\theta_{i}}$ for $j=1,\ldots , m$, $0<\theta_{i}\leqslant \infty$, $\boldsymbol{\alpha}^{(i)}=(\alpha_{1}^{(i)}, \ldots, \alpha_{m}^{(i)})$, $\mathbf{b}^{(i)}=(b_{1}^{(i)}, \ldots , b_{m}^{(i)})$, $i=0, 1$ and $1<p<\infty$, $1<\tau<\infty$.  Let us assume that one of the following conditions holds:

i) $\alpha_{j}^{(0)}>\alpha_{j}^{(1)}$ for $j=1,\ldots , m$, $0<\theta_{i}\leqslant \infty$, $i=0, 1$;

ii)  $\alpha_{j}^{(0)}=\alpha_{j}^{(1)}$, $b_{j}^{(1)}-\frac{1}{\theta_{1}}> b_{j}^{(0)}-\frac{1}{\theta_{0}}$ for $j=1,\ldots, m$;

iii) $\alpha_{j}^{(0)}=\alpha_{j}^{(1)}$, $b_{j}^{(1)}-\frac{1}{\theta_{1}}= b_{j}^{(0)}-\frac{1}{\theta_{0}}$ for $j=1,\ldots, m$.
Then 
 \begin{equation*}
 \text{Lip}_{p, \tau, \theta_{0}}^{(\boldsymbol{\alpha}^{(0)}, -\mathbf{b}^{(0)})}(\mathbb{T}^{m}) \hookrightarrow \text{Lip}_{p, \tau, \theta_{1}}^{(\boldsymbol{\alpha}^{(1)}, -\mathbf{b}^{(1)})}(\mathbb{T}^{m}).
 \end{equation*}
\end{theorem}
 \begin{rem}\label{rem1.3}
  Theorem 1.4 is analogous to Theorem 5.1 in \cite{14}.
  \end{rem}
Now, let us consider the case $b_{j}^{(0)}=b_{j}^{(1)}=b_{j}$ for $j=1,\ldots, m$.
 \begin{theorem}\label{th1.5} 
Let $0<p_{0}<p_{1}<\infty$, $1\leqslant\tau_{0}, \tau_{1}< \infty$, $0< \theta\leqslant \infty$,  $0<\alpha_{j}^{(0)}<\alpha_{j}^{(1)}<\infty$, $\alpha_{j}^{(0)}-\frac{1}{p_{0}}=\alpha_{j}^{(1)}-\frac{1}{p_{1}}$, $b_{j}> \frac{1}{\theta}$ for $j=1,\ldots, m$. Then  
 \begin{equation*}
 \text{Lip}_{p_{0}, \tau_{0}, \theta}^{(\boldsymbol{\alpha}^{(0)}, -\mathbf{b})}(\mathbb{T}^{m})
\hookrightarrow \text{Lip}_{p_{1}, \tau_{1}, \theta}^{(\boldsymbol{\alpha}^{(1)}, -\mathbf{b})}(\mathbb{T}^{m}).
 \end{equation*}
\end{theorem}
 \begin{rem}\label{rem1.4}
Theorem 1.5 is analogous to Theorem 5.3 in \cite{14}. 
  \end{rem}
 Next, consider the standard basis $\{\mathbf{e}_{\nu}\}_{\nu =1}^{m}$ in $\mathbb{R}^{m}$.
For brevity, let $v_{j}=-b_{j}+\frac{1}{\min\{\tau, \theta\}}$, $\beta_{j}=\alpha_{j}+\frac{1}{p_{0}}-\frac{1}{p}$ for  $j=1,\ldots, m$ and $\mathbf{v}=(v_{1}, \ldots ,v_{m})$, $\boldsymbol{\beta}=(\beta_{1}, \ldots ,\beta_{m})$.
Let us prove the optimality of the embeddings in Theorem 1.2, Theorem 1.3, Theorem 1.4.
\begin{theorem}\label{th1.6} 
Let $1<p_{0}<p<p_{1}<\infty$, $1<\tau_{0}, \tau_{1}, \tau<\infty$, $0<\theta\leqslant \infty$ and $\alpha_{j}> 0$, $b_{j}>\frac{1}{\theta}$ for  $j=1,\ldots, m$.

1. If $\min\{\tau, \theta\}=\tau$, then for any number $\varepsilon >0$ there exists a function $f_{1}\in S_{p_{0}, \tau_{0}, \tau}^{\boldsymbol{\beta}, \mathbf{v}-\varepsilon\mathbf{e}_{j_{0}}}B(\mathbb{T}^{m})$ such that $f_{1}\notin \text{Lip}_{p, \tau, \theta}^{(\boldsymbol{\alpha}, -\mathbf{b})}(\mathbb{T}^{m})$.

2. If $\max\{\tau, \theta\}=\tau$, then for any number $\varepsilon >0$ there exists a function $f_{2}\in \text{Lip}_{p, \tau, \theta}^{(\boldsymbol{\alpha}, -\mathbf{b})}(\mathbb{T}^{m})$ such that $f_{2}\notin S_{p_{1}, \tau_{1},  \tau}^{\mathbf{u}, \boldsymbol{\mu}+\varepsilon\mathbf{e}_{j_{0}}}B(\mathbb{T}^{m})$, where  $\mathbf{u}=(u_{1}, \ldots ,u_{m})$, $\boldsymbol{\mu}=(\mu_{1}, \ldots ,\mu_{m})$, $u_{j}=\alpha_{j}+\frac{1}{p_{1}}-\frac{1}{p}$, $\mu_{j}=-b_{j}+\frac{1}{\max\{\tau, \theta\}}$ for  $j=1,\ldots, m$.
\end{theorem}
\begin{rem}\label{rem1.5}
Theorem 1.6 shows the optimality of the embeddings \eqref{eq1.3}.
 \end{rem}
Let us set $v_{j}=-b_{j}+\frac{1}{\theta}$, $\beta_{j}=\alpha_{j}+\frac{1}{p_{0}}-\frac{1}{p}$ for $j=1,\ldots, m$ and $\mathbf{v}=(v_{1}, \ldots ,v_{m})$, $\boldsymbol{\beta}=(\beta_{1}, \ldots ,\beta_{m})$.
\begin{theorem}\label{th1.7} 
Let $1<p_{0}<p<p_{1}<\infty$, $1<\tau_{0}, \tau_{1},  \tau<\infty$, $0<r, \theta\leqslant \infty$ and $\alpha_{j}> 0$, $b_{j}>\frac{1}{\theta}$ for  $j=1,\ldots, m$. Then

i) the condition $r\leqslant \min\{\tau, \theta\}$ is sufficient and in the case $\min\{\tau, \theta\}=\tau$ is necessary for embedding
\begin{equation}\label{eq1.5} 
S_{p_{0}, \tau_{0},  r}^{\boldsymbol{\beta}, \mathbf{v}}B(\mathbb{T}^{m}) \hookrightarrow \text{Lip}_{p, \tau, \theta}^{(\boldsymbol{\alpha}, -\mathbf{b})}(\mathbb{T}^{m});
\end{equation}

ii) the condition $r\geqslant \max\{\tau, \theta\}$ is sufficient and in the case $\max\{\tau, \theta\}=\tau$ is necessary for embedding 
\begin{equation}\label{eq1.6} 
\text{Lip}_{p, \tau, \theta}^{(\boldsymbol{\alpha}, -\mathbf{b})}(\mathbb{T}^{m}) \hookrightarrow
S_{p_{1}, \tau_{1},  r}^{\mathbf{u}, \boldsymbol{\gamma}}B(\mathbb{T}^{m});
\end{equation}
where $\mathbf{u}=(u_{1}, \ldots ,u_{m})$, $\boldsymbol{\gamma}=(\gamma_{1}, \ldots ,\gamma_{m})$, $u_{j}=\alpha_{j}+\frac{1}{p_{1}}-\frac{1}{p}$, $\gamma_{j}=-b_{j}+\frac{1}{\theta}$ for  $j=1,\ldots, m$.
\end{theorem}
\begin{rem}\label{rem1.6}
Theorems 1.6 and 1.7 are analogues of Theorem 9. 1 in \cite{14}.
 Theorem 1.7 shows the optimality of embeddings  \eqref{eq1.4} when $\min\{\tau, \theta\}=\tau$ and $\max\{\tau, \theta\}=\tau$.
 \end{rem}

\begin{theorem}\label{th1.8} 
Let $1<p<\infty$, $1< \tau<\infty$, $0<\theta\leqslant \infty$ and $\alpha_{j}> 0$, $b_{j}>\frac{1}{\theta}$, $\xi_{j}\in \mathbb{R}$ for  $j=1,\ldots, m$ and $\boldsymbol{\xi}=(\xi_{1},\ldots,\xi_{m})$. Then

i) The condition $\min_{j=1,\ldots,m}\{\xi_{j}\}\geqslant \frac{1}{\min\{2, \tau, \theta\}}$ is sufficient and, in the case of $\min\{2, \tau, \theta\}=\tau$, necessary for inclusion 
\begin{equation}\label{eq1.7} 
S_{p, \tau, \theta}^{\boldsymbol{\alpha}, -\mathbf{b}+\boldsymbol{\xi}}B(\mathbb{T}^{m}) \hookrightarrow  \text{Lip}_{p, \tau, \theta}^{(\boldsymbol{\alpha}, -\mathbf{b})}(\mathbb{T}^{m}); 
\end{equation}
ii) the condition $\max_{j=1,\ldots,m}\{\xi_{j}\}\leqslant \frac{1}{\max\{2, \tau, \theta\}}$ is sufficient and in the case $\max\{2, \tau, \theta\}=\tau$ is necessary for embedding 
\begin{equation}\label{eq1.8}  
\text{Lip}_{p, \tau, \theta}^{(\boldsymbol{\alpha}, -\mathbf{b})}(\mathbb{T}^{m}) \hookrightarrow S_{p, \tau, \theta}^{\boldsymbol{\alpha}, -\mathbf{b}+\boldsymbol{\xi}}B(\mathbb{T}^{m});
\end{equation}

iii) The condition $q\leqslant \min\{2, \tau, \theta\}$ is sufficient and, in the case of $\min\{2, \tau, \theta\}=\tau$, necessary for inclusion 
\begin{equation}\label{eq1.9} 
S_{p, \tau, q}^{\boldsymbol{\alpha}, -\mathbf{b}+\frac{1}{\theta}\mathbf{e}}B(\mathbb{T}^{m}) \hookrightarrow   \text{Lip}_{p, \tau, \theta}^{(\boldsymbol{\alpha}, -\mathbf{b})}(\mathbb{T}^{m}); 
\end{equation}

iv) the condition $q\geqslant \max\{2, \tau, \theta\}$ is sufficient and in the case $\max\{2, \tau, \theta\}=\tau$ is necessary for embedding 
\begin{equation}\label{eq1.10} 
\text{Lip}_{p, \tau, \theta}^{(\boldsymbol{\alpha}, -\mathbf{b})}(\mathbb{T}^{m}) \hookrightarrow
S_{p, \tau, q}^{\boldsymbol{\alpha}, -\mathbf{b}+\frac{1}{\theta}\mathbf{e}}B(\mathbb{T}^{m}).     
\end{equation}
\end{theorem}

 \begin{rem}\label{rem1.7}
Theorem 1.8 is an analogue of Theorem 4. 1 in \cite{14}.
Theorem 1.8 shows the optimality of the embeddings \eqref{eq1.1} and \eqref{eq1.2}.
\end{rem}

  \smallskip
\setcounter{equation}{0}
\setcounter{lemma}{0}
\setcounter{theorem}{0}

\section{Auxiliary statements}\label{sec2}

\smallskip
First, we introduce additional notation and give auxiliary statements.
Denote by $e_{m}$ the set of indices $\{1, \dots, m\}$, its arbitrary subset by $e$ and the number of elements of $e$ by $|e|$.
Here $\mathbf{r} =(r_{1},\ldots,r_{m})$ is an element of an $m$--dimensional space with
non-negative coordinates, and $\mathbf{r}^{e} =(r_{1}^{e},\ldots,r_{m}^{e})$
is the vector with components $r_{j}^{e} = r_{j}$ for $j\in e$ and $r_{j}^{e}= 0$ for $j\notin e.$

Let $\mathbf{l} = (l_{1},\ldots,l_{m})$ be an element of an $m$--dimensional space with
positive integer coordinates and a nonempty set $e\subset e_{m}.$
We set
 $$
G_{\mathbf{l}}(e) = \{\mathbf{k} = (k_{1},\ldots,k_{m})\in \mathbb{Z}^{m} :
|k_{j}|\le l_{j}, j\in e \quad |k_{j}| > l_{j}, j\notin e \}.
 $$
We consider various partial sums of Fourier series:
 $$
S_{\mathbf{l}}(f, 2\pi\mathbf{x}) = S_{l_{1},\ldots,l_{m}}(f, 2\pi\mathbf{x}) =
\sum\limits_{|k_{1}|\le l_{1}}\ldots\sum\limits_{|k_{m}|\le l_{m} }
a_{\mathbf{k}} (f) e^{i\langle\mathbf{k} , 2\pi\mathbf{x}\rangle}
 $$
--is the partial sum with respect to all variables;
 $$
S_{l_{1}, \infty}(f, 2\pi\mathbf{x}) =
\sum\limits_{|k_{1}|\le l_{1}}\sum\limits_{k_{2}=-\infty }^{+\infty}\ldots\sum\limits_{k_{m}=-\infty}^{+\infty}a_{\mathbf{k}}(f) e^{i\langle\mathbf{k} , 2\pi\mathbf{x}\rangle}
 $$
--is the partial sum with respect to $ x_{j}\in [0, 1)$.

In the general case,
 $$
S_{\mathbf{{l}^{e}}, \mathbf{\infty}}(f, 2\pi\mathbf{x}) =
\sum\limits_{\mathbf{k} \in \prod_{j \in e} [-l_{j}, l_{j}]\times \mathbb{R}^{m-|e|} }
a_{\mathbf{k}}(f) e^{i\langle\mathbf{k} , 2\pi\mathbf{x}\rangle}
 $$
is the partial sum with respect to $x_{j}\in [0, 1)$ for $j\in e.$

For a given subset of $e\subset e_{m}$, we put
 $$
U_{\mathbf{l}}(f, 2\pi\mathbf{x}) = \sum\limits_{e\subset e_{m}, e \neq \emptyset} \;\; \sum\limits_{\mathbf{k} \in G_{\mathbf{l}}(e)} a_{\mathbf{k}}(f) e^{i\langle\mathbf{k} ,2\pi\mathbf{x}\rangle}.
 $$
In particular, for $m=2$, we have (see, for example, \cite{32})
 $$
U_{l_{1},l_{2}}(f, 2\pi\mathbf{x}) = S_{l_{1} ,\infty }(f, 2\pi\mathbf{x}) + S_{\infty,l_{2}  }(f, 2\pi\mathbf{x}) - S_{l_{1},l_{2}  }(f, 2\pi\mathbf{x}).
 $$
Below we present some properties of the mixed modulus of smoothness of a function, which are proved by well-known methods as in \cite{1}, \cite{13}, \cite{23}, \cite{32}, \cite{35}.

\begin{lemma}\label{lem2 1} (see \cite{6}).
Let $1 < p < +\infty$, $1\leqslant \tau<\infty$, $\alpha_{j}\in (0, \infty)$ for $j=1, \ldots , m$ and functions $f, g \in L_{p, \tau}(\mathbb{T}^{m})$. Then

1. $\omega_{\boldsymbol{\alpha}}(f, \boldsymbol{\delta}^{e})_{p, \tau}=\omega_{\boldsymbol{\alpha}}(f, \mathbf{0})_{p, \tau}=0$,

2. $\omega_{\boldsymbol{\alpha}}(f + g, \boldsymbol{\delta})_{p, \tau}\ll\omega_{\boldsymbol{\alpha}}(f, \boldsymbol{\delta})_{p, \tau} + \omega_{\boldsymbol{\alpha}}(g, \boldsymbol{\delta})_{p, \tau}$;
 
3. $\omega_{\boldsymbol{\alpha}}(f, \boldsymbol{\delta})_{p, \tau}\leqslant \omega_{\boldsymbol{\alpha}}(f, \mathbf{t})_{p, \tau} \, \, \text{for} \, \, 0\leqslant \delta_{j}<t_{j}, \mathbf{t}=(t_{1}, \ldots , t_{m})$; 

4. $\prod_{j=1}^{m}\delta_{j}^{-\alpha_{j}}\omega_{\boldsymbol{\alpha}}(f, \boldsymbol\delta)_{p, \tau}\leqslant \prod_{j=1}^{m}t_{j}^{-\alpha_{j}}\omega_{\boldsymbol{\alpha}}(f, \mathbf{t})_{p, \tau}$ for  $0<t_{j}\leqslant \delta_{j}\leqslant 1$ and  $j=1, \ldots , m$; 

5. $\omega_{\boldsymbol{\alpha}}(f, \lambda_{1}\delta_{1}, ...,\lambda_{m}\delta_{m})_{p, \tau} \ll \prod_{j=1}^{m}\lambda_{j}^{\alpha_{j}}\omega_{\boldsymbol{\alpha}}(f, \delta_{1}, ...,\delta_{m})_{p, \tau}$,
for numbers $\lambda_{j} \geq 1$, $j=1,...,m$;

6. $\omega_{\boldsymbol{\beta}}(f, \boldsymbol{\delta})_{p, \tau}\leqslant \omega_{\boldsymbol{\alpha}}(f, \boldsymbol{\delta})_{p, \tau}$ for $0<\alpha_{j}<\beta_{j}$, $j=1,...,m$, $\boldsymbol{\beta}=(\beta_{1}, \ldots,\beta_{m})$. 
\end{lemma}
Let $n_{j}\in \mathbb{N}$, $j=1,...,m$. We consider a trigonometric polynomial
 \begin{equation}\label{eq2.1}
T_{\mathbf{n}}(2\pi\mathbf{x}) = {\sum\limits_{k_{1}=-n_{1}}^{n_{1}}}' \ldots {\sum\limits_{k_{m}=-n_{m}}^{n_{m}}}'c_{\mathbf{k}}e^{i\langle\mathbf{k} , 2\pi\mathbf{x}\rangle}=\sum\limits_{\mathbf{0}<|\mathbf{k}|\leqslant\mathbf{n}}c_{\mathbf{k}}e^{i\langle\mathbf{k} , 2\pi\mathbf{x}\rangle}, \, \,  \mathbf{x}\in \mathbb{I}^{m}, 
\end{equation} 
where $|\mathbf{k}|=(|k_{1}|,..., |k_{m}|)$ and the notation $\mathbf{0}<|\mathbf{k}|\leqslant\mathbf{n} $ means that $|k_{j}|\leqslant |n_{j}|$ for $j=1,...,m$ and ${\sum_{k=-n}^{n}}'a_{k}$ means that $k\neq 0$. 

\begin{lemma}\label{lem2 2} (see \cite{6}). 
Let $1 < p < +\infty$, $1\leqslant \tau<\infty$,
$\alpha_{j}\in (0, \infty)$ for $j=1, \ldots , m$. Then, for the derivative $T_{\mathbf{n}}^{(\alpha_{1},...,\alpha_{m})}(2\pi\mathbf{x})$ 
of the trigonometric polynomial $T_{\mathbf{n}}(2\pi\mathbf{x})$ of the form \eqref{eq2.1} 
the following inequalities hold:
\begin{equation*}
 \omega_{\boldsymbol{\alpha}}(T_{\mathbf{n}}, \delta_{1}, ...,\delta_{m})_{p, \tau} \ll 
\prod_{j=1}^{m}\delta_{j}^{\alpha_{j}}\|T_{\mathbf{n}}^{(\alpha_{1},...,\alpha_{m})}\|_{p, \tau}.,
 \end{equation*}
  \begin{equation*}
  \|T_{\mathbf{n}}^{(\alpha_{1},...,\alpha_{m})}\|_{p, \tau}\ll \prod_{j=1}^{m}n_{j}^{\alpha_{j}}\|\Delta_{\frac{\pi}{n_{m}}}^{\alpha_{m}}(...(\Delta_{\frac{\pi}{n_{1}}}^{\alpha_{1}}T_{\mathbf{n}})...)\|_{p, \tau}.
   \end{equation*}
 \end{lemma}
For $e\subset e_{m}$, by $T_{\mathbf{n}^{e}, \infty}(\mathbf{x})$ we denote the trigonometric polynomial of order at most $n_{j}\in \mathbb{N}$ in the variable $x_{j}$ for $j\in e$. In the case $e=e_{m}$, the polynomial $T_{\mathbf{n}^{e}, \infty}(\mathbf{x})$ is defined by the equality \eqref{eq2.1}.
\begin{lemma}\label{lem2 3} (see \cite{6}).
Let $1 < p < +\infty$, $1\leqslant \tau<\infty$,
$\alpha_{j}\in (0, \infty)$ for $j\in e$. Then, for the derivative
of the trigonometric polynomial $T_{\mathbf{n}^{e}, \infty}(2\pi\mathbf{x})$,
the inequality holds
 \begin{equation*}
 \|T_{\mathbf{n}^{e}, \infty}^{(\boldsymbol{\alpha}^{e})}\|_{p, \tau}\ll \prod_{j\in e}n_{j}^{\alpha_{j}}\|\Delta_{\mathbf{h}^{e}(\mathbf{n})}^{\boldsymbol{\alpha}^{e}}T_{\mathbf{n}^{e}, \infty}\|_{p, \tau} 
\end{equation*}
for $\mathbf{h}^{e}(\mathbf{n})=(h_{1}^{e}(n_{1}),\ldots , h_{m}^{e}(n_{m}))$, $h_{j}^{e}(n_{j})=\frac{\pi}{n_{j}}$ for $j\in e$.
\end{lemma}
\begin{rem}\label{rem2.1}
In the case $\tau = p$ and $m=2$ Lemma 2.3 was previously proved in \cite{32}, and in the Lebesgue space with mixed norm in \cite{33}.
\end{rem}

\begin{lemma}\label{lem2 4} (Bernstein's inequality). 
Let $1 < p, \tau < +\infty,$ $\alpha_{j}\in \mathbb{Z}_{+}$, for $j=1, \ldots , m$. Then, for the trigonometric polynomial $T_{\mathbf{n}}$, the following inequality holds
\begin{equation*}
\|T_{\mathbf{n}}^{(\alpha_{1},...,\alpha_{m})}\|_{p, \tau}\ll\prod_{j=1}^{m}(n_{j}+1)^{\alpha_{j}} \|T_{\mathbf{n}}\|_{p, \tau}
\end{equation*}
\end{lemma}
 
\begin{lemma}\label{lem2 5} Let $1 < p < +\infty,$\,\,$1 < \tau < +\infty$ and $f\in {\mathring L}_{p, \tau}(\mathbb{T}^{m})$. Then
 \begin{equation*}
\|f - U_{l_{1},...,l_{m}} (f)\|_{p, \tau}\ll Y_{l_{1},...,l_{m}}(f)_{p, \tau}.
\end{equation*}
\end{lemma}

\begin{lemma}(direct and inverse theorems in \cite{5})\label{lem2 6}.    
Let $\alpha_{j}> 0$ for $j=1,\ldots,m$.
If $f\in \mathring{L}_{p, \tau}(\mathbb{T}^{m})$,  $1<p<+\infty,$
$1<\tau <+\infty$, then
 \begin{equation*}
Y_{\mathbf{n}}(f)_{p, \tau}\ll \omega_{\boldsymbol{\alpha}}
\Bigl(f,\frac{1}{n_{1}+1},...,\frac{1}{n_{m}+1}\Bigr)_{p, \tau}\ll\prod_{j=1}^{m}n_{j}^{-\alpha_{j}}\sum\limits_{\nu_{1}=1}^{n_{1}+1} \ldots \sum\limits_{\nu_{m}=1}^{n_{m}+1} \prod_{j=1}^{m}\nu_{j}^{\alpha_{j} - 1}Y_{\boldsymbol{\nu}}(f)_{p, \tau}.
 \end{equation*}
\end{lemma}
\begin{rem}\label{rem2.2}
In the case of $\tau = p$, Lemmas 2.3Ц2.5 are proven in \cite{29} and Lemma 2.1 is proven in \cite{23}.
\end{rem}

\begin{theorem}\label{th2.1} (see in \cite{6}). Let $\alpha_{j}> 0$ for $j=1, \ldots , m$ and $1< p<\infty$, $1< \tau < \infty$. Then, for a function $f\in L_{p, \tau}(\mathbb{T}^{m})$, the following relation holds:
 \begin{equation*}
 \omega_{\boldsymbol{\alpha}}\bigl(f,\frac{\pi}{n_{1}},...,\frac{\pi}{n_{m}}\bigr)_{p, \tau} \asymp 
\|f - U_{\mathbf{n}} (f)\|_{p, \tau} + \sum_{e\subseteq e_{m}, e\neq \emptyset} \prod_{j\in e}n_{j}^{-\alpha_{j}} \|S_{\mathbf{n}^{e}, \infty}^{(\boldsymbol{\alpha}^{e})}(f-S_{\infty, \mathbf{n}^{\hat{e}}}(f))\|_{p, \tau},
  \end{equation*}
where $\hat{e}$ is the complement of the set $e$. 
\end{theorem}

\begin{theorem}\label{th2.2} (see \cite{6}).
Let $\alpha_{j}> 0$ for $j=1, \ldots , m$ and $1< p<\infty$, $1< \tau < \infty$, $\beta=\min\{2, \tau\}$. Then, for a function $f\in L_{p, \tau}(\mathbb{T}^{m})$, the following inequality holds
 \begin{equation*}
 \omega_{\boldsymbol\alpha}\bigl(f,\frac{\pi}{n_{1}},...,\frac{\pi}{n_{m}}\bigr)_{p, \tau} \ll
 \prod_{j=1}^{m}n_{j}^{-\alpha_{j}}\Biggl(\sum\limits_{\nu_{1}=1}^{n_{1}+1} \ldots \sum\limits_{\nu_{m}=1}^{n_{m}+1} \prod_{j=1}^{m}\nu_{j}^{\beta\alpha_{j} - 1}Y_{\boldsymbol{\nu}}^{\beta}(f)_{p, \tau}\Biggr)^{1/\beta}, \,\, n_{j}\in \mathbb{N}.
 \end{equation*}
 \end{theorem}
\begin{rem} 
In the case of $\tau=p$, $m=2$, Theorem 2.1 and Theorem 2.2 are proven in \cite{32}. For $\tau\neq p$, these theorems are proven in \cite{1}.
\end{rem} 
  \smallskip

\begin{lemma}\label{lem2.7}
Let $\alpha_{j}> 0$,  $b_{j}>1/\theta$ for $j=1,\ldots ,m$ and 
$\boldsymbol{\alpha}=(\alpha_{1},\ldots ,\alpha_{m})$, $\mathbf{b}=(b_{1},\ldots ,b_{m})$, $1<p<\infty$, $1<\tau<\infty$, $0<\theta \leqslant \infty$. Then, for the function $f \in \text{Lip}_{p, \tau, \theta}^{(\boldsymbol{\alpha}, -\mathbf{b})}(\mathbb{T}^{m})$, the following relation holds
 \begin{equation*}
\|f\|_{\text{Lip}_{p, \tau, \theta}^{(\boldsymbol{\alpha}, -\mathbf{b})}}\asymp \|f\|_{p, \tau} +\bigl(\Omega_{\boldsymbol{\alpha}}(f)_{p, \tau, \theta} \bigr)^{\frac{1}{\theta}}.
\end{equation*}
Here and further  
 \begin{equation*}
\Omega_{\boldsymbol{\alpha}}(f)_{p, \tau, \theta}:=\sum\limits_{\nu_{m} =0}^{\infty}...\sum\limits_{\nu_{1} =0}^{\infty} \prod_{j=1}^{m}2^{\nu_{j}\alpha_{j}\theta}(\nu_{j} + 1)^{-\theta b_{j}} \omega_{\boldsymbol{\alpha}}^{\theta}\bigl(f,\frac{1}{2^{\nu_{1}}},...,\frac{1}{2^{\nu_{m}}}\bigr)_{p, \tau}
\end{equation*}
\end{lemma}
 \proof
Let $f \in \text{Lip}_{p, \tau, \theta}^{(\boldsymbol{\alpha}, -\mathbf{b})}(\mathbb{T}^{m})$.
Then, taking into account the relation
 \begin{equation*}
\int_{2^{-\nu}}^{2^{-(\nu-1)}} t^{-\alpha_{j}\theta}(1 - \log t)^{-b_j \theta}\frac{dt}{t} \,  \asymp 2^{\nu \alpha_{j} \theta}(\nu + 1)^{-b_j \theta}
\end{equation*}
  by the monotonicity property of the mixed modulus of smoothness, we have
  \begin{multline}\label{eq2.2}
I_{\boldsymbol{\alpha}, p, \tau, \theta}(f):=\int_{0}^{1}...\int_{0}^{1}\omega_{\boldsymbol{\alpha}}^{\theta}(f, \mathbf{t})_{p, \tau}\Bigl(\prod_{j=1}^{m}t_{j}^{-\alpha_{j}} (1-\log t_{j})^{- b_{j}}\Bigr)^{\theta} \prod_{j=1}^{m}\frac{dt_{j}}{t_{j}}
\\
  \ll \sum\limits_{\nu_{m} =1}^{\infty}...\sum\limits_{\nu_{1} =1}^{\infty}\omega_{\boldsymbol{\alpha}}^{\theta}(f, \frac{1}{2^{\nu_{1} - 1}}, \ldots , \frac{1}{2^{\nu_{m} - 1}})_{p, \tau}\int_{\frac{1}{2^{\nu_{1}}}}^{\frac{1}{2^{\nu_{1} - 1}}}...\int_{\frac{1}{2^{\nu_{m}}}}^{\frac{1}{2^{\nu_{m} - 1}}}\Bigl(\prod_{j=1}^{m}t_{j}^{-\alpha_{j}} (1-\log t_{j})^{- b_{j}}\Bigr)^{\theta} \prod_{j=1}^{m}\frac{dt_{j}}{t_{j}} 
  \\
\ll \sum\limits_{\nu_{m} =1}^{\infty}...\sum\limits_{\nu_{1} =1}^{\infty}\prod_{j=1}^{m}2^{\nu_{j} \alpha_{j} \theta}\nu_{j}^{- b_{j}\theta}\omega_{\boldsymbol{\alpha}}^{\theta}(f, \frac{1}{2^{\nu_{1} - 1}}, \ldots , \frac{1}{2^{\nu_{m} - 1}})_{p, \tau}, 
\end{multline}
   \begin{equation}\label{eq2.3}
I_{\boldsymbol{\alpha}, p, \tau, \theta}(f) 
 \gg \sum\limits_{\nu_{m} =1}^{\infty}...\sum\limits_{\nu_{1} =1}^{\infty}\prod_{j=1}^{m}2^{\nu_{j} \alpha_{j} \theta}\nu_{j}^{- b_{j}\theta}\omega_{\boldsymbol{\alpha}}^{\theta}(f, \frac{1}{2^{\nu_{1}}}, \ldots , \frac{1}{2^{\nu_{m}}})_{p, \tau}.
\end{equation}
From the inequalities \eqref{eq2.2} and \eqref{eq2.3} the statement of Lemma 2.7 follows.
\hfill $\Box$

\setcounter{equation}{0}
\setcounter{lemma}{0}
\setcounter{theorem}{0}


\section{Proof of Theorem 1.1}\label{sec3}

Using Theorem 2.1 with $n_{j}=2^{k_{j}}$, $j=1,\ldots ,m$, we have
 \begin{equation}\label{eq3.1}
\omega_{\boldsymbol\alpha}^{\theta}\bigl(f,\frac{1}{2^{\nu_{1}}},...,\frac{1}{2^{\nu_{m}}}\bigr)_{p, \tau} \gg \prod_{j=1}^{m}2^{-\nu_{j}\alpha_{j}}\|S_{2^{\nu_{1}},...,2^{\nu_{m}}}^{(\boldsymbol{\alpha})}(f)\|_{p, \tau}
\end{equation}
for the function $f\in L_{p, \tau}(\mathbb{T}^{m})$, $1<p<\infty$, $1<\tau<\infty$.
By Theorem 2.1 \cite{2}, we have
 \begin{equation*}
\|S_{2^{\nu_{1}},...,2^{\nu_{m}}}^{(\boldsymbol{\alpha})}(f)\|_{p, \tau}\asymp \Biggl\|\Biggl(\sum\limits_{s_{1}=1}^{\nu_{1}}...\sum\limits_{s_{m}=1}^{\nu_{m}}\prod_{j=1}^{m}2^{s_{j}\alpha_{j}2}|\delta_{\mathbf{s}}(f)|^{2}\Biggr)^{\frac{1}{2}}\Biggr\|_{p, \tau}.
\end{equation*}
Therefore, taking into account the boundedness of the partial sum operator in the Lorentz space 
$L_{p, \tau}(\mathbb{T}^{m})$, $1<p<\infty$, $1<\tau<\infty$ from \eqref{eq3.1} we obtain
\begin{equation}\label{eq3.2}
\omega_{\boldsymbol\alpha}^{\theta}\bigl(f,\frac{1}{2^{\nu_{1}}},...,\frac{1}{2^{\nu_{m}}}\bigr)_{p, \tau} \gg \prod_{j=1}^{m}2^{-\nu_{j}\alpha_{j}}\Biggl\|\Biggl(\sum\limits_{s_{1}=1}^{\nu_{1}}...\sum\limits_{s_{m}=1}^{\nu_{m}}\prod_{j=1}^{m}2^{s_{j}\alpha_{j}2}|\delta_{\mathbf{s}}(f)|^{2}\Biggr)^{\frac{1}{2}}\Biggr\|_{p, \tau}.
\end{equation}
Now, using inequality \eqref{eq3.2}, we have
\begin{equation}\label{eq3.3}
\Omega_{\boldsymbol{\alpha}}(f)_{p, \tau, \theta}\gg J_{p, \tau, \theta, \boldsymbol{\alpha}}(f)
\end{equation}
for the function $f\in \text{Lip}_{p, \tau, \theta}^{(\boldsymbol{\alpha}, -\mathbf{b})}(\mathbb{T}^{m})$, $1<p<\infty$, $1<\tau<\infty$.
Since  
\begin{equation}\label{eq3.4}
\sum\limits_{\nu=2^{l}}^{2^{l+1}}(\nu+)^{ -b}\asymp 2^{l(1-b)}, \, \, b\in \mathbb{R}, 
\end{equation}
then  
\begin{equation}\label{eq3.5} 
   \begin{gathered}
   \sum\limits_{\nu_{1} =2^{l_{1}}}^{2^{l_{1}+1}-1}...\sum\limits_{\nu_{m} =2^{l_{m}}}^{2^{l_{m}+1}-1}
\prod_{j=1}^{m}(\nu_{j} + 1)^{-\theta b_{j}}\Biggl\|\Biggl(\sum\limits_{s_{1}=1}^{\nu_{1}}...\sum\limits_{s_{m}=1}^{\nu_{m}}\prod_{j=1}^{m}2^{s_{j}\alpha_{j}2}|\delta_{\mathbf{s}}(f)|^{2}\Biggr)^{\frac{1}{2}}\Biggr\|_{p, \tau}^{\theta}
\\
\gg \Biggl\|\Biggl(\sum\limits_{s_{1}=1}^{2^{l_{1}}}...\sum\limits_{s_{m}=1}^{2^{l_{m}}}\prod_{j=1}^{m}2^{s_{j}\alpha_{j}2}|\delta_{\mathbf{s}}(f)|^{2}\Biggr)^{\frac{1}{2}}\Biggr\|_{p, \tau}^{\theta}\sum\limits_{\nu_{1} =2^{l_{1}}}^{2^{l_{1}+1}-1}...\sum\limits_{\nu_{m} =2^{l_{m}}}^{2^{l_{m}+1}-1}\prod_{j=1}^{m}(\nu_{j} + 1)^{-\theta b_{j}}
\\
\gg \prod_{j=1}^{m}2^{l_{j}(1-\theta b_{j})}\Biggl\|\Biggl(\sum\limits_{s_{1}=2^{l_{1}-1}}^{2^{l_{1}}-1}...\sum\limits_{s_{m}=2^{l_{m}-1}}^{2^{l_{m}}-1}\prod_{j=1}^{m}2^{s_{j}\alpha_{j}2}|\delta_{\mathbf{s}}(f)|^{2}\Biggr)^{\frac{1}{2}}\Biggr\|_{p, \tau}^{\theta}.
\end{gathered}
\end{equation}
Now, using Lemma 2.7 from inequalities \eqref{eq3.3} and \eqref{eq3.5}, we conclude that
\begin{equation}\label{eq3.6}
\|f\|_{p, \tau} + (J_{p, \tau, \theta, \boldsymbol{\alpha}}(f))^{1/\theta}
 \ll \|f\|_{\text{Lip}_{p, \tau, \theta}^{(\boldsymbol{\alpha}, -\mathbf{b})}}
\end{equation}
for the function $f\in \text{Lip}_{p, \tau, \theta}^{(\boldsymbol{\alpha}, -\mathbf{b})}(\mathbb{T}^{m})$, $1<p<\infty$, $1<\tau<\infty$.

Let us prove the opposite inequality to \eqref{eq3.6}.
Applying Theorem 2.1 and the inequality $(a+b)^{\theta}\leqslant C(\theta)(a^{\theta}+b^{\theta})$, $a, b\geqslant 0$, $0<\theta <\infty$ we have
\begin{multline}\label{eq3.7}
\Omega_{\boldsymbol{\alpha}}(f)_{p, \tau, \theta}\ll \sum\limits_{\nu_{m} =1}^{\infty}...\sum\limits_{\nu_{1} =1}^{\infty} \prod_{j=1}^{m}2^{\nu_{j}\alpha_{j}\theta}(\nu_{j} + 1)^{-\theta b_{j}}\Bigl(\prod_{j=1}^{m}2^{-\nu_{j}\alpha_{j}}\|S_{2^{\nu_{1}}-1,...,2^{\nu_{m}}-1}^{(\boldsymbol{\alpha})}(f)\|_{p, \tau} \Bigr)^{\theta}
\\
+\sum\limits_{\nu_{m} =1}^{\infty}...\sum\limits_{\nu_{1} =1}^{\infty} \prod_{j=1}^{m}2^{\nu_{j}\alpha_{j}\theta}(\nu_{j} + 1)^{-\theta b_{j}}  \sum_{e\subset e_{m}, e\neq \emptyset} \Bigl( \prod_{j\in e}2^{-\nu_{j}\alpha_{j}} \|S_{2^{\boldsymbol{\nu}^{e}}-\mathbf{1}^{e}, \infty}^{(\boldsymbol{\alpha}^{e})}(f-S_{\infty, 2^{\boldsymbol{\nu}^{\hat{e}}}-\mathbf{1}^{e}}(f))\|_{p, \tau} \Bigr)^{\theta}
\\
+\sum\limits_{\nu_{m} =1}^{\infty}...\sum\limits_{\nu_{1} =1}^{\infty} \prod_{j=1}^{m}2^{\nu_{j}\alpha_{j}\theta}(\nu_{j} + 1)^{-\theta b_{j}}\|f - U_{2^{\nu_{1}}-1,...,2^{\nu_{m}}-1} (f)\|_{p, \tau}.
\end{multline}
Let us estimate each term in \eqref{eq3.7}. By the definition of a rectangular partial sum of a Fourier series and by the Littlewood--Paley theorem in Lorentz space in \cite[Theorem 1.1]{2} we have
\begin{multline}\label{eq3.8}
   J_{1}:=\sum\limits_{\nu_{m} =1}^{\infty}...\sum\limits_{\nu_{1} =1}^{\infty} \prod_{j=1}^{m}(\nu_{j} + 1)^{-\theta b_{j}}\|S_{2^{\nu_{1}}-1,...,2^{\nu_{m}}-1}^{(\boldsymbol{\alpha})}(f)\|_{p, \tau}^{\theta}
\\
\ll \sum\limits_{\nu_{m} =1}^{\infty}...\sum\limits_{\nu_{1} =1}^{\infty} \prod_{j=1}^{m}(\nu_{j} + 1)^{-\theta b_{j}} \biggl\|\biggl(\sum\limits_{s_{1}=1}^{2^{l_{1}}-1}...\sum\limits_{s_{m}=1}^{2^{l_{m}}-1}\prod_{j=1}^{m}2^{s_{j}\alpha_{j}2}|\delta_{\mathbf{s}}(f)|^{2}\biggr)^{\frac{1}{2}}  \biggr\|_{p, \tau}^{\theta} 
\\
\ll \sum\limits_{l_{m} =1}^{\infty}...\sum\limits_{l_{1} =1}^{\infty}\biggl\|\biggl(\sum\limits_{s_{1}=1}^{2^{l_{1}}-1}...\sum\limits_{s_{m}=1}^{2^{l_{m}}-1}\prod_{j=1}^{m}2^{s_{j}\alpha_{j}2}|\delta_{\mathbf{s}}(f)|^{2}\biggr)^{\frac{1}{2}}  \biggr\|_{p, \tau}^{\theta}\sum\limits_{\nu_{1} =2^{l_{1}-1}}^{2^{l_{1}}-1}...\sum\limits_{\nu_{m} =2^{l_{m}-1}}^{2^{l_{m}}-1}\prod_{j=1}^{m}(\nu_{j} + 1)^{-\theta b_{j}}.
\end{multline}
Now, taking into account the relation \eqref{eq3.4} from \eqref{eq3.8} we obtain
\begin{equation}\label{eq3.9}
J_{1}\ll J_{p, \tau, \theta, \boldsymbol{\alpha}}(f).
\end{equation}
By virtue of Jensen's inequality \cite[Chapter 3, Section 3.3]{27} and the triangle inequality, we have
\begin{equation}\label{eq3.10} 
   \begin{gathered}
   \biggl\|\biggl(\sum\limits_{s_{1}=1}^{2^{l_{1}}-1}...\sum\limits_{s_{m}=1}^{2^{l_{m}}-1}\prod_{j=1}^{m}2^{s_{j}\alpha_{j}2}|\delta_{\mathbf{s}}(f)|^{2}\biggr)^{\frac{1}{2}}  \biggr\|_{p, \tau}
\\
\ll \biggl\|\sum\limits_{\nu_{1}=1}^{l_{1}}...\sum\limits_{\nu_{m}=1}^{l_{m}}\biggl(\sum\limits_{s_{1}=2^{\nu_{1}-1}}^{2^{\nu_{1}}-1}...\sum\limits_{s_{m}=2^{\nu_{m}-1}}^{2^{\nu_{m}}-1}\prod_{j=1}^{m}2^{s_{j}\alpha_{j}2}|\delta_{\mathbf{s}}(f)|^{2}\biggr)^{\frac{1}{2}}  \biggr\|_{p, \tau}
\\
\ll \sum\limits_{\nu_{1}=1}^{l_{1}}...\sum\limits_{\nu_{m}=1}^{l_{m}}\biggl\|\biggl(\sum\limits_{s_{1}=2^{\nu_{1}-1}}^{2^{\nu_{1}}-1}...\sum\limits_{s_{m}=2^{\nu_{m}-1}}^{2^{\nu_{m}}-1}\prod_{j=1}^{m}2^{s_{j}\alpha_{j}2}|\delta_{\mathbf{s}}(f)|^{2}\biggr)^{\frac{1}{2}}  \biggr\|_{p, \tau}.
\end{gathered}
\end{equation}
From inequalities \eqref{eq3.9} and \eqref{eq3.10}, we obtain
\begin{equation}\label{eq3.11}
J_{1}\ll \sum\limits_{l_{m} =1}^{\infty}...\sum\limits_{l_{1} =1}^{\infty}\prod_{j=1}^{m}2^{l_{j}(1-\theta b_{j})} \biggl(\sum\limits_{\nu_{1}=1}^{l_{1}}...\sum\limits_{\nu_{m}=1}^{l_{m}}\biggl\|\biggl(\sum\limits_{s_{1}=2^{\nu_{1}-1}}^{2^{\nu_{1}}-1}...\sum\limits_{s_{m}=2^{\nu_{m}-1}}^{2^{\nu_{m}}-1}\prod_{j=1}^{m}2^{s_{j}\alpha_{j}2}|\delta_{\mathbf{s}}(f)|^{2}\biggr)^{\frac{1}{2}}  \biggr\|_{p, \tau} \biggr)^{\theta}.
\end{equation}
Since $1-\theta b_{j}< 0$ for $j=1,\ldots ,m$, then
\begin{equation}\label{eq3.12}
\sum\limits_{\nu_{j} =l_{j}}^{\infty}2^{\nu_{j}(1-\theta b_{j})} \ll 2^{l_{j}(1-\theta b_{j})}. 
\end{equation}

Therefore, using Hardy's generalized inequality (see \cite[Lemma 3.2]{32}, \cite[Lemma B.2]{15}) from \eqref{eq3.11}, we obtain 
 \begin{equation}\label{eq3.13}
J_{1}\ll  \sum\limits_{l_{m} =1}^{\infty}...\sum\limits_{l_{1} =1}^{\infty}\prod_{j=1}^{m}2^{l_{j}(1-\theta b_{j})}\biggl\|\biggl(\sum\limits_{s_{1}=2^{\nu_{1}-1}}^{2^{\nu_{1}}-1}...\sum\limits_{s_{m}=2^{\nu_{m}-1}}^{2^{\nu_{m}}-1}\prod_{j=1}^{m}2^{s_{j}\alpha_{j}2}|\delta_{\mathbf{s}}(f)|^{2}\biggr)^{\frac{1}{2}}  \biggr\|_{p, \tau}^{\theta}
\end{equation}
for  $1<p<\infty$, $1<\tau<\infty$, $0< \theta< \infty$, $b_{j}> \frac{1}{\theta}$ for $j=1,\ldots ,m$.
Now we will evaluate 
 \begin{equation}\label{eq3.14}
J_{2}(f):=\sum\limits_{\nu_{m} =1}^{\infty}...\sum\limits_{\nu_{1} =1}^{\infty} \prod_{j\in \hat{e}}2^{\nu_{j}\alpha_{j}\theta}\prod_{j=1}^{m}(\nu_{j} + 1)^{-\theta b_{j}}\|S_{2^{\boldsymbol{\nu}^{e}}-\mathbf{1}^{e}, \infty}^{(\boldsymbol{\alpha}^{e})}(f-S_{\infty, 2^{\boldsymbol{\nu}^{\hat{e}}}-\mathbf{1}^{\hat{e}}}(f))\|_{p, \tau}^{\theta}. 
\end{equation}
By Theorem 2.1 in \cite{1}, we have 
 \begin{equation}\label{eq3.15}
\|S_{2^{\boldsymbol{\nu}^{e}}-\mathbf{1}^{e}, \infty}^{(\boldsymbol{\alpha}^{e})}(f-S_{\infty, 2^{\boldsymbol{\nu}^{\hat{e}}}-\mathbf{1}^{e}}(f))\|_{p, \tau} = \Bigl\|\sum\limits_{\mathbf{s}\in G_{\boldsymbol{\nu}}(e)} \delta_{\mathbf{s}}^{(\boldsymbol{\alpha}^{e})}(f)\Bigr\|_{p, \tau}\asymp \Bigl\|\Bigl(\sum\limits_{\mathbf{s}\in G_{\boldsymbol{\nu}}(e)} \prod_{j\in e}2^{s_{j}\alpha_{j}2}|\delta_{\mathbf{s}}(f)|^{2}\Bigr)^{\frac{1}{2}}\Bigr\|_{p, \tau} 
\end{equation}
for  $1<p<\infty$, $1<\tau<\infty$.
From inequalities \eqref{eq3.14} and \eqref{eq3.15}, it follows that 
 \begin{equation}\label{eq3.16}
J_{2}(f) \ll \sum\limits_{\nu_{m} =1}^{\infty}...\sum\limits_{\nu_{1} =1}^{\infty} \prod_{j\in \hat{e}}2^{\nu_{j}\alpha_{j}\theta}\prod_{j=1}^{m}(\nu_{j} + 1)^{-\theta b_{j}}\Bigl\|\Bigl(\sum\limits_{\mathbf{s}\in G_{\boldsymbol{\nu}}(e)} \prod_{j\in e}2^{s_{j}\alpha_{j}2}|\delta_{\mathbf{s}}(f)|^{2}\Bigr)^{\frac{1}{2}}\Bigr\|_{p, \tau}^{\theta}
\end{equation}
for 
 $1<p<\infty$, $1<\tau<\infty$.
Let the number $\varepsilon>0$. Since $\alpha_{j}>0$, the sequence $\{2^{\nu_{j}\alpha_{j}\theta}(\nu_{j} + 1)^{-\theta b_{j}+(\frac{1}{\theta}-\varepsilon)\theta}\}$ increases. Therefore
\begin{equation}\label{eq3.17} 
   \begin{gathered}
   \prod_{j\in \hat{e}}2^{\nu_{j}\alpha_{j}\theta}(\nu_{j} + 1)^{-\theta b_{j}} \Bigl\|\Bigl(\sum\limits_{\mathbf{s}\in G_{\boldsymbol{\nu}}(e)} \prod_{j\in e}2^{s_{j}\alpha_{j}2}|\delta_{\mathbf{s}}(f)|^{2}\Bigr)^{\frac{1}{2}}\Bigr\|_{p, \tau}^{\theta}
\\
\ll \prod_{j\in \hat{e}}(\nu_{j} + 1)^{\varepsilon\theta-1}\Bigl\|\Bigl(\sum\limits_{\mathbf{s}\in G_{\boldsymbol{\nu}}(e)} \prod_{j=1}^{m}2^{s_{j}\alpha_{j}2}\prod_{j\in \hat{e}}(\nu_{j} + 1)^{2(-b_{j}+\frac{1}{\theta}-\varepsilon)}|\delta_{\mathbf{s}}(f)|^{2}\Bigr)^{\frac{1}{2}}\Bigr\|_{p, \tau}^{\theta}.
\end{gathered}
\end{equation}
Now, from inequalities \eqref{eq3.16} and \eqref{eq3.17}, it follows that
\begin{equation}\label{eq3.18} 
   \begin{gathered}
   J_{2}(f) \ll \sum\limits_{\nu_{m} =1}^{\infty}...\sum\limits_{\nu_{1} =1}^{\infty} \prod_{j\in \hat{e}}(\nu_{j} + 1)^{\varepsilon\theta-1}\prod_{j\in e}(\nu_{j} + 1)^{-\theta b_{j}}
\\
\times
\Bigl\|\Bigl(\sum\limits_{\mathbf{s}\in G_{\boldsymbol{\nu}}(e)} \prod_{j=1}^{m}2^{s_{j}\alpha_{j}2}\prod_{j\in \hat{e}}(\nu_{j} + 1)^{2(-b_{j}+\frac{1}{\theta}-\varepsilon)}|\delta_{\mathbf{s}}(f)|^{2}\Bigr)^{\frac{1}{2}}\Bigr\|_{p, \tau}^{\theta}
\end{gathered}
\end{equation}
for 
 $1<p<\infty$, $1<\tau<\infty$.
Further, taking into account the relation \eqref{eq3.4} from \eqref{eq3.18} we will have
 \begin{equation*}  
   \begin{gathered}
   J_{2}(f) \ll \sum\limits_{l_{m} =0}^{\infty}...\sum\limits_{l_{1} =0}^{\infty}\sum\limits_{\nu_{1}=2^{l_{1}}}^{2^{l_{1}+1}-1}...\sum\limits_{s_{m}=2^{l_{m}}}^{2^{l_{m}+1}-1}\prod_{j\in \hat{e}}(\nu_{j} + 1)^{\varepsilon\theta-1}\prod_{j\in e}(\nu_{j} + 1)^{-\theta b_{j}}
\\
\times
\Bigl\|\Bigl(\sum\limits_{\mathbf{s}\in G_{\boldsymbol{\nu}}(e)} \prod_{j=1}^{m}2^{s_{j}\alpha_{j}2}\prod_{j\in \hat{e}}(s_{j} + 1)^{2(-b_{j}+\frac{1}{\theta}-\varepsilon)}|\delta_{\mathbf{s}}(f)|^{2}\Bigr)^{\frac{1}{2}}\Bigr\|_{p, \tau}^{\theta}
\end{gathered}
\end{equation*}
\begin{equation}\label{eq3.19} 
   \begin{gathered}
\ll \sum\limits_{\mathbf{l}\in \mathbb{Z}_{+}^{m}}\prod_{j\in e}2^{l_{j}(1-\theta b_{j})}\prod_{j\in \hat{e}}2^{l_{j}\varepsilon\theta}\Bigl\|\Bigl(\sum\limits_{\mathbf{s}\in G_{2^{\mathbf{l}}}(e)} \prod_{j=1}^{m}2^{s_{j}\alpha_{j}2}\prod_{j\in \hat{e}}(s_{j} + 1)^{2(-b_{j}+\frac{1}{\theta}-\varepsilon)}|\delta_{\mathbf{s}}(f)|^{2}\Bigr)^{\frac{1}{2}}\Bigr\|_{p, \tau}^{\theta}.
\end{gathered}
\end{equation}
According to Jensen's inequality in 
\cite[Chapter 3, Section 3.3]{27} and the triangle inequality, we have
\begin{equation*}  
   \begin{gathered}
   \Bigl\|\Bigl(\sum\limits_{\mathbf{s}\in G_{2^{\mathbf{l}}}(e)} \prod_{j=1}^{m}2^{s_{j}\alpha_{j}2}\prod_{j\in \hat{e}}(s_{j} + 1)^{2(-b_{j}+\frac{1}{\theta}-\varepsilon)}|\delta_{\mathbf{s}}(f)|^{2}\Bigr)^{\frac{1}{2}}\Bigr\|_{p, \tau}
\\
 = \Bigl\|\Bigl(\sum\limits_{\boldsymbol{\nu}\in G_{\mathbf{l}}(e)}\sum\limits_{s_{1}=2^{\nu_{1}}}^{2^{\nu_{1}+1}-1}...\sum\limits_{s_{m}=2^{\nu_{m}}}^{2^{\nu_{m}+1}-1} \prod_{j=1}^{m}2^{s_{j}\alpha_{j}2}\prod_{j\in \hat{e}}(s_{j} + 1)^{2(-b_{j}+\frac{1}{\theta}-\varepsilon)}|\delta_{\mathbf{s}}(f)|^{2}\Bigr)^{\frac{1}{2}}\Bigr\|_{p, \tau}
 \\
 \ll \Bigl\|\sum\limits_{\boldsymbol{\nu}\in G_{\mathbf{l}}(e)}\Bigl(\sum\limits_{s_{1}=2^{\nu_{1}}}^{2^{\nu_{1}+1}-1}...\sum\limits_{s_{m}=2^{\nu_{m}}}^{2^{\nu_{m}+1}-1} \prod_{j=1}^{m}2^{s_{j}\alpha_{j}2}\prod_{j\in \hat{e}}(s_{j} + 1)^{2(-b_{j}+\frac{1}{\theta}-\varepsilon)}|\delta_{\mathbf{s}}(f)|^{2}\Bigr)^{\frac{1}{2}}\Bigr\|_{p, \tau} 
\end{gathered}
\end{equation*}
\begin{equation}\label{eq3.20}
 \ll \sum\limits_{\boldsymbol{\nu}\in G_{\mathbf{l}}(e)} \Bigl\|\Bigl(\sum\limits_{s_{1}=2^{\nu_{1}}}^{2^{\nu_{1}+1}-1}...\sum\limits_{s_{m}=2^{\nu_{m}}}^{2^{\nu_{m}+1}-1} \prod_{j=1}^{m}2^{s_{j}\alpha_{j}2}\prod_{j\in \hat{e}}(s_{j} + 1)^{2(-b_{j}+\frac{1}{\theta}-\varepsilon)}|\delta_{\mathbf{s}}(f)|^{2}\Bigr)^{\frac{1}{2}}\Bigr\|_{p, \tau}.
\end{equation}
Now from the inequalities \eqref{eq3.19} and \eqref{eq3.20} we obtain
\begin{equation}\label{eq3.21} 
   \begin{gathered}
   J_{2}(f) \ll \sum\limits_{l_{m} =0}^{\infty}...\sum\limits_{l_{1} =0}^{\infty}\prod_{j\in e}2^{l_{j}(1-\theta b_{j})}\prod_{j\in \hat{e}}2^{l_{j}\varepsilon\theta}
\\
\times
\Bigl(\sum\limits_{\boldsymbol{\nu}\in G_{\mathbf{l}}(e)} \Bigl\|\Bigl(\sum\limits_{s_{1}=2^{\nu_{1}}}^{2^{\nu_{1}+1}-1}...\sum\limits_{s_{m}=2^{\nu_{m}}}^{2^{\nu_{m}+1}-1} \prod_{j=1}^{m}2^{s_{j}\alpha_{j}2}\prod_{j\in \hat{e}}(s_{j} + 1)^{2(-b_{j}+\frac{1}{\theta}-\varepsilon)}|\delta_{\mathbf{s}}(f)|^{2}\Bigr)^{\frac{1}{2}}\Bigr\|_{p, \tau}\Bigr)^{\theta}
\end{gathered}
\end{equation}
for   $1<p<\infty$, $1<\tau<\infty$, $0<\theta<\infty$. 
Since $\varepsilon > 0$, then
 \begin{equation*}
\sum\limits_{\nu_{j}=0}^{l_{j}}2^{\nu_{j}\varepsilon\theta} \ll 2^{l_{j}\varepsilon\theta}, \, \, j=1,\ldots,m. 
\end{equation*}
Therefore, taking into account this estimate and inequality \eqref{eq3.12} to the right-hand side \eqref{eq3.21} applying Hardy's inequality (see \cite[Lemma 3.2]{32}, \cite[Lemma B.2]{15}) we will have
\begin{equation}\label{eq3.22} 
   \begin{gathered}
   J_{2}(f) \ll \sum\limits_{l_{m} =0}^{\infty}...\sum\limits_{l_{1} =0}^{\infty}\prod_{j\in e}2^{l_{j}(1-\theta b_{j})}\prod_{j\in \hat{e}}2^{l_{j}\varepsilon\theta} 
\\
\times
\Bigl\|\Bigl(\sum\limits_{s_{1}=2^{l_{1}}}^{2^{l_{1}+1}-1}...\sum\limits_{s_{m}=2^{l_{m}}}^{2^{l_{m}+1}-1} \prod_{j=1}^{m}2^{s_{j}\alpha_{j}2}\prod_{j\in \hat{e}}(s_{j} + 1)^{2(-b_{j}+\frac{1}{\theta}-\varepsilon)}|\delta_{\mathbf{s}}(f)|^{2}\Bigr)^{\frac{1}{2}}\Bigr\|_{p, \tau}^{\theta} \ll J_{p, \tau, \theta, \boldsymbol\alpha}(f)
\end{gathered}
\end{equation}
for   $1<p<\infty$, $1<\tau<\infty$, $0<\theta<\infty$.
Now, we will evaluate
\begin{equation*} 
   \begin{gathered}
   J_{3}(f):=\sum\limits_{\nu_{m} =1}^{\infty}...\sum\limits_{\nu_{1} =1}^{\infty} \prod_{j=1}^{m}2^{\nu_{j}\alpha_{j}\theta}(\nu_{j} + 1)^{-\theta b_{j}}\|f - U_{2^{\nu_{1}}-1,...,2^{\nu_{m}}-1} (f)\|_{p, \tau} 
\\
= \sum\limits_{\nu_{m} =1}^{\infty}...\sum\limits_{\nu_{1} =1}^{\infty} \prod_{j=1}^{m}2^{\nu_{j}\alpha_{j}\theta}(\nu_{j} + 1)^{-\theta b_{j}}\Bigl\|\sum\limits_{s_{m}=\nu_{m} +1}^{\infty}...\sum\limits_{s_{1}=\nu_{1}+1}^{\infty} \delta_{\mathbf{s}}(f) \Bigr\|_{p, \tau}^{\theta}.
\end{gathered}
\end{equation*}
According to the Littlewood--Paley theorem in the Lorentz space in \cite{1} and taking into account that the sequence $\{2^{\nu_{j}\alpha_{j}\theta}(\nu_{j} + 1)^{-\theta b_{j}\theta}\}$ increases, we obtain
\begin{equation}\label{eq3.23} 
   \begin{gathered}
   J_{3}(f) \ll \sum\limits_{\nu_{m} =1}^{\infty}...\sum\limits_{\nu_{1} =1}^{\infty} \prod_{j=1}^{m}2^{\nu_{j}\alpha_{j}\theta}(\nu_{j} + 1)^{-\theta b_{j}}\Bigl\|\Bigl(\sum\limits_{s_{m}=\nu_{m} +1}^{\infty}...\sum\limits_{s_{1}=\nu_{1}+1}^{\infty} |\delta_{\mathbf{s}}(f)|^{2}\Bigr)^{\frac{1}{2}} \Bigr\|_{p, \tau}^{\theta}
\\
\ll  \sum\limits_{\nu_{m} =1}^{\infty}...\sum\limits_{\nu_{1} =1}^{\infty} \Bigl\|\Bigl(\sum\limits_{s_{m}=\nu_{m} +1}^{\infty}...\sum\limits_{s_{1}=\nu_{1}+1}^{\infty}\prod_{j=1}^{m}2^{s_{j}\alpha_{j}2}(s_{j} + 1)^{-2b_{j}} |\delta_{\mathbf{s}}(f)|^{2}\Bigr)^{\frac{1}{2}} \Bigr\|_{p, \tau}^{\theta}
\\
\ll \sum\limits_{l_{1} =0}^{\infty}...\sum\limits_{l_{m} =0}^{\infty}\sum\limits_{\nu_{1}=2^{l_{1}}}^{2^{l_{1}+1}-1}...\sum\limits_{\nu_{m}=2^{l_{m}}}^{2^{l_{m}+1}-1}\Bigl\|\Bigl(\sum\limits_{s_{m}=\nu_{m}}^{\infty}...\sum\limits_{s_{1}=\nu_{1}}^{\infty}\prod_{j=1}^{m}2^{s_{j}\alpha_{j}2}(s_{j} + 1)^{-2b_{j}} |\delta_{\mathbf{s}}(f)|^{2}\Bigr)^{\frac{1}{2}} \Bigr\|_{p, \tau}^{\theta}
\\
\ll \sum\limits_{l_{1} =0}^{\infty}...\sum\limits_{l_{m} =0}^{\infty}\prod_{j=1}^{m}2^{l_{j}}\Bigl\|\Bigl(\sum\limits_{s_{m}=2^{l_{m}}}^{\infty}...\sum\limits_{s_{1}=2^{l_{1}}}^{\infty}\prod_{j=1}^{m}2^{s_{j}\alpha_{j}2}(s_{j} + 1)^{-2b_{j}} |\delta_{\mathbf{s}}(f)|^{2}\Bigr)^{\frac{1}{2}} \Bigr\|_{p, \tau}^{\theta}.
\end{gathered}
\end{equation}
Further, using Jensen's inequality, then applying the triangle inequality and Hardy's inequality (see  \cite[Lemma B.2]{15},  \cite[Lemma 3.2]{32}) from \eqref{eq3.23} we obtain
\begin{equation*} 
   \begin{gathered}
   J_{3}(f) \ll \sum\limits_{l_{1} =0}^{\infty}...\sum\limits_{l_{m} =0}^{\infty}\prod_{j=1}^{m}2^{l_{j}}\Bigl\|\Bigl(\sum\limits_{\nu_{1}=l_{1}}^{\infty}...\sum\limits_{\nu_{m}=l_{m}}^{\infty}  \sum\limits_{s_{1}=2^{\nu_{1}}}^{2^{\nu_{1}+1}-1}...\sum\limits_{s_{m}=2^{\nu_{m}}}^{2^{\nu_{m}+1}-1}\prod_{j=1}^{m}2^{s_{j}\alpha_{j}2}(s_{j} + 1)^{-2b_{j}} |\delta_{\mathbf{s}}(f)|^{2}\Bigr)^{\frac{1}{2}}\Bigr\|_{p, \tau}^{\theta}
\\
\ll \sum\limits_{l_{1} =0}^{\infty}...\sum\limits_{l_{m} =0}^{\infty}\prod_{j=1}^{m}2^{l_{j}}\Bigl\|\sum\limits_{\nu_{1}=l_{1}}^{\infty}...\sum\limits_{\nu_{m}=l_{m}}^{\infty} \Bigl( \sum\limits_{s_{1}=2^{\nu_{1}}}^{2^{\nu_{1}+1}-1}...\sum\limits_{s_{m}=2^{\nu_{m}}}^{2^{\nu_{m}+1}-1}\prod_{j=1}^{m}2^{s_{j}\alpha_{j}2}(s_{j} + 1)^{-2b_{j}} |\delta_{\mathbf{s}}(f)|^{2}\Bigr)^{\frac{1}{2}}\Bigr\|_{p, \tau}^{\theta}
\\
\ll \sum\limits_{l_{1} =0}^{\infty}...\sum\limits_{l_{m} =0}^{\infty}\prod_{j=1}^{m}2^{l_{j}}\Bigl( 
\sum\limits_{\nu_{1}=l_{1}}^{\infty}...\sum\limits_{\nu_{m}=l_{m}}^{\infty}\Bigl\|\Bigl( \sum\limits_{s_{1}=2^{\nu_{1}}}^{2^{\nu_{1}+1}-1}...\sum\limits_{s_{m}=2^{\nu_{m}}}^{2^{\nu_{m}+1}-1}\prod_{j=1}^{m}2^{s_{j}\alpha_{j}2}(s_{j} + 1)^{-2b_{j}} |\delta_{\mathbf{s}}(f)|^{2}\Bigr)^{\frac{1}{2}}\Bigr\|_{p, \tau}\Bigr)^{\theta}
\\
\ll \sum\limits_{l_{1} =0}^{\infty}...\sum\limits_{l_{m} =0}^{\infty}\prod_{j=1}^{m}2^{l_{j}}\Bigl\|\Bigl( \sum\limits_{s_{1}=2^{l_{1}}}^{2^{l_{1}+1}-1}...\sum\limits_{s_{m}=2^{l_{m}}}^{2^{l_{m}+1}-1}\prod_{j=1}^{m}2^{s_{j}\alpha_{j}2}(s_{j} + 1)^{-2b_{j}} |\delta_{\mathbf{s}}(f)|^{2}\Bigr)^{\frac{1}{2}}\Bigr\|_{p, \tau}^{\theta}
\ll  J_{p, \tau, \theta, \boldsymbol{\alpha}}(f).
\end{gathered}
\end{equation*}
Thus we have proved that  
 \begin{equation}\label{eq3.24}
J_{3}(f) \ll \sum\limits_{l_{1} =0}^{\infty}...\sum\limits_{l_{m} =0}^{\infty}\prod_{j=1}^{m}2^{l_{j}(1-b_{j}\theta)} \Bigl\|\Bigl( \sum\limits_{s_{1}=2^{l_{1}}}^{2^{l_{1}+1}-1}...\sum\limits_{s_{m}=2^{l_{m}}}^{2^{l_{m}+1}-1}\prod_{j=1}^{m}2^{s_{j}\alpha_{j}2} |\delta_{\mathbf{s}}(f)|^{2}\Bigr)^{\frac{1}{2}}\Bigr\|_{p, \tau}^{\theta}.
\end{equation}
Now, from inequalities \eqref{eq3.6}, \eqref{eq3.13}, \eqref{eq3.22}, and \eqref{eq3.24}, it follows that
 \begin{equation}\label{eq3.25}
\Omega_{\boldsymbol{\alpha}}(f)_{p, \tau, \theta}\ll J_{p, \tau, \theta, \boldsymbol{\alpha}}(f)
\end{equation}
for  $1<p<\infty$, $1<\tau<\infty$, $0<\theta<\infty$.
From Lemma 2.7 and inequality \eqref{eq3.25}, we obtain
 \begin{equation*}
\|f\|_{\text{Lip}_{p, \tau, \theta}^{(\boldsymbol{\alpha}, -\mathbf{b})}}\ll \|f\|_{p, \tau} + \Bigl(\sum\limits_{l_{1} =0}^{\infty}...\sum\limits_{l_{m} =0}^{\infty}\prod_{j=1}^{m}2^{l_{j}(1-b_{j}\theta)} \Bigl\|\Bigl( \sum\limits_{s_{1}=2^{l_{1}}}^{2^{l_{1}+1}-1}...\sum\limits_{s_{m}=2^{l_{m}}}^{2^{l_{m}+1}-1}\prod_{j=1}^{m}2^{s_{j}\alpha_{j}2} |\delta_{\mathbf{s}}(f)|^{2}\Bigr)^{\frac{1}{2}}\Bigr\|_{p, \tau}^{\theta} \Bigr)^{\frac{1}{\theta}} 
\end{equation*}
for  $1<p<\infty$, $1<\tau<\infty$, $0<\theta<\infty$.
\hfill $\Box$

\setcounter{equation}{0}
\setcounter{lemma}{0}
\setcounter{theorem}{0}

\section{Proofs of Theorem 1.2 and Theorem 1.3}\label{sec4}

Let $\mathbf{e}=\sum_{\nu=1}^{m}\mathbf{e}_{\nu}=(1,...,1)$, where $\{\mathbf{e}_{\nu}\}_{\nu=1}^{m}$ is a standard basis in $\mathbb{R}^{m}$.

\textbf{Proof of Theorem 1.2.} 
Let $f\in S_{p, \tau, \theta}^{\boldsymbol{\alpha}, -\mathbf{b}+\frac{1}{\min\{2, \tau, \theta\}}\mathbf{e}}B(\mathbb{T}^{m})$.
By Theorem 3.2 \cite{2}, the following inequality holds
 \begin{equation*}
\omega_{\boldsymbol{\alpha}}\Bigl(f,\frac{\pi}{n_{1}},...,\frac{\pi}{n_{m}}\Bigr)_{p, \tau} \ll
 \prod_{j=1}^{m}n_{j}^{-\alpha_{j}}\Biggl(\sum\limits_{\nu_{1}=1}^{n_{1}+1} \ldots \sum\limits_{\nu_{m}=1}^{n_{m}+1} \prod_{j=1}^{m}\nu_{j}^{\beta\alpha_{j} - 1}Y_{\boldsymbol{\nu}}^{\beta}(f)_{p, \tau}\Biggr)^{1/\beta}, \,\, n_{j}\in \mathbb{N},
\end{equation*}
where  $\beta = \min\{2, \tau\}$.
Due to the decreasing order of each index of the sequence $\{Y_{\boldsymbol{\nu}}(f)_{p, \tau}\}$
and the property of the mixed modulus of smoothness, from this inequality we obtain
\begin{equation}\label{eq4.1}
\omega_{\boldsymbol{\alpha}}\Bigl(f,\frac{1}{2^{\nu_{1}}},...,\frac{1}{2^{\nu_{m}}}\Bigr)_{p, \tau} \ll
 \prod_{j=1}^{m}2^{-\nu_{j}\alpha_{j}}\Biggl(\sum\limits_{k_{1}=0}^{\nu_{1}} \ldots \sum\limits_{k_{m}=0}^{\nu_{m}} \prod_{j=1}^{m}2^{k_{j}\beta\alpha_{j}}Y_{2^{k_{1}},\ldots,2^{k_{m}}}^{\beta}(f)_{p, \tau}\Biggr)^{1/\beta}.
\end{equation}
Now, using Lemma 2.7 and inequality \eqref{eq4.1}, we obtain
\begin{equation}\label{eq4.2}
\|f\|_{\text{Lip}_{p, \tau, \theta}^{(\boldsymbol{\alpha}, -\mathbf{b})}}\ll \|f\|_{p, \tau} +\Biggl(\sum\limits_{\boldsymbol{\nu}\in \mathbb{Z}_{+}^{m}}  \prod_{j=1}^{m}(\nu_{j} + 1)^{-\theta b_{j}} \Biggl(\sum\limits_{k_{1}=0}^{\nu_{1}} \ldots \sum\limits_{k_{m}=0}^{\nu_{m}} \prod_{j=1}^{m}2^{k_{j}\beta\alpha_{j}}Y_{2^{k_{1}},\ldots,2^{k_{m}}}^{\beta}(f)_{p, \tau}\Biggr)^{\frac{\theta}{\beta}} \Biggr)^{\frac{1}{\theta}}.
\end{equation}
If $\frac{\theta}{\beta}> 1$, then according to Hardy's inequality in \cite[Theorem 346]{20} taking into account that $b_{j}\theta >1$ for $j=1,\ldots ,m$ we will have
 \begin{equation}\label{eq4.3}
   \begin{gathered} 
\sum\limits_{\nu_{m} =0}^{\infty}...\sum\limits_{\nu_{1} =0}^{\infty} \prod_{j=1}^{m}(\nu_{j} + 1)^{-\theta b_{j}} \Biggl(\sum\limits_{k_{1}=0}^{\nu_{1}} \ldots \sum\limits_{k_{m}=0}^{\nu_{m}} \prod_{j=1}^{m}2^{k_{j}\beta\alpha_{j}}Y_{2^{k_{1}},\ldots,2^{k_{m}}}^{\beta}(f)_{p, \tau}\Biggr)^{\frac{\theta}{\beta}}
\\
\ll \sum\limits_{\nu_{m} =0}^{\infty}...\sum\limits_{\nu_{1} =0}^{\infty} \prod_{j=1}^{m}2^{\nu_{j}\theta\alpha_{j}}(\nu_{j} + 1)^{\theta(\frac{1}{\beta} - b_{j})}Y_{2^{k_{1}},\ldots,2^{k_{m}}}^{\theta}(f)_{p, \tau}.
\end{gathered}
\end{equation}
Now, from inequalities \eqref{eq4.2} and \eqref{eq4.3}, it follows that
\begin{equation}\label{eq4.4}
\|f\|_{\text{Lip}_{p, \tau, \theta}^{(\boldsymbol{\alpha}, -\mathbf{b})}}\ll \|f\|_{p, \tau} +\Biggl(\sum\limits_{\nu_{m} =0}^{\infty}...\sum\limits_{\nu_{1} =0}^{\infty} \prod_{j=1}^{m}2^{\nu_{j}\theta\alpha_{j}}(\nu_{j} + 1)^{\theta(\frac{1}{\beta} - b_{j})}Y_{2^{k_{1}},\ldots,2^{k_{m}}}^{\theta}(f)_{p, \tau}\Biggr)^{\frac{1}{\theta}}
\end{equation}
in the case $\beta < \theta \leqslant \infty$. 

If $\frac{\theta}{\beta}\leqslant 1$, then applying Jensen's inequality (see \cite[Ch.~3, Sec.~3]{27}) and changing the order of summation, taking into account that  $b_{j}\theta >1$ for $j=1,\ldots ,m$ , we obtain
 \begin{equation}\label{eq4.5}
   \begin{gathered} 
\sum\limits_{\nu_{m} =0}^{\infty}...\sum\limits_{\nu_{1} =0}^{\infty} \prod_{j=1}^{m}(\nu_{j} + 1)^{-\theta b_{j}} \Biggl(\sum\limits_{k_{1}=0}^{\nu_{1}} \ldots \sum\limits_{k_{m}=0}^{\nu_{m}} \prod_{j=1}^{m}2^{k_{j}\beta\alpha_{j}}Y_{2^{k_{1}},\ldots,2^{k_{m}}}^{\beta}(f)_{p, \tau}\Biggr)^{\frac{\theta}{\beta}}
\\
\ll \sum\limits_{\nu_{m} =0}^{\infty}...\sum\limits_{\nu_{1} =0}^{\infty} \prod_{j=1}^{m}(\nu_{j} + 1)^{-\theta b_{j}} \sum\limits_{k_{1}=0}^{\nu_{1}} \ldots \sum\limits_{k_{m}=0}^{\nu_{m}} \prod_{j=1}^{m}2^{k_{j}\theta\alpha_{j}}Y_{2^{k_{1}},\ldots,2^{k_{m}}}^{\theta}(f)_{p, \tau}
\\
= \sum\limits_{k_{m} =0}^{\infty}...\sum\limits_{k_{1} =0}^{\infty}\prod_{j=1}^{m}2^{k_{j}\theta\alpha_{j}}Y_{2^{k_{1}},\ldots,2^{k_{m}}}^{\theta}(f)_{p, \tau}\sum\limits_{\nu_{m} =k_{m}}^{\infty}...\sum\limits_{\nu_{1} =k_{1}}^{\infty} \prod_{j=1}^{m}(\nu_{j} + 1)^{-\theta b_{j}}
\\
\ll \sum\limits_{k_{m} =0}^{\infty}...\sum\limits_{k_{1} =0}^{\infty}\prod_{j=1}^{m}2^{k_{j}\theta\alpha_{j}}(k_{j} + 1)^{1-\theta b_{j}}Y_{2^{k_{1}},\ldots,2^{k_{m}}}^{\theta}(f)_{p, \tau}.
\end{gathered}
\end{equation}
Now, from inequalities \eqref{eq4.2} and \eqref{eq4.5}, it follows that
 \begin{equation}\label{eq4.6}
\|f\|_{\text{Lip}_{p, \tau, \theta}^{(\boldsymbol{\alpha}, -\mathbf{b})}}\ll \|f\|_{p, \tau} +\Biggl(\sum\limits_{k_{m} =0}^{\infty}...\sum\limits_{k_{1} =0}^{\infty}\prod_{j=1}^{m}2^{k_{j}\beta\alpha_{j}}(k_{j} + 1)^{1-\theta b_{j}}Y_{2^{k_{1}},\ldots,2^{k_{m}}}^{\theta}(f)_{p, \tau}\Biggr)^{\frac{1}{\theta}}
\end{equation} 
in the case $0< \theta \leqslant \beta$.
Further, using Lemma 2.5 from the inequalities \eqref{eq4.4} and \eqref{eq4.6} we obtain
 \begin{equation*}
\|f\|_{\text{Lip}_{p, \tau, \theta}^{(\boldsymbol{\alpha}, -\mathbf{b})}}\ll \|f\|_{p, \tau} +\Biggl(\sum\limits_{k_{m} =0}^{\infty}...\sum\limits_{k_{1} =0}^{\infty}\prod_{j=1}^{m}2^{k_{j}\beta\alpha_{j}}(k_{j} + 1)^{(\frac{1}{\min\{2, \tau, \theta\}} - b_{j})\theta}\omega_{\mathbf{r}}^{\theta}(f, \frac{1}{2^{k_{1}}}, \ldots , \frac{1}{2^{k_{m}}})_{p, \tau}\Biggr)^{\frac{1}{\theta}}
\end{equation*} 
for  $r_{j}> \alpha_{j}$, $j=1,\ldots,m$, $\overline{r}=(r_{1},\ldots, r_{m})$.
Due to the property of the mixed smoothness modulus, it follows that
    \begin{equation*}
\|f\|_{\text{Lip}_{p, \tau, \theta}^{(\boldsymbol{\alpha}, -\mathbf{b})}}\ll \|f\|_{p, \tau} +\Biggl(\int_{0}^{1}...\int_{0}^{1}\omega_{\mathbf{r}}^{\theta}(f, \mathbf{t})_{p, \tau}\prod_{j=1}^{m}\frac{(1-\log t_{j})^{(\frac{1}{\min\{2, \tau, \theta\}}- b_{j})\theta}}{t_{j}^{\alpha_{j}\theta+1}}dt_{j}
\Biggr)^{\frac{1}{\theta}}
\end{equation*} 
for  $r_{j}> \alpha_{j}$, $j=1,\ldots,m$. Hence  
 \begin{equation}\label{eq4.7}
 S_{p, \tau, \theta}^{\boldsymbol{\alpha}, -\mathbf{b}+\frac{1}{\min\{2, \tau, \theta\}}\mathbf{e}}B(\mathbb{T}^{m}) \hookrightarrow \text{Lip}_{p, \tau, \theta}^{(\boldsymbol{\alpha}, -\mathbf{b})}(\mathbb{T}^{m}).
 \end{equation}
 Let $f \in \text{Lip}_{p, \tau, \theta}^{(\boldsymbol{\alpha}, -\mathbf{b})}(\mathbb{T}^{m})$. By Theorem 4.1 in \cite{1}, the following inequality holds
 \begin{equation}\label{eq4.8}
 \prod_{j=1}^{m}n_{j}^{-\alpha_{j}}\Biggl(\sum\limits_{\nu_{1}=1}^{n_{1}+1} \ldots \sum\limits_{\nu_{m}=1}^{n_{m}+1} \prod_{j=1}^{m}\nu_{j}^{\sigma\alpha_{j} - 1}Y_{\boldsymbol{\nu}}^{\beta}(f)_{p, \tau}\Biggr)^{1/\sigma} \ll \omega_{\boldsymbol{\alpha}}\Bigl(f,\frac{\pi}{n_{1}},...,\frac{\pi}{n_{m}}\Bigr)_{p, \tau}, \,\, n_{j}\in \mathbb{N},
\end{equation}
where  $\sigma = \max\{2, \tau\}$, $\alpha_{j}\in \mathbb{N}$. 
Due to the decay of each index of the sequence $\{Y_{\boldsymbol{\nu}}(f)_{p, \tau}\}$
and the property of the mixed smoothness modulus from  inequality \eqref{eq4.8}, we obtain
\begin{equation}\label{eq4.9}
  \prod_{j=1}^{m}2^{-\nu_{j}\alpha_{j}}\Biggl(\sum\limits_{k_{1}=0}^{\nu_{1}} \ldots \sum\limits_{k_{m}=0}^{\nu_{m}} \prod_{j=1}^{m}2^{k_{j}\sigma\alpha_{j}}Y_{2^{k_{1}},\ldots,2^{k_{m}}}^{\sigma}(f)_{p, \tau}\Biggr)^{1/\sigma} \ll \omega_{\boldsymbol\alpha}\Bigl(f,\frac{1}{2^{\nu_{1}}},...,\frac{1}{2^{\nu_{m}}}\Bigr)_{p, \tau}.
\end{equation}
From Lemma 2.7 and inequality \eqref{eq4.9}, it follows that 
 \begin{equation}\label{eq4.10}
  \begin{gathered}
\|f\|_{\text{Lip}_{p, \tau, \theta}^{(\boldsymbol{\alpha}, -\mathbf{b})}}\gg \|f\|_{p, \tau} 
\\
+\Biggl(\sum\limits_{\nu_{m} =0}^{\infty}...\sum\limits_{\nu_{1} =0}^{\infty} \prod_{j=1}^{m}(\nu_{j} + 1)^{-\theta b_{j}} \Biggl(\sum\limits_{k_{1}=0}^{\nu_{1}} \ldots \sum\limits_{k_{m}=0}^{\nu_{m}} \prod_{j=1}^{m}2^{k_{j}\sigma\alpha_{j}}Y_{2^{k_{1}},\ldots,2^{k_{m}}}^{\sigma}(f)_{p, \tau}\Biggr)^{\frac{\theta}{\sigma}} \Biggr)^{\frac{1}{\theta}}.
\end{gathered}
\end{equation}
If $\frac{\sigma}{\theta}\leqslant 1$, then using Jensen's inequality (see \cite[Ch.~3, Sec.~3]{27}), then changing the order of summation, taking into account that  $b_{j}\theta >1$ for $j=1,\ldots ,m$  we obtain
   \begin{equation}\label{eq4.11}
   \begin{gathered}
\sum\limits_{\nu_{m} =0}^{\infty}...\sum\limits_{\nu_{1} =0}^{\infty} \prod_{j=1}^{m}(\nu_{j} + 1)^{-\theta b_{j}} \Biggl(\sum\limits_{k_{1}=0}^{\nu_{1}} \ldots \sum\limits_{k_{m}=0}^{\nu_{m}} \prod_{j=1}^{m}2^{k_{j}\sigma\alpha_{j}}Y_{2^{k_{1}},\ldots,2^{k_{m}}}^{\sigma}(f)_{p, \tau}\Biggr)^{\frac{\theta}{\sigma}}
\\
\gg \sum\limits_{\nu_{m} =0}^{\infty}...\sum\limits_{\nu_{1} =0}^{\infty} \prod_{j=1}^{m}(\nu_{j} + 1)^{-\theta b_{j}} \sum\limits_{k_{1}=0}^{\nu_{1}} \ldots \sum\limits_{k_{m}=0}^{\nu_{m}} \prod_{j=1}^{m}2^{k_{j}\theta\alpha_{j}}Y_{2^{k_{1}},\ldots,2^{k_{m}}}^{\theta}(f)_{p, \tau}
\\
= \sum\limits_{k_{m} =0}^{\infty}...\sum\limits_{k_{1} =0}^{\infty}\prod_{j=1}^{m}2^{k_{j}\theta\alpha_{j}}Y_{2^{k_{1}},\ldots,2^{k_{m}}}^{\theta}(f)_{p, \tau}\sum\limits_{\nu_{m} =k_{m}}^{\infty}...\sum\limits_{\nu_{1} =k_{1}}^{\infty} \prod_{j=1}^{m}(\nu_{j} + 1)^{-\theta b_{j}}
\\
\gg
 \sum\limits_{k_{m} =0}^{\infty}...\sum\limits_{k_{1} =0}^{\infty}\prod_{j=1}^{m}2^{k_{j}\theta\alpha_{j}}(k_{j} + 1)^{1-\theta b_{j}}Y_{2^{k_{1}},\ldots,2^{k_{m}}}^{\theta}(f)_{p, \tau}.
\end{gathered}
\end{equation}
Now, from inequalities \eqref{eq4.10} and \eqref{eq4.11}, it follows that
 \begin{equation}\label{eq4.12}
 \|f\|_{\text{Lip}_{p, \tau, \theta}^{(\boldsymbol{\alpha}, -\mathbf{b})}}\gg \|f\|_{p, \tau}
 + \biggl(\sum\limits_{k_{m} =0}^{\infty}...\sum\limits_{k_{1} =0}^{\infty}\prod_{j=1}^{m}2^{k_{j}\theta\alpha_{j}}(k_{j} + 1)^{1-\theta b_{j}}Y_{2^{k_{1}},\ldots,2^{k_{m}}}^{\theta}(f)_{p, \tau}\biggr)^{1/\theta}
\end{equation}
in the case  $\sigma\leqslant\theta $. 
Furthermore, it is known that from the direct and inverse theorems of approximation by angle (see Lemma 2.4 and Lemma 2.5), the following relationship follows:
  \begin{equation}\label{eq4.13}
   \begin{gathered}
\sum\limits_{k_{m} =0}^{\infty}...\sum\limits_{k_{1} =0}^{\infty}\prod_{j=1}^{m}2^{k_{j}\theta\alpha_{j}}(k_{j} + 1)^{1-\theta b_{j}}Y_{2^{k_{1}},\ldots,2^{k_{m}}}^{\theta}(f)_{p, \tau}
\\
\sum\limits_{k_{m} =0}^{\infty}...\sum\limits_{k_{1} =0}^{\infty}\prod_{j=1}^{m}2^{k_{j}\theta\alpha_{j}}(k_{j} + 1)^{1-\theta b_{j}}\omega_{\mathbf{r}}^{\theta}(f, 2^{-k_{1}},\ldots, 2^{-k_{m}})_{p, \tau},
\end{gathered}
\end{equation}
where  $\mathbf{r}=(r_{1},\ldots,r_{m})$, $r_{j}>\alpha_{j}$, $j=1,\ldots ,m$.
Now, from inequalities \eqref{eq4.12}, \eqref{eq4.13}, and Lemma 2.7, it follows that
  \begin{equation}\label{eq4.14}
   \begin{gathered} 
 \|f\|_{\text{Lip}_{p, \tau, \theta}^{(\boldsymbol{\alpha}, -\mathbf{b})}}\gg \|f\|_{p, \tau}
+\biggl(\sum\limits_{k_{m} =0}^{\infty}...\sum\limits_{k_{1} =0}^{\infty}\prod_{j=1}^{m}2^{k_{j}\theta\alpha_{j}}(k_{j} + 1)^{1-\theta b_{j}}\omega_{\mathbf{r}}^{\theta}(f, 2^{-k_{1}},\ldots, 2^{-k_{m}})_{p, \tau} \biggr)^{1/\theta}
\\
\gg \|f\|_{p, \tau} + \biggl(\int_{0}^{1}...\int_{0}^{1}\omega_{\mathbf{r}}^{\theta}(f, \mathbf{t})_{p, \tau}\prod_{j=1}^{m}\frac{(1-\log t_{j})^{(\frac{1}{\theta}- b_{j})\theta}}{t_{j}^{\alpha_{j}\theta+1}}dt_{j}  \biggr)^{1/\theta}
\end{gathered}
\end{equation}
in the case  $\sigma\leqslant\theta$, $r_{j}>\alpha_{j}$, $j=1,\ldots ,m$.
Inequality \eqref{eq4.14} means that
  \begin{equation*}
\text{Lip}_{p, \tau, \theta}^{(\boldsymbol{\alpha}, -\mathbf{b})}(\mathbb{T}^{m}) \hookrightarrow 
 S_{p, \tau, \theta}^{\boldsymbol{\alpha}, -\mathbf{b}+\frac{1}{\max\{2, \tau, \theta\}}\overline{e}}B(\mathbb{T}^{m}), 
\end{equation*}  
in the case  $\sigma\leqslant\theta$.
 
 Let $0<\theta < \sigma=\max\{2, \tau\}$.  Then, considering that $b_{j}\theta> 1$ for $j=1,\ldots ,m$ and $\frac{\theta}{\sigma}$ using Hardy's inequality in \cite[Theorem 346]{20}, we will have
 \begin{equation}\label{eq4.15}
   \begin{gathered}
\sum\limits_{\nu_{m} =0}^{\infty}...\sum\limits_{\nu_{1} =0}^{\infty} \prod_{j=1}^{m}(\nu_{j} + 1)^{-\theta b_{j}} \Biggl(\sum\limits_{k_{1}=0}^{\nu_{1}} \ldots \sum\limits_{k_{m}=0}^{\nu_{m}} \prod_{j=1}^{m}2^{k_{j}\sigma\alpha_{j}}Y_{2^{k_{1}},\ldots,2^{k_{m}}}^{\sigma}(f)_{p, \tau}\Biggr)^{\frac{\theta}{\sigma}}
\\
\gg \sum\limits_{\nu_{m} =0}^{\infty}...\sum\limits_{\nu_{1} =0}^{\infty} \prod_{j=1}^{m}2^{\nu_{j}\theta\alpha_{j}}(\nu_{j} + 1)^{\theta(\frac{1}{\sigma} - b_{j})}Y_{2^{k_{1}},\ldots,2^{k_{m}}}^{\theta}(f)_{p, \tau}.
\end{gathered}
\end{equation}
Now, from inequalities \eqref{eq4.10}, \eqref{eq4.15}, and relation \eqref{eq4.13}, it follows that
 \begin{equation*} 
 \|f\|_{\text{Lip}_{p, \tau, \theta}^{(\boldsymbol{\alpha}, -\mathbf{b})}}\gg \|f\|_{p, \tau}
+ \biggl(\int_{0}^{1}...\int_{0}^{1}\omega_{\mathbf{r}}^{\theta}(f, \mathbf{t})_{p, \tau}\prod_{j=1}^{m}\frac{(1-\log t_{j})^{(\frac{1}{\sigma}- b_{j})\theta}}{t_{j}^{\alpha_{j}\theta+1}}dt_{j}  \biggr)^{1/\theta}
\end{equation*}
in the case $0<\theta < \sigma=\max\{2, \tau\}$. Hence 
 \begin{equation}\label{eq4.16}
 \text{Lip}_{p, \tau, \theta}^{(\boldsymbol{\alpha}, -\mathbf{b})}(\mathbb{T}^{m}) \hookrightarrow 
 S_{p, \tau, \theta}^{\boldsymbol{\alpha}, -\mathbf{b}+\frac{1}{\max\{2, \tau, \theta\}}\mathbf{e}}B(\mathbb{T}^{m}),
\end{equation} 
in the case  $0<\theta < \sigma=\max\{2, \tau\}$. This proves the inequality \eqref{eq1.1}. 
 
 Let us prove the inequality \eqref{eq1.2}. In the case $\min\{2, \tau, \theta\}=\theta<\beta=\min\{2, \tau\}$, it follows from the embedding \eqref{eq4.7} that
   \begin{equation}\label{eq4.17}
S_{p, \tau, \min\{2, \tau, \theta\}}^{\boldsymbol{\alpha}, -\mathbf{b}+\frac{1}{\theta}\mathbf{e}}B(\mathbb{T}^{m}) = S_{p, \tau, \theta}^{\boldsymbol{\alpha}, -\mathbf{b}+\frac{1}{\theta}\mathbf{e}}B(\mathbb{T}^{m}) \hookrightarrow \text{Lip}_{p, \tau, \theta}^{(\boldsymbol{\alpha}, -\mathbf{b})}(\mathbb{T}^{m}).
 \end{equation}
Let $\min\{2, \tau, \theta\}=\beta=\min\{2, \tau\}$. Then $\beta\leqslant \theta$. Therefore, by Jensen's inequality, we have
     \begin{equation*} 
   \begin{gathered}
 \biggl(\sum\limits_{k_{m} =0}^{\infty}...\sum\limits_{k_{1} =0}^{\infty}\prod_{j=1}^{m}2^{k_{j}\theta\alpha_{j}}(k_{j} + 1)^{1-\theta b_{j}}\omega_{\mathbf{r}}^{\theta}(f, 2^{-k_{1}},\ldots, 2^{-k_{m}})_{p, \tau} \biggr)^{1/\theta}
 \\
  \ll \biggl(\sum\limits_{k_{m} =0}^{\infty}...\sum\limits_{k_{1} =0}^{\infty}\prod_{j=1}^{m}2^{k_{j}\beta\alpha_{j}}(k_{j} + 1)^{(\frac{1}{\theta} - b_{j})\beta}\omega_{\mathbf{r}}^{\beta}(f, 2^{-k_{1}},\ldots, 2^{-k_{m}})_{p, \tau} \biggr)^{1/\beta}.
\end{gathered}
\end{equation*}
This means that  $S_{p, \tau, \min\{2, \tau, \theta\}}^{\boldsymbol{\alpha}, -\mathbf{b}+\frac{1}{\theta}\mathbf{e}}B(\mathbb{T}^{m})\hookrightarrow S_{p, \tau,  \theta}^{\boldsymbol{\alpha}, -\mathbf{b}+\frac{1}{\theta}\mathbf{e}}B(\mathbb{T}^{m})$,
 in the case   $\min\{2, \tau, \theta\}=\beta$. 
Therefore, from \eqref{eq4.17} it follows that
    $S_{p, \tau, \min\{2, \tau, \theta\}}^{\boldsymbol{\alpha}, -\mathbf{b}+\frac{1}{\theta}\mathbf{e}}B(\mathbb{T}^{m}) \hookrightarrow \text{Lip}_{p, \tau, \theta}^{(\boldsymbol{\alpha}, -\mathbf{b})}(\mathbb{T}^{m})$ and in the case   $\min\{2, \tau, \theta\}=\beta$.
  
Let $\max\{2, \tau, \theta\}=\theta$. Then $\sigma =\max\{2, \tau\}\leqslant \theta$.
Hence 
    \begin{equation}\label{eq4.18}
 \text{Lip}_{p, \tau, \theta}^{(\boldsymbol{\alpha}, -\mathbf{b})}(\mathbb{T}^{m}) \hookrightarrow S_{p, \tau, \max\{2, \tau, \theta\}}^{\boldsymbol{\alpha}, -\mathbf{b}+\frac{1}{\theta}\mathbf{e}}B(\mathbb{T}^{m}).
 \end{equation}
 Let $\max\{2, \tau, \theta\}=\sigma=\max\{2, \tau\}$. Then $0<\theta <\sigma$.
  Therefore, according to Jensen's inequality, we have
   \begin{equation*} 
   \begin{gathered}  
  \biggl(\sum\limits_{k_{m} =0}^{\infty}...\sum\limits_{k_{1} =0}^{\infty}\prod_{j=1}^{m}2^{k_{j}\sigma\alpha_{j}}(k_{j} + 1)^{(\frac{1}{\theta} - b_{j})\sigma}\omega_{\mathbf{r}}^{\sigma}(f, 2^{-k_{1}},\ldots, 2^{-k_{m}})_{p, \tau} \biggr)^{1/\sigma}
  \\
  \leqslant  \biggl(\sum\limits_{k_{m} =0}^{\infty}...\sum\limits_{k_{1} =0}^{\infty}\prod_{j=1}^{m}2^{k_{j}\theta\alpha_{j}}(k_{j} + 1)^{1-\theta b_{j}}\omega_{\mathbf{r}}^{\theta}(f, 2^{-k_{1}},\ldots, 2^{-k_{m}})_{p, \tau} \biggr)^{1/\theta}.
 \end{gathered}
\end{equation*}
This means that  $S_{p, \tau,  \theta}^{\boldsymbol{\alpha}, -\mathbf{b}+\frac{1}{\theta}\mathbf{e}}B(\mathbb{T}^{m}) \hookrightarrow S_{p, \tau,  \sigma}^{\boldsymbol{\alpha}, -\mathbf{b}+\frac{1}{\theta}\mathbf{e}}B(\mathbb{T}^{m})$. 
Therefore, it follows from \eqref{eq4.18} that $\text{Lip}_{p, \tau, \theta}^{(\boldsymbol{\alpha}, -\mathbf{b})}(\mathbb{T}^{m}) \hookrightarrow S_{p, \tau, \max\{2, \tau, \theta\}}^{\boldsymbol{\alpha}, -\mathbf{b}+\frac{1}{\theta}\mathbf{e}}B(\mathbb{T}^{m})$ 
and in the case $\max\{2, \tau, \theta\}=\sigma=\max\{2, \tau\}$. 
  Theorem 1.2 is proven.

\smallskip

\textbf{Proof of Theorem 1.3.} 
Let  $f \in S_{p_{0}, \tau_{0}, \theta}^{\boldsymbol{\alpha}+(\frac{1}{p_{0}}-\frac{1}{p})\mathbf{e}, -\mathbf{b}+\frac{1}{\min\{\tau, \theta\}}\mathbf{e}}B(\mathbb{T}^{m})$. Then   
\begin{equation}\label{eq4.19}
\|f\|_{p_{0}, \tau_{0}}+ \biggl(\sum\limits_{\mathbf{s}\in \mathbb{Z}_{+}^{m}}\prod_{j=1}^{m}2^{s_{j}(\alpha_{j}+\frac{1}{p_{0}}-\frac{1}{p})\theta}(s_{j}+1)^{(-b_{j}+\frac{1}{\min\{\tau, \theta\}})\theta}\|\delta_{\mathbf{s}}(f)\|_{p_{0}, \tau_{0}}^{\theta} \biggr)^{1/\theta} <+\infty. 
\end{equation}
According to Theorem 3.1 in \cite{2}, the following inequality holds:
 \begin{equation}\label{eq4.20}
\|f\|_{p, \tau} \ll \biggl(\sum\limits_{\mathbf{s}\in \mathbb{Z}_{+}^{m}}\prod_{j=1}^{m}2^{s_{j}(\frac{1}{p_{0}}-\frac{1}{p})\tau}\|\delta_{\mathbf{s}}(f)\|_{p_{0}, \tau_{0}}^{\tau} \biggr)^{1/\tau}
\end{equation}
for the function  $f\in L_{p_{0}, \tau_{0}}(\mathbb{T}^{m})$, $1<p_{0}<p<\infty$, $1<\tau_{0}, \tau<\infty$. 

If $0<\theta \leqslant \tau$, then according to Jensen's inequality in \cite[Ch.~3, Sec.~3]{27} from \eqref{eq4.20}, we will have
\begin{equation}\label{eq4.21}
\|f\|_{p, \tau} \ll \biggl(\sum\limits_{\mathbf{s}\in \mathbb{Z}_{+}^{m}}\prod_{j=1}^{m}2^{s_{j}(\frac{1}{p_{0}}-\frac{1}{p})\theta}\|\delta_{\mathbf{s}}(f)\|_{p_{0}, \tau_{0}}^{\theta} \biggr)^{1/\theta}.
\end{equation}
Since $\alpha_{j}> 0$ for $j=1,\ldots ,m$,  from inequality \eqref{eq4.21} it follows  that
\begin{equation}\label{eq4.22}
\|f\|_{p, \tau} \ll \biggl(\sum\limits_{\mathbf{s}\in \mathbb{Z}_{+}^{m}}\prod_{j=1}^{m}2^{s_{j}(\alpha_{j}+\frac{1}{p_{0}}-\frac{1}{p})\theta}(s_{j}+1)^{(-b_{j}+\frac{1}{\theta})\theta}\|\delta_{\mathbf{s}}(f)\|_{p_{0}, \tau_{0}}^{\theta} \biggr)^{1/\theta}
\end{equation}
in the case  $0<\theta \leqslant \tau<\infty$.

Let $1<\tau<\theta \leqslant \infty$. Then, applying H\"{o}lder's inequality with $\eta=\frac{\theta}{\tau}$, $\frac{1}{\eta}+\frac{1}{\eta'}=1$ from \eqref{eq4.20}, we obtain
  \begin{equation}\label{eq4.23} 
   \begin{gathered}
\|f\|_{p, \tau} \ll \biggl(\sum\limits_{\mathbf{s}\in \mathbb{Z}_{+}^{m}}\prod_{j=1}^{m}2^{s_{j}(\alpha_{j}+\frac{1}{p_{0}}-\frac{1}{p})\theta}(s_{j}+1)^{(-b_{j}+\frac{1}{\tau})\theta}\|\delta_{\mathbf{s}}(f)\|_{p_{0}, \tau_{0}}^{\theta} \biggr)^{1/\theta}
\\
\times \biggl(\sum\limits_{\mathbf{s}\in \mathbb{Z}_{+}^{m}}\prod_{j=1}^{m}2^{-s_{j}\alpha_{j}\tau\eta'}(s_{j}+1)^{(b_{j}-\frac{1}{\tau})\tau\eta'}\biggr)^{\frac{1}{\tau\eta'}}.
 \end{gathered}
\end{equation}
Since $\alpha_{j}> 0$ for $j=1,\ldots ,m$, then 
\begin{equation*}
\sum\limits_{\mathbf{s}\in \mathbb{Z}_{+}^{m}}\prod_{j=1}^{m}2^{-s_{j}\alpha_{j}\tau\eta'}(s_{j}+1)^{(b_{j}-\frac{1}{\tau})\tau\eta'}<+\infty.
\end{equation*}
Therefore,  from \eqref{eq4.23} it follows that $S_{p_{0}, \tau_{0}, \theta}^{\boldsymbol{\alpha}+(\frac{1}{p_{0}}-\frac{1}{p})\mathbf{e}, -\mathbf{b}+\frac{1}{\tau}\mathbf{e}}B(\mathbb{T}^{m})\hookrightarrow L_{p, \tau}(\mathbb{T}^{m})$ in the case $1<\tau<\theta \leqslant \infty$. And in the case $0<\theta\leqslant\tau< \infty$,  from \eqref{eq4.23} it follows that $S_{p_{0}, \tau_{0}, \theta}^{\boldsymbol{\alpha}+(\frac{1}{p_{0}}-\frac{1}{p})\mathbf{e}, -\mathbf{b}+\frac{1}{\theta}\mathbf{e}}B(\mathbb{T}^{m})\hookrightarrow L_{p, \tau}(\mathbb{T}^{m})$.
Further, using Theorem 1.1 in \cite{2} and inequality \eqref{eq4.20}, we obtain
 \begin{equation}\label{eq4.24} 
   \begin{gathered}
\Biggl\|\Biggl(\sum\limits_{s_{m}=2^{l_{m}}}^{2^{l_{m}+1}-1}...\sum\limits_{s_{1}=2^{l_{1}}}^{2^{l_{1}+1}-1}\prod_{j=1}^{m}2^{s_{j}\alpha_{j}2}|\delta_{\mathbf{s}}(f)|^{2}\Biggr)^{\frac{1}{2}}\Biggr\|_{p, \tau} \asymp \Biggl\|\sum\limits_{s_{m}=2^{l_{m}}}^{2^{l_{m}+1}-1}...\sum\limits_{s_{1}=2^{l_{1}}}^{2^{l_{1}+1}-1}\prod_{j=1}^{m}2^{s_{j}\alpha_{j}}\delta_{\mathbf{s}}(f)\Biggr\|_{p, \tau}
\\
\ll \Biggl(\sum\limits_{s_{m}=2^{l_{m}}}^{2^{l_{m}+1}-1}...\sum\limits_{s_{1}=2^{l_{1}}}^{2^{l_{1}+1}-1} \prod_{j=1}^{m}2^{s_{j}(\alpha_{j}+
\frac{1}{p_{0}}-\frac{1}{p})\tau}\|\delta_{\mathbf{s}}(f)\|_{p_{0}, \tau_{0}}^{\tau}\Biggr)^{\frac{1}{\tau}}.
 \end{gathered}
\end{equation}
Now, applying Theorem 1.1 and \eqref{eq4.24}, we obtain
\begin{equation}\label{eq4.25} 
   \begin{gathered}
\|f\|_{\text{Lip}_{p, \tau, \theta}^{(\boldsymbol{\alpha}, -\mathbf{b})}}\ll \|f\|_{p, \tau}
\\
 + \biggl(\sum\limits_{l_{m} =0}^{\infty}...\sum\limits_{l_{1} =0}^{\infty}\prod_{j=1}^{m}2^{l_{j}(1-b_{j}\theta)}\Biggl(\sum\limits_{s_{m}=2^{l_{m}}}^{2^{l_{m}+1}-1}...\sum\limits_{s_{1}=2^{l_{1}}}^{2^{l_{1}+1}-1} \prod_{j=1}^{m}2^{s_{j}(\alpha_{j}+
\frac{1}{p_{0}}-\frac{1}{p})\tau}\|\delta_{\mathbf{s}}(f)\|_{p_{0}, \tau_{0}}^{\tau}\Biggr)^{\frac{\theta}{\tau}}\biggr)^{1/\theta}. 
\end{gathered}
\end{equation}
If $0<\theta \leqslant \tau < \infty$, then by Jensen's inequality from \eqref{eq4.25}  we obtain  
\begin{equation}\label{eq4.26} 
   \begin{gathered}
\|f\|_{\text{Lip}_{p, \tau, \theta}^{(\boldsymbol{\alpha}, -\mathbf{b})}}\ll \|f\|_{p, \tau}
\\
 +
\biggl(\sum\limits_{l_{m} =0}^{\infty}...\sum\limits_{l_{1} =0}^{\infty}\prod_{j=1}^{m}2^{l_{j}(1-b_{j}\theta)}\sum\limits_{s_{m}=2^{l_{m}}}^{2^{l_{m}+1}-1}...\sum\limits_{s_{1}=2^{l_{1}}}^{2^{l_{1}+1}-1} \prod_{j=1}^{m}2^{s_{j}(\alpha_{j}+
\frac{1}{p_{0}}-\frac{1}{p})\theta}\|\delta_{\mathbf{s}}(f)\|_{p_{0}, \tau_{0}}^{\theta}\biggr)^{1/\theta} 
\\
\ll \|f\|_{p, \tau} + \biggl(\sum\limits_{\mathbf{s}\in \mathbb{Z}_{+}^{m}}\prod_{j=1}^{m}2^{s_{j}(\alpha_{j}+\frac{1}{p_{0}}-\frac{1}{p})\theta}(s_{j}+1)^{(\frac{1}{\theta}-b_{j})\theta}\|\delta_{\mathbf{s}}(f)\|_{p_{0}, \tau_{0}}^{\theta} \biggr)^{1/\theta}.
\end{gathered}
\end{equation}
From the inequalities \eqref{eq4.22}, \eqref{eq4.23} and \eqref{eq4.26} it follows that
\begin{equation*}
\|f\|_{\text{Lip}_{p, \tau, \theta}^{(\boldsymbol{\alpha}, -\mathbf{b})}}\ll \|f\|_{p_{0}, \tau_{0}} + \biggl(\sum\limits_{\mathbf{s}\in \mathbb{Z}_{+}^{m}}\prod_{j=1}^{m}2^{s_{j}(\alpha_{j}+\frac{1}{p_{0}}-\frac{1}{p})\theta}(s_{j}+1)^{(\frac{1}{\theta}-b_{j})\theta}\|\delta_{\mathbf{s}}(f)\|_{p_{0}, \tau_{0}}^{\theta} \biggr)^{1/\theta}
\end{equation*}
in the case   $0<\theta \leqslant \tau < \infty$ for the function  
  $f \in S_{p_{0}, \tau_{0}, \theta}^{\boldsymbol{\alpha}+(\frac{1}{p_{0}}-\frac{1}{p})\mathbf{e}, -\mathbf{b}+\frac{1}{\min\{\tau, \theta\}}\mathbf{e}}B(\mathbb{T}^{m})$. Hence            
   $S_{p_{0}, \tau_{0}, \theta}^{\boldsymbol{\alpha}+(\frac{1}{p_{0}}-\frac{1}{p})\mathbf{e}, -\mathbf{b}+\frac{1}{\min\{\tau, \theta\}}\mathbf{e}}B(\mathbb{T}^{m}) \hookrightarrow \text{Lip}_{p, \tau, \theta}^{(\boldsymbol{\alpha}, -\mathbf{b})}(\mathbb{T}^{m})$ in the case 
    $0<\theta \leqslant \tau < \infty$.

Let $1<\tau<\theta \leqslant \infty$ and set $\eta=\frac{\theta}{\tau}$, $\frac{1}{\eta}+\frac{1}{\eta'}=1$ (in the case of $\theta =\infty$, we assume that $\eta'=1$). Then, applying H\"{o}lder's inequality, we have 
\begin{equation}\label{eq4.27} 
   \begin{gathered}
   \Biggl(\sum\limits_{s_{m}=2^{l_{m}}}^{2^{l_{m}+1}-1}...\sum\limits_{s_{1}=2^{l_{1}}}^{2^{l_{1}+1}-1} \prod_{j=1}^{m}2^{s_{j}(\alpha_{j}+
\frac{1}{p_{0}}-\frac{1}{p})\tau}\|\delta_{\mathbf{s}}(f)\|_{p_{0}, \tau_{0}}^{\tau}\Biggr)^{\frac{1}{\tau}}
\\
 \ll \prod_{j=1}^{m}2^{l_{j}(\frac{\theta}{\tau}-1)}\biggl(\sum\limits_{s_{m}=2^{l_{m}}}^{2^{l_{m}+1}-1}...\sum\limits_{s_{1}=2^{l_{1}}}^{2^{l_{1}+1}-1} \prod_{j=1}^{m}2^{s_{j}(\alpha_{j}+
\frac{1}{p_{0}}-\frac{1}{p})\theta}\|\delta_{\mathbf{s}}(f)\|_{p_{0}, \tau_{0}}^{\theta}\biggr)^{1/\theta}.
\end{gathered}
\end{equation}
Now, from inequalities \eqref{eq4.22}, \eqref{eq4.23}, \eqref{eq4.25}, and \eqref{eq4.27}, it follows that
 \begin{equation}\label{eq4.28}
\|f\|_{\text{Lip}_{p, \tau, \theta}^{(\boldsymbol{\alpha}, -\mathbf{b})}}\ll \|f\|_{p_{0}, \tau_{0}} + \biggl(\sum\limits_{\mathbf{s}\in \mathbb{Z}_{+}^{m}}\prod_{j=1}^{m}2^{s_{j}(\alpha_{j}+\frac{1}{p_{0}}-\frac{1}{p})\theta}(s_{j}+1)^{(\frac{1}{\tau}-b_{j})\theta}\|\delta_{\mathbf{s}}(f)\|_{p_{0}, \tau_{0}}^{\theta} \biggr)^{1/\theta}
\end{equation}
for the function  $f \in S_{p_{0}, \tau_{0}, \theta}^{\boldsymbol\alpha+(\frac{1}{p_{0}}-\frac{1}{p})\mathbf{e}, -\mathbf{b}+\frac{1}{\min\{\tau, \theta\}}\mathbf{e}}B(\mathbb{T}^{m})$,
in the case  $1<\tau<\theta \leqslant \infty$. Therefore 
   $S_{p_{0}, \tau_{0}, \theta}^{\boldsymbol{\alpha}+(\frac{1}{p_{0}}-\frac{1}{p})\mathbf{e}, -\mathbf{b}+\frac{1}{\tau}\mathbf{e}}B(\mathbb{T}^{m}) \hookrightarrow \text{Lip}_{p, \tau, \theta}^{(\boldsymbol{\alpha}, -\mathbf{b})}(\mathbb{T}^{m})$ in the case 
    $1<\tau<\theta \leqslant \infty$. This proves the left-hand side of \eqref{eq1.3}.
    
Let  $f \in \text{Lip}_{p, \tau, \theta}^{(\boldsymbol{\alpha}, -\mathbf{b})}(\mathbb{T}^{m})$.
Since $1<p<p_{1}<\infty$, using Theorem 3.1 in \cite{2} and the Nikol'skii inequality for trigonometric polynomials in Lorentz space in \cite[p. 7]{3}, we obtain
 \begin{equation}\label{eq4.29} 
   \begin{gathered}
   \Biggl\|\sum\limits_{s_{m}=2^{l_{m}}}^{2^{l_{m}+1}-1}...\sum\limits_{s_{1}=2^{l_{1}}}^{2^{l_{1}+1}-1}\prod_{j=1}^{m}2^{s_{j}\alpha_{j}}\delta_{\mathbf{s}}(f)\Biggr\|_{p, \tau} \gg 
\Biggl(\sum\limits_{s_{m}=2^{l_{m}}}^{2^{l_{m}+1}-1}...\sum\limits_{s_{1}=2^{l_{1}}}^{2^{l_{1}+1}-1} \prod_{j=1}^{m}2^{s_{j}(\frac{1}{\lambda}-\frac{1}{p})\tau}\|\prod_{j=1}^{m}2^{s_{j}\alpha_{j}}\delta_{\mathbf{s}}(f)\|_{\lambda}^{\tau}\Biggr)^{\frac{1}{\tau}}
\\
\gg \Biggl(\sum\limits_{s_{m}=2^{l_{m}}}^{2^{l_{m}+1}-1}...\sum\limits_{s_{1}=2^{l_{1}}}^{2^{l_{1}+1}-1} \prod_{j=1}^{m}2^{s_{j}(\frac{1}{p_{1}}-\frac{1}{p})\tau}\|\prod_{j=1}^{m}2^{s_{j}\alpha_{j}}\delta_{\mathbf{s}}(f)\|_{p_{1}, \tau_{1}}^{\tau}\Biggr)^{\frac{1}{\tau}}
\end{gathered}
\end{equation}
for  $\lambda\in (p, \, p_{1})$. From inequality \eqref{eq4.29}, it follows that 
 \begin{equation}\label{eq4.30} 
   \begin{gathered}
\biggl(\sum\limits_{l_{m} =0}^{\infty}...\sum\limits_{l_{1} =0}^{\infty}\prod_{j=1}^{m}2^{l_{j}(1-b_{j}\theta)}\Biggl(\sum\limits_{s_{m}=2^{l_{m}}}^{2^{l_{m}+1}-1}...\sum\limits_{s_{1}=2^{l_{1}}}^{2^{l_{1}+1}-1} \prod_{j=1}^{m}2^{s_{j}(\alpha_{j}+\frac{1}{p_{1}}-\frac{1}{p})\tau}\|\delta_{\mathbf{s}}(f)\|_{p_{1}, \tau_{1}}^{\tau}\Biggr)^{\frac{\theta}{\tau}}\biggr)^{1/\theta}
\\
\ll  \biggl(\sum\limits_{l_{m} =0}^{\infty}...\sum\limits_{l_{1} =0}^{\infty}\prod_{j=1}^{m}2^{l_{j}(1-b_{j}\theta)} \Biggl\|\sum\limits_{s_{m}=2^{l_{m}}}^{2^{l_{m}+1}-1}...\sum\limits_{s_{1}=2^{l_{1}}}^{2^{l_{1}+1}-1}\prod_{j=1}^{m}2^{s_{j}\alpha_{j}}\delta_{\mathbf{s}}(f)\Biggr\|_{p, \tau}^{\theta}\biggr)^{1/\theta}
\end{gathered}
\end{equation}
for  $1<p<p_{1}<\infty$, $1<\tau, \tau_{1}<\infty$, $0<\theta\leqslant \infty$.

If $1<\tau\leqslant \theta\leqslant \infty$, then by Jensen's inequality we have
  \begin{equation}\label{eq4.31} 
   \begin{gathered}
   \sum\limits_{s_{m}=2^{l_{m}}}^{2^{l_{m}+1}-1}...\sum\limits_{s_{1}=2^{l_{1}}}^{2^{l_{1}+1}-1} \prod_{j=1}^{m}2^{s_{j}(\alpha_{j}+\frac{1}{p_{1}}-\frac{1}{p})\theta}\|\delta_{\mathbf{s}}(f)\|_{p_{1}, \tau_{1}}^{\theta}
\\
 \leqslant \Biggl(\sum\limits_{s_{m}=2^{l_{m}}}^{2^{l_{m}+1}-1}...\sum\limits_{s_{1}=2^{l_{1}}}^{2^{l_{1}+1}-1} \prod_{j=1}^{m}2^{s_{j}(\alpha_{j}+\frac{1}{p_{1}}-\frac{1}{p})\tau}\|\delta_{\mathbf{s}}(f)\|_{p_{1}, \tau_{1}}^{\tau}\Biggr)^{\frac{\theta}{\tau}}.
\end{gathered}
\end{equation}
Now, from inequalities \eqref{eq4.30} and \eqref{eq4.31}, it follows that
\begin{equation}\label{eq4.32} 
   \begin{gathered}
   \biggl(\sum\limits_{\mathbf{s}\in \mathbb{Z}_{+}^{m}}\prod_{j=1}^{m}2^{s_{j}(\alpha_{j}+\frac{1}{p_{1}}-\frac{1}{p})\theta}(s_{j}+1)^{(\frac{1}{\tau}-b_{j})\theta}\|\delta_{\mathbf{s}}(f)\|_{p_{1}, \tau_{1}}^{\theta} \biggr)^{1/\theta} 
\\
\ll  \biggl(\sum\limits_{l_{m} =0}^{\infty}...\sum\limits_{l_{1} =0}^{\infty}\prod_{j=1}^{m}2^{l_{j}(\frac{1}{\theta}-b_{j})\theta} \Biggl\|\sum\limits_{s_{m}=2^{l_{m}}}^{2^{l_{m}+1}-1}...\sum\limits_{s_{1}=2^{l_{1}}}^{2^{l_{1}+1}-1}\prod_{j=1}^{m}2^{s_{j}\alpha_{j}}\delta_{\mathbf{s}}(f)\Biggr\|_{p, \tau}^{\theta}\biggr)^{1/\theta}
\end{gathered}
\end{equation}
in the case   $1<\tau\leqslant \theta\leqslant \infty$.
Therefore, according to Theorem 1.1, we have 
\begin{equation*}
\text{Lip}_{p, \tau, \theta}^{(\boldsymbol{\alpha}, -\mathbf{b})}(\mathbb{T}^{m})\hookrightarrow S_{p_{1}, \tau_{1}, \theta}^{\boldsymbol{\alpha}+(\frac{1}{p_{1}}-\frac{1}{p})\mathbf{e}, -\mathbf{b}+\frac{1}{\theta}\mathbf{e}}B(\mathbb{T}^{m})
\end{equation*} 
in the case $1<\tau\leqslant \theta\leqslant \infty$.
  
Let $0<\theta<\tau <\infty$. Then, applying H\"{o}lder's inequality with $\eta=\frac{\tau}{\theta}$, $\frac{1}{\eta}+\frac{1}{\eta'}=1$, we have
 \begin{equation}\label{eq4.33} 
   \begin{gathered}
   \Biggl(\sum\limits_{s_{m}=2^{l_{m}}}^{2^{l_{m}+1}-1}...\sum\limits_{s_{1}=2^{l_{1}}}^{2^{l_{1}+1}-1} \prod_{j=1}^{m}2^{s_{j}(\alpha_{j}+
\frac{1}{p_{1}}-\frac{1}{p})\theta}\|\delta_{\mathbf{s}}(f)\|_{p_{1}, \tau_{1}}^{\theta}\Biggr)^{\frac{1}{\theta}}
\\
 \ll \prod_{j=1}^{m}2^{l_{j}(\frac{1}{\theta}-\frac{1}{\tau})}\biggl(\sum\limits_{s_{m}=2^{l_{m}}}^{2^{l_{m}+1}-1}...\sum\limits_{s_{1}=2^{l_{1}}}^{2^{l_{1}+1}-1} \prod_{j=1}^{m}2^{s_{j}(\alpha_{j}+
\frac{1}{p_{1}}-\frac{1}{p})\tau}\|\delta_{\mathbf{s}}(f)\|_{p_{1}, \tau_{1}}^{\tau}\biggr)^{1/\tau}.
\end{gathered}
\end{equation}
Now from the inequalities \eqref{eq4.30} and \eqref{eq4.34} it follows that
\begin{equation*}
\begin{gathered}
\biggl(\sum\limits_{\mathbf{s}\in \mathbb{Z}_{+}^{m}}\prod_{j=1}^{m}2^{s_{j}(\alpha_{j}+\frac{1}{p_{1}}-\frac{1}{p})\theta}(s_{j}+1)^{(\frac{1}{\tau}-b_{j})\theta}\|\delta_{\mathbf{s}}(f)\|_{p_{1}, \tau_{1}}^{\theta} \biggr)^{1/\theta} 
\\
\ll  \biggl(\sum\limits_{l_{m} =0}^{\infty}...\sum\limits_{l_{1} =0}^{\infty}\prod_{j=1}^{m}2^{l_{j}(\frac{1}{\theta}-b_{j})\theta} \Biggl\|\sum\limits_{s_{m}=2^{l_{m}}}^{2^{l_{m}+1}-1}...\sum\limits_{s_{1}=2^{l_{1}}}^{2^{l_{1}+1}-1}\prod_{j=1}^{m}2^{s_{j}\alpha_{j}}\delta_{\mathbf{s}}(f)\Biggr\|_{p, \tau}^{\theta}\biggr)^{1/\theta}
\end{gathered}
\end{equation*}
in the case   $0<\theta<\tau <\infty$.
Hence,  
 $\text{Lip}_{p, \tau, \theta}^{(\boldsymbol{\alpha}, -\mathbf{b})}(\mathbb{T}^{m})\hookrightarrow S_{p_{1}, \tau_{1}, \theta}^{\boldsymbol{\alpha}+(\frac{1}{p_{1}}-\frac{1}{p})\mathbf{e}, -\mathbf{b}+\frac{1}{\tau}\mathbf{e}}B(\mathbb{T}^{m})$ in the case 
  $0<\theta<\tau <\infty$. This proves the inequality \eqref{eq1.3}.
   
Let us prove the embeddings \eqref{eq1.4}.
  Let 
 $f\in S_{p_{0}, \tau_{0}, \min\{\tau, \theta\}}^{\boldsymbol{\alpha}+(\frac{1}{p_{0}}-\frac{1}{p})\mathbf{e}, -\mathbf{b}+\frac{1}{\theta}\mathbf{e}}B(\mathbb{T}^{m})$.
If $0<\theta \leqslant \tau<\infty$, then from the left inclusion in \eqref{eq1.3} it follows that
 \begin{equation*}
S_{p_{0}, \tau_{0}, \min\{\tau, \theta\}}^{\boldsymbol{\alpha}+(\frac{1}{p_{0}}-\frac{1}{p})\mathbf{e}, -\mathbf{b}+\frac{1}{\theta}\mathbf{e}}B(\mathbb{T}^{m})=S_{p_{0}, \tau_{0}, \theta}^{\boldsymbol{\alpha}+(\frac{1}{p_{0}}-\frac{1}{p})\mathbf{e}, -\mathbf{b}+\frac{1}{\theta}\mathbf{e}}B(\mathbb{T}^{m}) \hookrightarrow \text{Lip}_{p, \tau, \theta}^{(\boldsymbol{\alpha}, -\mathbf{b})}(\mathbb{T}^{m}).
\end{equation*}
Let $1< \tau<\theta \leqslant \infty$. Since $0<\frac{\tau}{\theta}<1$, by Jensen's inequality and Theorem 1.1, we obtain
\begin{equation}\label{eq4.34} 
   \begin{gathered}
\|f\|_{\text{Lip}_{p, \tau, \theta}^{(\boldsymbol{\alpha}, -\mathbf{b})}}\ll \|f\|_{p, \tau}
+ \biggl(\sum\limits_{\mathbf{l}\in \mathbb{Z}_{+}^{m}}\prod_{j=1}^{m}2^{l_{j}(\frac{1}{\theta}-b_{j})\tau} \Biggl\|\sum\limits_{s_{m}=2^{l_{m}}}^{2^{l_{m}+1}-1}...\sum\limits_{s_{1}=2^{l_{1}}}^{2^{l_{1}+1}-1}\prod_{j=1}^{m}2^{s_{j}\alpha_{j}}\delta_{\mathbf{s}}(f)\Biggr\|_{p, \tau}^{\tau}\biggr)^{1/\tau}.
\end{gathered}
\end{equation}
Further, using inequality \eqref{eq4.24} from \eqref{eq4.34}, we obtain
\begin{equation*}
\|f\|_{\text{Lip}_{p, \tau, \theta}^{(\boldsymbol{\alpha}, -\mathbf{b})}}\ll \|f\|_{p_{0}, \tau_{0}} + \biggl(\sum\limits_{\mathbf{s}\in \mathbb{Z}_{+}^{m}}\prod_{j=1}^{m}2^{s_{j}(\alpha_{j}+\frac{1}{p_{0}}-\frac{1}{p})\tau}(s_{j}+1)^{(\frac{1}{\theta}-b_{j})\tau}\|\delta_{\mathbf{s}}(f)\|_{p_{0}, \tau_{0}}^{\tau} \biggr)^{1/\tau} 
\end{equation*}
in the case  $1< \tau<\theta \leqslant \infty$.
Hence, $S_{p_{0}, \tau_{0}, \tau}^{\boldsymbol{\alpha}+(\frac{1}{p_{0}}-\frac{1}{p})\mathbf{e}, -\mathbf{b}+\frac{1}{\theta}\mathbf{e}}B(\mathbb{T}^{m}) \hookrightarrow \text{Lip}_{p, \tau, \theta}^{(\boldsymbol{\alpha}, -\mathbf{b})}(\mathbb{T}^{m})$ in the case  $1< \tau<\theta \leqslant \infty$.

Let 
 $f\in \text{Lip}_{p, \tau, \theta}^{(\boldsymbol{\alpha}, -\mathbf{b})}(\mathbb{T}^{m})$.
If $1< \tau<\theta \leqslant \infty$, then from the right inclusion in \eqref{eq1.3}  
it follows that
$\text{Lip}_{p, \tau, \theta}^{(\boldsymbol{\alpha}, -\mathbf{b})}(\mathbb{T}^{m}) \hookrightarrow 
S_{p_{1}, \tau_{1},  \theta}^{\boldsymbol{\alpha}+(\frac{1}{p_{1}}-\frac{1}{p})\mathbf{e}, -\mathbf{b}+\frac{1}{\theta}\mathbf{e}}B(\mathbb{T}^{m})=S_{p_{1}, \tau_{1}, \min\{\tau, \theta\}}^{\boldsymbol{\alpha}+(\frac{1}{p_{1}}-\frac{1}{p})\mathbf{e}, -\mathbf{b}+\frac{1}{\theta}\mathbf{e}}B(\mathbb{T}^{m}).
$
Let $0<\theta \leqslant \tau<\infty$. Then, by Jensen's inequality, we have
\begin{equation}\label{eq4.35} 
   \begin{gathered}
   \biggl(\sum\limits_{l_{m} =0}^{\infty}...\sum\limits_{l_{1} =0}^{\infty}\prod_{j=1}^{m}2^{l_{j}(\frac{1}{\theta}-b_{j})\tau} \Biggl\|\sum\limits_{s_{m}=2^{l_{m}}}^{2^{l_{m}+1}-1}...\sum\limits_{s_{1}=2^{l_{1}}}^{2^{l_{1}+1}-1}\prod_{j=1}^{m}2^{s_{j}\alpha_{j}}\delta_{\mathbf{s}}(f)\Biggr\|_{p, \tau}^{\tau}\biggr)^{1/\tau}
\\
\ll
\biggl(\sum\limits_{l_{m} =0}^{\infty}...\sum\limits_{l_{1} =0}^{\infty}\prod_{j=1}^{m}2^{l_{j}(\frac{1}{\theta}-b_{j})\theta} \Biggl\|\sum\limits_{s_{m}=2^{l_{m}}}^{2^{l_{m}+1}-1}...\sum\limits_{s_{1}=2^{l_{1}}}^{2^{l_{1}+1}-1}\prod_{j=1}^{m}2^{s_{j}\alpha_{j}}\delta_{\mathbf{s}}(f)\Biggr\|_{p, \tau}^{\theta}\biggr)^{1/\theta}.
\end{gathered}
\end{equation}
Now, from inequalities \eqref{eq4.29} and \eqref{eq4.35}, it follows that
\begin{equation}\label{eq4.36} 
   \begin{gathered}
   \biggl(\sum\limits_{\mathbf{s}\in \mathbb{Z}_{+}^{m}}\prod_{j=1}^{m}2^{s_{j}(\alpha_{j}+\frac{1}{p_{1}}-\frac{1}{p})\tau}(s_{j}+1)^{(\frac{1}{\theta}-b_{j})\tau}\|\delta_{\mathbf{s}}(f)\|_{p_{1}, \tau_{1}}^{\tau} \biggr)^{1/\tau} 
\\
\ll  \biggl(\sum\limits_{l_{m} =0}^{\infty}...\sum\limits_{l_{1} =0}^{\infty}\prod_{j=1}^{m}2^{l_{j}(\frac{1}{\theta}-b_{j})\tau} \Biggl\|\sum\limits_{s_{m}=2^{l_{m}}}^{2^{l_{m}+1}-1}...\sum\limits_{s_{1}=2^{l_{1}}}^{2^{l_{1}+1}-1}\prod_{j=1}^{m}2^{s_{j}\alpha_{j}}\delta_{\mathbf{s}}(f)\Biggr\|_{p, \tau}^{\tau}\biggr)^{1/\tau}.
\end{gathered}
\end{equation}
Furthermore, from inequalities \eqref{eq4.35}, \eqref{eq4.36} and Theorem 1.1, it follows that
\begin{equation*}
\text{Lip}_{p, \tau, \theta}^{(\boldsymbol{\alpha}, -\mathbf{b})}(\mathbb{T}^{m}) \hookrightarrow 
S_{p_{1}, \tau_{1},  \tau}^{\boldsymbol{\alpha}+(\frac{1}{p_{1}}-\frac{1}{p})\mathbf{e}, -\mathbf{b}+\frac{1}{\theta}\mathbf{e}}B(\mathbb{T}^{m}),
\end{equation*}
in the case  $0<\theta \leqslant \tau<\infty$. Theorem 1.3 is proven. 
  \hfill $\Box$

\setcounter{equation}{0}
\setcounter{lemma}{0}
\setcounter{theorem}{0}

\section{Proofs of Theorem 1.4 and Theorem 1.5}\label{5}

\textbf{Proof of Theorem 1.4.} 
Let $f \in \text{Lip}_{p, \tau, \theta_{0}}^{(\boldsymbol{\alpha}^{(0)}, -\mathbf{b}^{(0)})}(\mathbb{T}^{m})$.  For the Nikol'skii--Besov space, it is well known that
 \begin{equation}\label{eq5.1}
S_{p, \tau, \theta_{0}}^{\boldsymbol{\alpha}^{(0)}, \boldsymbol{\xi}^{(0)}}B(\mathbb{T}^{m}) \hookrightarrow S_{p, \tau, \theta_{1}}^{\boldsymbol{\alpha}^{(1)}, \boldsymbol{\xi}^{(1)}}B(\mathbb{T}^{m})
\end{equation}
for  $0< \alpha_{j}^{(1)}<\alpha_{j}^{(0)}<\infty$, $\xi_{j}^{(0)}, \xi_{j}^{(1)}\in \mathbb{R}$, $0<\theta_{i}\leqslant \infty$, $i=1, 2$ and $j=1,\ldots , m$.

This follows from the well-known relation
   \begin{equation}\label{eq5.2}
\|f\|_{S_{p, \tau, \theta}^{\mathbf{r}, \mathbf{b}}B} \asymp \biggl(\sum\limits_{\mathbf{s}\in \mathbb{Z}_{+}^{m}}\prod_{j=1}^{m}2^{s_{j}r_{j}\theta}(s_{j}+1)^{b_{j}\theta}\|\delta_{\mathbf{s}}(f)\|_{p_{1}, \tau_{1}}^{\theta} \biggr)^{1/\theta} 
\end{equation}
for $r_{j}>0$, $b_{j}\in \mathbb{R}$, $j=1,\ldots , m$, $1<p<\infty$, $1<\tau<\infty$, $0< \theta \leqslant \infty$, because if $\alpha_{j}^{(0)}>\alpha_{j}^{(1)}$ for $j=1,\ldots , m$, then
 \begin{equation*}
\prod_{j=1}^{m}2^{s_{j}\alpha_{j}^{(0)}}(s_{j}+1)^{b_{j}^{(0)}} > \prod_{j=1}^{m}2^{s_{j}\alpha_{j}^{(1)}}(s_{j}+1)^{b_{j}^{(1)}}.
\end{equation*}
From the embeddings \eqref{eq1.1} and \eqref{eq5.1} it follows that
\begin{equation*}
\text{Lip}_{p, \tau, \theta_{0}}^{(\boldsymbol{\alpha}^{(0)}, -\mathbf{b}^{(0)})}(\mathbb{T}^{m}) \hookrightarrow S_{p, \tau, \theta_{0}}^{\boldsymbol{\alpha}^{(0)}, -\mathbf{b}^{(0)}+\frac{1}{\max\{2, \tau, \theta_{0}\}}\mathbf{e}}B(\mathbb{T}^{m}) 
\end{equation*}
\begin{equation*}
\hookrightarrow
S_{p, \tau, \theta_{1}}^{\boldsymbol{\alpha}^{(1)}, -\mathbf{b}^{(1)}+\frac{1}{\min\{2, \tau, \theta_{1}\}}\mathbf{e}}B(\mathbb{T}^{m}) \hookrightarrow \text{Lip}_{p, \tau, \theta_{1}}^{(\boldsymbol{\alpha}^{(1)}, -\mathbf{b}^{(1)})}(\mathbb{T}^{m}) .
\end{equation*}
Let $\alpha_{j}^{(0)}=\alpha_{j}^{(1)}=\alpha_{j}$ for $j=1,\ldots, m$. First, suppose that $0<\theta_{1}<\theta_{0}\leqslant \infty$ and $b_{j}^{(1)}-\frac{1}{\theta_{1}}> b_{j}^{(0)}-\frac{1}{\theta_{0}}$ for $j=1,\ldots, m$. Applying H\"{o}lder's inequality with $\eta=\frac{\theta_{0}}{\theta_{1}}$, $\frac{1}{\eta}+\frac{1}{\eta^{'}}=1$, we obtain
\begin{equation}\label{eq5.3} 
   \begin{gathered}
   \Biggl(\int_{0}^{1}...\int_{0}^{1}\omega_{\boldsymbol{\alpha}}^{\theta}(f, \mathbf{t})_{p, \tau}\Bigl(\prod_{j=1}^{m}t_{j}^{-\alpha_{j}} (1-\log t_{j})^{- b_{j}^{(1)}}\Bigr)^{\theta_{1}} \prod_{j=1}^{m}\frac{dt_{j}}{t_{j}}\Biggr)^{\frac{1}{\theta_{1}}}
\\
\ll \Biggl(\int_{0}^{1}...\int_{0}^{1}\omega_{\boldsymbol{\alpha}}^{\theta}(f, \mathbf{t})_{p, \tau}\Bigl(\prod_{j=1}^{m}t_{j}^{-\alpha_{j}} (1-\log t_{j})^{- b_{j}^{(0)}}\Bigr)^{\theta_{0}} \prod_{j=1}^{m}\frac{dt_{j}}{t_{j}}\Biggr)^{\frac{1}{\theta_{0}}}
\\
\times
\Biggl(\int_{0}^{1}...\int_{0}^{1}\prod_{j=1}^{m}(1-\log t_{j})^{(b_{j}^{(0)}- b_{j}^{(1)})\theta_{1}\eta^{'}}\frac{dt_{j}}{t_{j}} \Biggr)^{\frac{1}{\theta_{1}\eta^{'}}}.
\end{gathered}
\end{equation}
Since $b_{j}^{(1)}-\frac{1}{\theta_{1}}> b_{j}^{(0)}-\frac{1}{\theta_{0}}$, then $1+(b_{j}^{(0)}- b_{j}^{(1)})\theta_{1}\eta^{'}<0$. Therefore
\begin{equation}\label{eq5.4}
\int_{0}^{1}(1-\log t_{j})^{(b_{j}^{(0)}- b_{j}^{(1)})\theta_{1}\eta^{'}}\frac{dt_{j}}{t_{j}}<+\infty, \, \, j=1,\ldots, m.
\end{equation}
 From the inequalities \eqref{eq5.3} and \eqref{eq5.4} it follows that
$\text{Lip}_{p, \tau, \theta_{0}}^{(\boldsymbol{\alpha}^{(0)}, -\mathbf{b}^{(0)})}(\mathbb{T}^{m})
\hookrightarrow \text{Lip}_{p, \tau, \theta_{1}}^{(\boldsymbol{\alpha}^{(1)}, -\mathbf{b}^{(1)})}(\mathbb{T}^{m})$ 
in the case $0<\theta_{1}<\theta_{0}\leqslant \infty$ and $b_{j}^{(1)}-\frac{1}{\theta_{1}}> b_{j}^{(0)}-\frac{1}{\theta_{0}}$ for $j=1,\ldots, m$.

Let $\alpha_{j}^{(0)}=\alpha_{j}^{(1)}=\alpha_{j}$, $0<\theta_{0}\leqslant\theta_{1}\leqslant \infty$ and $b_{j}^{(1)}-\frac{1}{\theta_{1}}\geqslant b_{j}^{(0)}-\frac{1}{\theta_{0}}$ for $j=1,\ldots, m$. Then by Jensen's inequality \cite{27} we have
\begin{equation}\label{eq5.5} 
   \begin{gathered}
   \Biggl(\sum\limits_{\nu_{m} =0}^{\infty}...\sum\limits_{\nu_{1} =0}^{\infty} \prod_{j=1}^{m}2^{\nu_{j}\alpha_{j}\theta_{1}}(\nu_{j} + 1)^{-\theta_{1} b_{j}^{(1)}} \omega_{\boldsymbol{\alpha}}^{\theta_{1}}\bigl(f,\frac{1}{2^{\nu_{1}}},...,\frac{1}{2^{\nu_{m}}}\bigr)_{p, \tau} \Biggr)^{\frac{1}{\theta_{1}}}
\\
\leqslant \Biggl(\sum\limits_{\nu_{m} =0}^{\infty}...\sum\limits_{\nu_{1} =0}^{\infty} \prod_{j=1}^{m}2^{\nu_{j}\alpha_{j}\theta_{0}}(\nu_{j} + 1)^{-\theta_{0} b_{j}^{(0)}}(\nu_{j} + 1)^{\theta_{0}(b_{j}^{(0)} - b_{j}^{(1)})} \omega_{\boldsymbol{\alpha}}^{\theta_{0}}\bigl(f,\frac{1}{2^{\nu_{1}}},...,\frac{1}{2^{\nu_{m}}}\bigr)_{p, \tau} \Biggr)^{\frac{1}{\theta_{0}}}.
\end{gathered}
\end{equation}
Since $0<\theta_{0}\leqslant\theta_{1}\leqslant \infty$ and $b_{j}^{(1)}-\frac{1}{\theta_{1}}\geqslant b_{j}^{(0)}-\frac{1}{\theta_{0}}$ for $j=1,\ldots, m$, then $b_{j}^{(0)}-b_{j}^{(1)}\leqslant \frac{1}{\theta_{0}} - \frac{1}{\theta_{1}}\leqslant 0$ for $j=1,\ldots, m$.
Therefore, it follows from \eqref{eq5.5} that
$\Omega_{\boldsymbol{\alpha}}^{\frac{1}{\theta_{1}}}(f)_{p, \tau, \theta_{1}}\ll \Omega_{\boldsymbol{\alpha}}^{\frac{1}{\theta_{0}}}(f)_{p, \tau, \theta_{0}}$. 
Therefore, according to Lemma 2.7 we will have that
$\text{Lip}_{p, \tau, \theta_{0}}^{(\boldsymbol{\alpha}^{(0)}, -\mathbf{b}^{(0)})}(\mathbb{T}^{m})
\hookrightarrow \text{Lip}_{p, \tau, \theta_{1}}^{(\boldsymbol{\alpha}^{(1)}, -\mathbf{b}^{(1)})}(\mathbb{T}^{m})$ in the case 
$\alpha_{j}^{(0)}=\alpha_{j}^{(1)}=\alpha_{j}$,  $0<\theta_{0}\leqslant\theta_{1}\leqslant \infty$ and $b_{j}^{(1)}-\frac{1}{\theta_{1}}\geqslant b_{j}^{(0)}-\frac{1}{\theta_{0}}$ for $j=1,\ldots, m$.
This proves Theorem 1.4.  
\hfill $\Box$

\textbf{Proof of Theorem 1.5.}
Due to the decreasing of the sequence $\{Y_{\boldsymbol{\nu}}(f)_{p, \tau}\}$ with respect to each index and the property of the mixed modulus of smoothness from Lemma 2.5 for $\alpha_{j}=\alpha_{j}^{(1)}$, $p=p_{1}$, $\tau=\tau_{1}$, we obtain
\begin{equation}\label{eq5.6}
\omega_{\boldsymbol{\alpha}^{(1)}}\Bigl(f,\frac{1}{2^{\nu_{1}}},...,\frac{1}{2^{\nu_{m}}}\Bigr)_{p_{1}, \tau_{1}} \ll
 \prod_{j=1}^{m}2^{-\nu_{j}\alpha_{j}^{(1)}}\sum\limits_{l_{1}=0}^{\nu_{1}} \ldots \sum\limits_{l_{m}=0}^{\nu_{m}} \prod_{j=1}^{m}2^{l_{j}\alpha_{j}^{(1)}}Y_{2^{l_{1}},\ldots,2^{l_{m}}}(f)_{p_{1}, \tau_{1}}
\end{equation}
for the function  $f\in L_{p_{1}, \tau_{1}}(\mathbb{T}^{m})$. 
On the other hand, it can be proven that
   \begin{equation}\label{eq5.7}
Y_{2^{l_{1}},\ldots,2^{l_{m}}}(f)_{p_{1}, \tau_{1}}\ll \prod_{j=1}^{m}2^{l_{j}(\frac{1}{p_{0}}-\frac{1}{p_{1}})}Y_{2^{l_{1}},\ldots,2^{l_{m}}}(f)_{p_{0}, \tau_{0}} + \sum\limits_{k_{1}=l_{1}}^{\infty} \ldots \sum\limits_{k_{m}=l_{m}}^{\infty}\prod_{j=1}^{m}2^{k_{j}(\frac{1}{p_{0}}-\frac{1}{p_{1}})}Y_{2^{k_{1}},\ldots,2^{k_{m}}}(f)_{p_{0}, \tau_{0}}.
\end{equation}
Now, using Lemma 2.4 from inequalities \eqref{eq5.6} and \eqref{eq5.7}, we obtain
\begin{equation*}
   \begin{gathered}
   \omega_{\boldsymbol\alpha^{(1)}}\Bigl(f,\frac{1}{2^{\nu_{1}}},...,\frac{1}{2^{\nu_{m}}}\Bigr)_{p_{1}, \tau_{1}} \ll \biggl(\prod_{j=1}^{m}2^{-\nu_{j}\alpha_{j}^{(1)}}\sum\limits_{l_{1}=0}^{\nu_{1}} \ldots \sum\limits_{l_{m}=0}^{\nu_{m}} \prod_{j=1}^{m}2^{l_{j}(\alpha_{j}^{(1)}+\frac{1}{p_{0}}-\frac{1}{p_{1}})}Y_{2^{l_{1}},\ldots,2^{l_{m}}}(f)_{p_{0}, \tau_{0}}
\\
 + \prod_{j=1}^{m}2^{-\nu_{j}\alpha_{j}^{(1)}}\sum\limits_{l_{1}=0}^{\nu_{1}} \ldots \sum\limits_{l_{m}=0}^{\nu_{m}} \prod_{j=1}^{m}2^{l_{j}\alpha_{j}^{(1)}}\sum\limits_{k_{1}=l_{1}}^{\infty} \ldots \sum\limits_{k_{m}=l_{m}}^{\infty}\prod_{j=1}^{m}2^{k_{j}(\frac{1}{p_{0}}-\frac{1}{p_{1}})}Y_{2^{k_{1}},\ldots,2^{k_{m}}}(f)_{p_{0}, \tau_{0}} \biggr)
 \end{gathered}
\end{equation*}
 \begin{equation}\label{eq5.8}
 \ll \sum\limits_{k_{1}=\nu_{1}}^{\infty} \ldots \sum\limits_{k_{m}=\nu_{m}}^{\infty}\prod_{j=1}^{m}2^{k_{j}(\frac{1}{p_{0}}-\frac{1}{p_{1}})} \omega_{\boldsymbol{\alpha}^{(0)}}\Bigl(f,\frac{1}{2^{k_{1}}},...,\frac{1}{2^{k_{m}}}\Bigr)_{p_{0}, \tau_{0}}
\end{equation}
for the function $f\in L_{p_{0}, \tau_{0}}(\mathbb{T}^{m})$.
According to inequality \eqref{eq5.8}, we will have
\begin{equation}\label{eq5.9} 
   \begin{gathered}
   \sum\limits_{\nu_{m} =0}^{\infty}...\sum\limits_{\nu_{1} =0}^{\infty} \prod_{j=1}^{m}2^{\nu_{j}\alpha_{j}^{(1)}\theta}(\nu_{j} + 1)^{-\theta b_{j}} \omega_{\boldsymbol{\alpha}^{(1)}}^{\theta}\bigl(f, \frac{1}{2^{\nu_{1}}},...,\frac{1}{2^{\nu_{m}}}\bigr)_{p_{1}, \tau_{1}}
\\
\ll \sum\limits_{\boldsymbol{\nu} \in \mathbb{Z}_{+}^{m}} \prod_{j=1}^{m}2^{\nu_{j}\alpha_{j}^{(1)}\theta}(\nu_{j} + 1)^{-\theta b_{j}} \biggl(\sum\limits_{k_{1}=\nu_{1}}^{\infty} \ldots \sum\limits_{k_{m}=\nu_{m}}^{\infty}\prod_{j=1}^{m}2^{k_{j}(\frac{1}{p_{0}}-\frac{1}{p_{1}})} \omega_{\boldsymbol{\alpha}^{(0)}}\Bigl(f,\frac{1}{2^{k_{1}}},...,\frac{1}{2^{k_{m}}}\Bigr)_{p_{0}, \tau_{0}}\biggr)^{\theta}.
\end{gathered}
\end{equation}
Since $\alpha_{j}^{(1)}>0$ for $j=1,\ldots, m$, then
 \begin{equation*}
\sum\limits_{k_{j} =0}^{\nu_{j}}2^{k_{j}\alpha_{j}^{(1)}\theta}(k_{j} + 1)^{-\theta b_{j}}\ll 2^{\nu_{j}\alpha_{j}^{(1)}\theta}(\nu_{j} + 1)^{-\theta b_{j}}.
\end{equation*}
Therefore, using Hardy's generalized inequality (see \cite[Lemma B.2]{15},  \cite[Lemma 3.2]{32}) from \eqref{eq5.9}, we obtain
\begin{equation*} 
   \begin{gathered}
   \sum\limits_{\nu_{m} =0}^{\infty}...\sum\limits_{\nu_{1} =0}^{\infty} \prod_{j=1}^{m}2^{\nu_{j}\alpha_{j}^{(1)}\theta}(\nu_{j} + 1)^{-\theta b_{j}} \omega_{\boldsymbol{\alpha}^{(1)}}^{\theta}\bigl(f, \frac{1}{2^{\nu_{1}}},...,\frac{1}{2^{\nu_{m}}}\bigr)_{p_{1}, \tau_{1}}
\\
\ll \sum\limits_{\nu_{m} =0}^{\infty}...\sum\limits_{\nu_{1} =0}^{\infty} \prod_{j=1}^{m}2^{\nu_{j}(\alpha_{j}^{(1)}+\frac{1}{p_{0}}-\frac{1}{p_{1}})\theta}(\nu_{j} + 1)^{-\theta b_{j}}\omega_{\boldsymbol{\alpha}^{(0)}}^{\theta}\Bigl(f,\frac{1}{2^{\nu_{1}}},...,\frac{1}{2^{\nu_{m}}}\Bigr)_{p_{0}, \tau_{0}}.
\end{gathered}
\end{equation*}
Since $\alpha_{j}^{(1)}-\frac{1}{p_{1}}=\alpha_{j}^{(1)}-\frac{1}{p_{0}}$ for $j=1,\ldots, m$, this inequality can be rewritten as follows
\begin{equation*} 
   \begin{gathered}
\sum\limits_{\nu_{m} =0}^{\infty}...\sum\limits_{\nu_{1} =0}^{\infty} \prod_{j=1}^{m}2^{\nu_{j}\alpha_{j}^{(1)}\theta}(\nu_{j} + 1)^{-\theta b_{j}} \omega_{\boldsymbol{\alpha}^{(1)}}^{\theta}\bigl(f, \frac{1}{2^{\nu_{1}}},...,\frac{1}{2^{\nu_{m}}}\bigr)_{p_{1}, \tau_{1}}
\\
\ll \sum\limits_{\nu_{m} =0}^{\infty}...\sum\limits_{\nu_{1} =0}^{\infty} \prod_{j=1}^{m}2^{\nu_{j}\alpha_{j}^{(0)}\theta}(\nu_{j} + 1)^{-\theta b_{j}}\omega_{\boldsymbol{\alpha}^{(0)}}^{\theta}\Bigl(f,\frac{1}{2^{\nu_{1}}},...,\frac{1}{2^{\nu_{m}}}\Bigr)_{p_{0}, \tau_{0}}
\end{gathered}
\end{equation*}
for the function $f\in L_{p_{0}, \tau_{0}}(\mathbb{T}^{m})$.
Therefore, by Lemma 2.7, the inclusion 
\begin{equation*}
\text{Lip}_{p_{0}, \tau_{0}, \theta}^{(\boldsymbol{\alpha}^{(0)}, -\mathbf{b})}(\mathbb{T}^{m}) \hookrightarrow \text{Lip}_{p_{1}, \tau_{1}, \theta}^{(\boldsymbol{\alpha}^{(1)}, -\mathbf{b})}(\mathbb{T}^{m})
\end{equation*}
 holds. Theorem 1.5 is proved.
\hfill $\Box$

\setcounter{equation}{0}
\setcounter{lemma}{0}
\setcounter{theorem}{0}

\section{Proofs of Theorems 1.6Ц1.8}

Further, we will often use the following statement:

\begin{lemma}\label{lem6.1}
Let $1<p, \tau < \infty$, $\delta\in \mathbb{R}$ and
\begin{equation*}
g(\mathbf{x})=\sum\limits_{s_{j_{0}} =1}^{\infty}2^{-s_{j_{0}}(\alpha_{j_{0}}+1-\frac{1}{p})}(s_{j_{0}} +1)^{-\delta}\prod_{j\neq j_{0}}\cos x_{j}\sum\limits_{k_{j_{0}}=2^{s_{j_{0}}-1}}^{2^{s_{j_{0}}}-1}\cos k_{j_{0}}x_{j_{0}}.
\end{equation*}
Then 
\begin{equation*}
\Biggl\|\sum\limits_{s_{m}=2^{l_{m}}}^{2^{l_{m}+1}-1}...\sum\limits_{s_{1}=2^{l_{1}}}^{2^{l_{1}+1}-1}\prod_{j=1}^{m}2^{s_{j}\alpha_{j}}\delta_{\mathbf{s}}(g)\Biggr\|_{p, \tau} \asymp 2^{-l_{j_{0}}(\delta-\frac{1}{\tau})}.
\end{equation*}
\end{lemma}
\proof
First, let us choose a number $p_{2}>p>1$. Then, according to Theorem 3.2 in \cite{2} (with $q$ replaced by $p$ and $\lambda=p_{2}$), we will have
\begin{multline}\label{eq6.1}  
\Biggl\|\sum\limits_{s_{m}=2^{l_{m}}}^{2^{l_{m}+1}-1}...\sum\limits_{s_{1}=2^{l_{1}}}^{2^{l_{1}+1}-1}\prod_{j=1}^{m}2^{s_{j}\alpha_{j}}\delta_{\mathbf{s}}(f_{1})\Biggr\|_{p, \tau}\gg 
\Biggl(\sum\limits_{s_{m}=2^{l_{m}}}^{2^{l_{m}+1}-1}...\sum\limits_{s_{1}=2^{l_{1}}}^{2^{l_{1}+1}-1}\prod_{j=1}^{m}2^{s_{j}(\alpha_{j}+ \frac{1}{p_{2}}-\frac{1}{p})\tau} \|\delta_{\mathbf{s}}(f_{1})\|_{p_{2}}^{\tau}\Biggr)^{\frac{1}{\tau}}
\\ 
=C\Biggl(\sum\limits_{s_{j_{0}}=2^{l_{j_{0}}}}^{2^{l_{j_{0}}+1}-1}2^{s_{j_{0}}(\frac{1}{p_{2}}-\frac{1}{p})\tau}2^{-s_{j_{0}}(1-\frac{1}{p})\tau}(s_{j_{0}} +1)^{-\delta\tau}\Bigl\|\sum\limits_{k_{j_{0}}=2^{s_{j_{0}}-1}}^{2^{s_{j_{0}}}-1}\cos k_{j_{0}}x_{j_{0}} \Bigr\|_{p_{2}}^{\tau} \Biggr)^{\frac{1}{\tau}}
\\
\gg \Biggl(\sum\limits_{s_{j_{0}}=2^{l_{j_{0}}}}^{2^{l_{j_{0}}+1}-1}(s_{j_{0}} +1)^{-\delta\tau} \Biggr)^{\frac{1}{\tau}} \gg 2^{-l_{j_{0}}(\delta-\frac{1}{\tau})}. 
\end{multline}
 Let $1< p_{3}<p<\infty$. Then, by Theorem 3.1 in \cite{2} and the estimate of the Dirichlet kernel norm in Lorentz space in \cite{2}, we have
\begin{multline}\label{eq6.2} 
   \Biggl\|\sum\limits_{s_{m}=2^{l_{m}}}^{2^{l_{m}+1}-1}...\sum\limits_{s_{1}=2^{l_{1}}}^{2^{l_{1}+1}-1}\prod_{j=1}^{m}2^{s_{j}\alpha_{j}}\delta_{\mathbf{s}}(f_{2})\Biggr\|_{p, \tau}\ll
\Biggl(\sum\limits_{s_{m}=2^{l_{m}}}^{2^{l_{m}+1}-1}...\sum\limits_{s_{1}=2^{l_{1}}}^{2^{l_{1}+1}-1}\prod_{j=1}^{m}2^{s_{j}(\alpha_{j}+ \frac{1}{p_{3}}-\frac{1}{p})\tau} \|\delta_{\mathbf{s}}(f_{2})\|_{p_{3}}^{\tau}\Biggr)^{\frac{1}{\tau}}
\\
=C\Biggl(\sum\limits_{s_{j_{0}}=2^{l_{j_{0}}}}^{2^{l_{j_{0}}+1}-1}2^{s_{j_{0}}(\frac{1}{p_{3}}-\frac{1}{p})\tau}2^{-s_{j_{0}}(1-\frac{1}{p})\tau}(s_{j_{0}} +1)^{-\delta\tau}\Bigl\|\sum\limits_{k_{j_{0}}=2^{s_{j_{0}}-1}}^{2^{s_{j_{0}}}-1}\cos k_{j_{0}}x_{j_{0}} \Bigr\|_{p_{3}}^{\tau} \Biggr)^{\frac{1}{\tau}}
\\
\ll \Biggl(\sum\limits_{s_{j_{0}}=2^{l_{j_{0}}}}^{2^{l_{j_{0}}+1}-1}(s_{j_{0}} +1)^{-\delta\tau} \Biggr)^{\frac{1}{\tau}} \ll 2^{l_{j_{0}}(\frac{1}{\tau}-\delta)}. 
\end{multline}
\hfill $\Box$

\textbf{Proof of Theorem 1.6.}
 Let $1<\tau<\theta \leqslant\infty$, i.e.,  $\min\{\tau, \theta\}=\tau$. Let's choose a number Let's choose a number $\delta$ such that $-b_{j_{0}}+\frac{1}{\tau}+\frac{1}{\theta}-\varepsilon<\delta<-b_{j_{0}}+\frac{1}{\tau}+\frac{1}{\theta}$.
Let's consider the function 
\begin{equation*}
f_{1}(\mathbf{x}) = \sum\limits_{s_{j_{0}} =1}^{\infty}2^{-s_{j_{0}}(\alpha_{j_{0}}+1-\frac{1}{p})}(s_{j_{0}} +1)^{-\delta}\prod_{j\neq j_{0}}\cos x_{j}\sum\limits_{k_{j_{0}}=2^{s_{j_{0}}-1}}^{2^{s_{j_{0}}}-1}\cos k_{j_{0}}x_{j_{0}} .
\end{equation*}
According to the Dirichlet kernel norm in Lorentz space \cite{2}, we have
\begin{equation}\label{eq6.3}  
   \begin{gathered}
   \sum\limits_{\mathbf{s}\in \mathbb{Z}_{+}^{m}}\prod_{j=1}^{m}2^{s_{j}\beta\theta}(1 + s_{j})^{(v_{j}-\varepsilon)\theta}\|\delta_{\mathbf{s}}(f_{1})\|_{p_{0}, \tau_{0}}^{\theta}\ll \sum\limits_{s_{j_{0}} =1}^{\infty}2^{s_{j_{0}}\beta\theta}(1 + s_{j_{0}})^{(v_{j}-\varepsilon)\theta}
\\
\times
2^{-s_{j_{0}}(\alpha_{j_{0}}+1-\frac{1}{p})\theta}(s_{j_{0}} +1)^{-\delta\theta}
 \Bigl\|\sum\limits_{k_{j_{0}}=2^{s_{j_{0}}-1}}^{2^{s_{j_{0}}}-1}\cos k_{j_{0}}x_{j_{0}} \Bigr\|_{p_{0}, \tau_{0}}^{\theta} \ll \sum\limits_{s_{j_{0}} =1}^{\infty}(s_{j_{0}} +1)^{(v_{j_{0}}-\varepsilon-\delta)\theta}
\end{gathered}
\end{equation}
Since $(v_{j_{0}}-\varepsilon-\delta)\theta>1$ is chosen, the series
\begin{equation*}
\sum\limits_{s_{j_{0}} =1}^{\infty}(s_{j_{0}} +1)^{(v_{j_{0}}-\varepsilon-\delta)\theta}
\end{equation*}
converges. 
Hence, from \eqref{eq6.3} we obtain that the function $C_{1}f_{1}\in S_{p_{0}, \tau_{0},  \theta}^{\boldsymbol{\beta}, \mathbf{v}-\varepsilon\mathbf{e}_{j_{0}}}B(\mathbb{T}^{m})$, for some $C_{1}> 0$.
 Therefore, from Theorem 1.1 and formula \eqref{eq6.1}, it follows that
\begin{equation}\label{eq6.4} 
\|f\|_{\text{Lip}_{p, \tau, \theta}^{(\boldsymbol{\alpha}, -\mathbf{b})}} \gg \biggl(\sum\limits_{l_{j_{0}} =0}^{\infty}\frac{1}{2^{l_{j_{0}}(\delta+b_{j_{0}}-\frac{1}{\theta}-\frac{1}{\tau})\theta}} \biggr)^{\frac{1}{\theta}}.
\end{equation}
Since $\delta+b_{j_{0}}-\frac{1}{\theta}-\frac{1}{\tau}<0$, then the right-hand side of \eqref{eq6.4} diverges. Therefore, from \eqref{eq6.4} we obtain that $f_{1}\notin \text{Lip}_{p, \tau, \theta}^{(\boldsymbol{\alpha}, -\mathbf{b})}$. The first statement is proven for $1<\tau<\theta \leqslant\infty$.

 Let us prove the second statement in Theorem 1.6.
 Let $\mathbf{u}=(u_{1}, \ldots ,u_{m})$, $\boldsymbol{\mu}=(\mu_{1}, \ldots ,\mu_{m})$, $u_{j}=\alpha_{j}+\frac{1}{p_{1}}-\frac{1}{p}$, $\mu_{j}=-b_{j}+\frac{1}{\max\{\tau, \theta\}}$ for  $j=1,\ldots, m$.
Let $0<\theta<\tau <\infty$, i.e.,  $\max\{\tau, \theta\}=\tau$.
Let us choose a number  $\delta$ such that $-b_{j_{0}}+\frac{1}{\tau}+\frac{1}{\theta}<\delta<\min\{\frac{1}{\tau}, -b_{j_{0}}+\frac{1}{\tau}+\frac{1}{\theta}+\varepsilon\}$.
Let's consider the function 
\begin{equation*}
f_{2}(\mathbf{x}) = \sum\limits_{s_{j_{0}} =1}^{\infty}2^{-s_{j_{0}}(\alpha_{j_{0}}+1-\frac{1}{p})}(s_{j_{0}} +1)^{-\delta}\prod_{j\neq j_{0}}\cos x_{j}\sum\limits_{k_{j_{0}}=2^{s_{j_{0}}-1}}^{2^{s_{j_{0}}}-1}\cos k_{j_{0}}x_{j_{0}} .
\end{equation*}
From Theorem 1.1 and Lemma 6.1, it follows that
\begin{equation*} 
\|f_{2}\|_{\text{Lip}_{p, \tau, \theta}^{(\boldsymbol{\alpha}, -\mathbf{b})}} \ll \biggl(\sum\limits_{l_{j_{0}} =0}^{\infty}\frac{1}{2^{l_{j_{0}}(\delta+b_{j_{0}}-\frac{1}{\theta}-\frac{1}{\tau})\theta}} \biggr)^{\frac{1}{\theta}}.
\end{equation*}
Since $\delta+b_{j_{0}}-\frac{1}{\theta}-\frac{1}{\tau}>0$, the last series converges.
Hence, $f_{2}\in \text{Lip}_{p, \tau, \theta}^{(\boldsymbol\alpha, -\mathbf{b})}(\mathbb{T}^{m})$.  
According to the estimate of the Dirichlet norm in Lorentz space in \cite{2}, we have
\begin{equation}\label{eq6.5} 
   \begin{gathered}
   \sum\limits_{\mathbf{s}\in \mathbb{Z}_{+}^{m}}\prod_{j=1}^{m}2^{s_{j}u_{j}\theta}(1 + s_{j})^{(\mu_{j}+\varepsilon)\theta}\|\delta_{\mathbf{s}}(f_{2})\|_{p_{1}, \tau_{1}}^{\theta}\gg \sum\limits_{s_{j_{0}} =1}^{\infty}2^{s_{j_{0}}u_{j_{0}}\theta}(1 + s_{j_{0}})^{(\mu_{j}+\varepsilon)\theta}
\\
\times
2^{-s_{j_{0}}(\alpha_{j_{0}}+1-\frac{1}{p})\theta}(s_{j_{0}} +1)^{-\delta\theta}
 \Bigl\|\sum\limits_{k_{j_{0}}=2^{s_{j_{0}}-1}}^{2^{s_{j_{0}}}-1}\cos k_{j_{0}}x_{j_{0}} \Bigr\|_{p_{1}, \tau_{1}}^{\theta} \gg \sum\limits_{s_{j_{0}} =1}^{\infty}\frac{1}{(s_{j_{0}} +1)^{(b_{j_{0}}+\delta-\frac{1}{\tau}-\varepsilon)\theta}}.
\end{gathered}
\end{equation}
Since $\delta< -b_{j_{0}}+\frac{1}{\tau}+\frac{1}{\theta}+\varepsilon$, then 
$(b_{j_{0}}+\delta-\frac{1}{\tau}-\varepsilon)\theta<1$. This means that the series on the right-hand side of inequality \eqref{eq6.5} diverges.
Hence, from \eqref{eq6.5} we obtain that $f_{2}\notin S_{p_{1}, \tau_{1}, \theta}^{\mathbf{u}, \boldsymbol{\mu}+\varepsilon\mathbf{e}_{j_{0}}}B(\mathbb{T}^{m})$ for $0<\theta<\tau <\infty$. 
\hfill $\Box$

\textbf{ Proof of Theorem 1.7.}
In Theorem 1.3 it is proved that
\begin{equation}\label{eq6.6}
S_{p, \tau, \min\{2, \tau, \theta\}}^{\boldsymbol{\beta}, \mathbf{v}}B(\mathbb{T}^{m}) \hookrightarrow \text{Lip}_{p, \tau, \theta}^{(\boldsymbol{\alpha}, -\mathbf{b})}(\mathbb{T}^{m}).
\end{equation}
  It is known that
$S_{p, \tau, r}^{\boldsymbol{\beta}, \mathbf{v}}B(\mathbb{T}^{m}) \hookrightarrow S_{p, \tau, \min\{2, \tau, \theta\}}^{\boldsymbol{\beta}, \mathbf{v}}B(\mathbb{T}^{m})$ for $0<r\leqslant \min\{\tau, \theta\}$. Therefore, from the inclusion \eqref{eq6.6} it follows that $S_{p, \tau, r} ^{\boldsymbol{\beta}, \mathbf{v}}B(\mathbb{T}^{m}) \hookrightarrow \text{Lip}_{p, \tau, \theta}^{(\boldsymbol{\alpha}, -\mathbf{b})}(\mathbb{T}^{m})$ for $0<r\leqslant \min\{\tau, \theta\}$. Sufficiency in \eqref{eq1.5} is proven.

Let's prove the necessity.
Suppose that there exists $r> \min\{\tau, \theta\}$ such that 
\begin{equation}\label{eq6.7}
S_{p_{0}, \tau_{0}, r}^{\boldsymbol{\beta}, \mathbf{v}}B(\mathbb{T}^{m}) \hookrightarrow \text{Lip}_{p, \tau, \theta}^{(\boldsymbol{\alpha}, -\mathbf{b})}(\mathbb{T}^{m}).
\end{equation}
Let $\min\{\tau, \theta\}=\tau$.
First, we choose a number $\delta$ such that $-b_{j_{0}}+\frac{1}{\theta}+\frac{1}{r}<\delta< -b_{j_{0}}+\frac{1}{\tau}+\frac{1}{\theta}$, and then consider the function
\begin{equation*}
f_{3}(\mathbf{x}) = \sum\limits_{s_{j_{0}} =1}^{\infty}2^{-s_{j_{0}}(\alpha_{j_{0}}+1-\frac{1}{p})}(s_{j_{0}} +1)^{-\delta}\prod_{j\neq j_{0}}\cos x_{j}\sum\limits_{k_{j_{0}}=2^{s_{j_{0}}-1}}^{2^{s_{j_{0}}}-1}\cos k_{j_{0}}x_{j_{0}}.
\end{equation*}
As in the proof of inequality \eqref{eq6.3}, it is easy to see that
\begin{equation}\label{eq6.8}
\sum\limits_{\mathbf{s}\in \mathbb{Z}_{+}^{m}}\prod_{j=1}^{m}2^{s_{j}\beta_{j}r}(1 + s_{j})^{v_{j}r}\|\delta_{\mathbf{s}}(f_{3})\|_{p_{0}, \tau_{0}}^{r}\ll  \sum\limits_{s_{j_{0}} =1}^{\infty}\frac{1}{(s_{j_{0}} +1)^{(\delta-v_{j_{0}})r}}.
\end{equation}
Since $(\delta-v_{j_{0}})r=(\delta+b_{j}-\frac{1}{\theta})r>1$, it follows from inequality \eqref{eq6.8} that the function $f_{3}\in S_{p_{0}, \tau_{0}, r}^{\boldsymbol{\beta}, \mathbf{v}}B(\mathbb{T}^{m})$.
Now, using Theorem 1.1 and Lemma 6.1, we have
\begin{equation}\label{eq6.9}  
\|f_{3}\|_{\text{Lip}_{p, \tau, \theta}^{(\boldsymbol{\alpha}, -\mathbf{b})}} \gg \biggl(\sum\limits_{l_{j_{0}} =0}^{\infty}\frac{1}{2^{l_{j_{0}}(\delta+b_{j_{0}}-\frac{1}{\theta}-\frac{1}{\tau})\theta}} \biggr)^{\frac{1}{\theta}}.
\end{equation}
Since $\delta+b_{j_{0}}-\frac{1}{\theta}-\frac{1}{\tau}< 0$, it follows from inequality \eqref{eq6.9} that $f_{3}\notin \text{Lip}_{p, \tau, \theta}^{(\boldsymbol{\alpha}, -\mathbf{b})}(\mathbb{T}^{m})$. This contradicts assumption \eqref{eq6.7}. Statement \eqref{eq1.5} is proven for $\min\{\tau, \theta\}=\tau$.

In Theorem 1.3 it is proved that
\begin{equation}\label{eq6.10}  
\text{Lip}_{p, \tau, \theta}^{(\boldsymbol{\alpha}, -\mathbf{b})} \hookrightarrow
S_{p_{1}, \tau_{1}, \max\{\tau, \theta\}}^{\mathbf{u}, \boldsymbol{\gamma}}B(\mathbb{T}^{m}).
\end{equation}
If  $r\geqslant \max\{\tau, \theta\}$, then $S_{p_{1}, \tau_{1}, \max\{\tau, \theta\}}^{\mathbf{u}, \boldsymbol{\gamma}}B(\mathbb{T}^{m}) \hookrightarrow S_{p_{1}, \tau_{1}, r}^{\mathbf{u}, \boldsymbol{\gamma}}B(\mathbb{T}^{m})$. Therefore, from inequality \eqref{eq6.10} it follows that $\text{Lip}_{p, \tau, \theta}^{(\boldsymbol{\alpha}, -\mathbf{b})} \hookrightarrow S_{p_{1}, \tau_{1}, r}^{\mathbf{u}, \boldsymbol{\gamma}}B(\mathbb{T}^{m})$.  In \eqref{eq1.6}, sufficiency is proven.

Let us prove the necessity of the condition $r\geqslant \max\{\tau, \theta\}$. Suppose that there exists $r< \max\{\tau, \theta\}=\tau$ such that
\begin{equation}\label{eq6.11}  
\text{Lip}_{p, \tau, \theta}^{(\boldsymbol\alpha, -\mathbf{b})} \hookrightarrow S_{p_{1}, \tau_{1}, r}^{\mathbf{u}, \boldsymbol{\gamma}}B(\mathbb{T}^{m}).
\end{equation}
Let $\max\{\tau, \theta\}=\tau$. Let us choose a number $\delta$ such that $-b_{j_{0}}+\frac{1}{\tau}+\frac{1}{\theta}<\delta< -b_{j_{0}}+\frac{1}{\theta}+\frac{1}{r}$ and consider the function
\begin{equation*}
f_{4}(\mathbf{x}) = \sum\limits_{s_{j_{0}} =1}^{\infty}2^{-s_{j_{0}}(\alpha_{j_{0}}+1-\frac{1}{p})}(s_{j_{0}} +1)^{-\delta}\prod_{j\neq j_{0}}\cos x_{j}\sum\limits_{k_{j_{0}}=2^{s_{j_{0}}-1}}^{2^{s_{j_{0}}}-1}\cos k_{j_{0}}x_{j_{0}}.
\end{equation*}
 Now, using Theorem 1.1 and Lemma 6.1, we have
\begin{equation}\label{eq6.12}  
\|f_{4}\|_{\text{Lip}_{p, \tau, \theta}^{(\boldsymbol{\alpha}, -\mathbf{b})}} \ll \biggl(\sum\limits_{l_{j_{0}} =0}^{\infty}\frac{1}{2^{l_{j_{0}}(\delta+b_{j_{0}}-\frac{1}{\theta}-\frac{1}{\tau})\theta}} \biggr)^{\frac{1}{\theta}}.
\end{equation}
Since $\delta+b_{j_{0}}-\frac{1}{\theta}-\frac{1}{\tau}> 0$, it follows from inequality \eqref{eq6.10} that $f_{4}\in \text{Lip}_{p, \tau, \theta}^{(\boldsymbol{\alpha}, -\mathbf{b})}(\mathbb{T}^{m})$.
Applying the estimate of the Dirichlet kernel norm in Lorentz space in \cite{2}, we have
\begin{equation}\label{eq6.13}  
\sum\limits_{\mathbf{s}\in \mathbb{Z}_{+}^{m}}\prod_{j=1}^{m}2^{s_{j}u_{j}r}(1 + s_{j})^{\gamma_{j}r}\|\delta_{\mathbf{s}}(f_{1})\|_{p_{1}, \tau_{1}}^{r}\gg  \sum\limits_{s_{j_{0}} =1}^{\infty}\frac{1}{(s_{j_{0}} +1)^{(\delta+b_{j_{0}}-\frac{1}{\theta})r}}.
\end{equation}
Since $(\delta+b_{j_{0}}-\frac{1}{\theta})r<1$, the series on the right-hand side of \eqref{eq6.13} diverges. Therefore, the function $f_{4}\notin S_{p_{1}, \tau_{1}, r}^{\mathbf{u}, \boldsymbol{\gamma}}B(\mathbb{T}^{m})$. This contradicts assumption \eqref{eq6.11} when $\max\{\tau, \theta\}=\tau$.
\hfill $\Box$

\textbf{Proof of Theorem 1.8.} In Theorem 1.2 it is proved that
  \begin{equation}\label{eq6.14}  
S_{p, \tau, \theta}^{\boldsymbol{\alpha}, -\mathbf{b}+\frac{1}{\min\{2, \tau, \theta\}}\mathbf{e}}B(\mathbb{T}^{m}) \hookrightarrow \text{Lip}_{p, \tau, \theta}^{(\boldsymbol{\alpha}, -\mathbf{b})}(\mathbb{T}^{m}).
\end{equation}
It is known that if $\min_{j=1,\ldots,m}\{\xi_{j}\}\geqslant \frac{1}{\min\{2, \tau, \theta\}}$, then $S_{p, \tau, \theta}^{\boldsymbol{\alpha}, -\mathbf{b}+\boldsymbol{\xi}}B(\mathbb{T}^{m}) \hookrightarrow S_{p, \tau, \theta}^{\boldsymbol{\alpha}, -\mathbf{b}+\frac{1}{\min\{2, \tau, \theta\}}\mathbf{e}}B(\mathbb{T}^{m})$. Therefore, it follows from the embedding \eqref{eq6.14} that $S_{p, \tau, \theta}^{\boldsymbol{\alpha}, -\mathbf{b}+\boldsymbol{\xi}}B(\mathbb{T}^{m}) \hookrightarrow \text{Lip}_{p, \tau, \theta}^{(\boldsymbol{\alpha}, -\mathbf{b})}(\mathbb{T}^{m})$ for $\min_{j=1,\ldots,m}\{\xi_{j}\}\geqslant \frac{1}{\min\{2, \tau, \theta\}}$. In \eqref{eq1.7}, sufficiency is proved.

Let's prove the necessity. Suppose that there exists a positive number $\xi_{j_{0}}<\frac{1}{\min\{2, \tau, \theta\}}$ and $\xi_{j}\geqslant\frac{1}{\min\{2, \tau, \theta\}}$ for $j\neq j_{0}$ such that 
\begin{equation}\label{eq6.15}  
S_{p, \tau, \theta}^{\boldsymbol{\alpha}, -\mathbf{b}+\boldsymbol{\xi}_{0}}B(\mathbb{T}^{m}) \hookrightarrow \text{Lip}_{p, \tau, \theta}^{(\boldsymbol{\alpha}, -\mathbf{b})}(\mathbb{T}^{m}),
\end{equation}
where  $\boldsymbol{\xi}_{0}=(\xi_{1},\ldots,\xi_{j_{0}-1}, \xi_{j_{0}}, \xi_{j_{0}+1,\ldots}, \xi_{m} )$. 

Let $\min\{2, \tau, \theta\}=\tau$. Then $\frac{1}{\tau}-\xi_{j_{0}}>0$. Therefore, $\frac{1}{\theta}<\frac{1}{\theta}+\frac{1}{\tau}-\xi_{j_{0}}$. Let us choose a number $\delta$ such that $\frac{1}{\theta}<\delta<\frac{1}{\theta}+\frac{1}{\tau}-\xi_{j_{0}}$. Consider the function
\begin{equation*} 
   \begin{gathered}
   f_{\boldsymbol{\xi}_{0}}(\mathbf{x}) = \sum\limits_{s_{j_{0}} =1}^{\infty}2^{-s_{j_{0}}(\alpha_{j_{0}}+1-\frac{1}{p})}(s_{j_{0}} +1)^{-(\xi_{j_{0}}-b_{j_{0}})-\delta}\sum\limits_{k_{j_{0}}=2^{s_{j_{0}}-1}}^{2^{s_{j_{0}}}-1}\cos k_{j_{0}}x_{j_{0}}
\\
\times\prod_{j\neq j_{0}}\sum\limits_{s_{j} =1}^{\infty}2^{-s_{j}(\alpha_{j}+1-\frac{1}{p})}(s_{j} +1)^{-(\xi_{j}-b_{j})-t}\sum\limits_{k_{j}=2^{s_{j}-1}}^{2^{s_{j}}-1}\cos k_{j}x_{j},
\end{gathered}
\end{equation*}
where the number $t>\frac{1}{\theta}$.
 Then, according to the estimate of the norm of the Dirichlet kernel in the Lorentz space in \cite{2} we have
\begin{equation}\label{eq6.16} 
\sum\limits_{\mathbf{s}\in \mathbb{Z}_{+}^{m}}\prod_{j=1}^{m}2^{s_{j}\alpha_{j}\theta}(1 + s_{j})^{(-b_{j}+\xi_{j})\theta}\|\delta_{\mathbf{s}}(f_{\boldsymbol{\xi}_{0}})\|_{p, \tau}^{\theta}\ll  \sum\limits_{s_{j_{0}} =1}^{\infty}\frac{1}{(s_{j_{0}} +1)^{\delta\theta}}\prod_{j\neq j_{0}}\sum\limits_{s_{j} =1}^{\infty}\frac{1}{(s_{j} +1)^{t\theta}}.
\end{equation}
Since $\delta\theta>1$ and $t\theta>1$, the series on the right-hand side of \eqref{eq6.16} converge. Therefore, the function $f_{\boldsymbol{\xi}_{0}}\in S_{p, \tau, \theta}^{\boldsymbol{\alpha}, -\mathbf{b}+\boldsymbol{\xi}_{0}}B(\mathbb{T}^{m})$.  Let us choose a number $p_{2}>p>1$. Then, according to Theorem 3.2 in \cite{2} (with $q$ replaced by $p$ and $\lambda=p_{2}$), we will have (see \eqref{eq6.1})
\begin{equation*} 
   \begin{gathered}
\Biggl\|\sum\limits_{s_{m}=2^{l_{m}}}^{2^{l_{m}+1}-1}...\sum\limits_{s_{1}=2^{l_{1}}}^{2^{l_{1}+1}-1}\prod_{j=1}^{m}2^{s_{j}\alpha_{j}}\delta_{\mathbf{s}}(f_{\boldsymbol{\xi}_{0}})\Biggr\|_{p, \tau}\gg \Biggl(\sum\limits_{s_{m}=2^{l_{m}}}^{2^{l_{m}+1}-1}...\sum\limits_{s_{1}=2^{l_{1}}}^{2^{l_{1}+1}-1}\prod_{j=1}^{m}2^{s_{j}(\alpha_{j}+ \frac{1}{p_{2}}-\frac{1}{p})\tau} \|\delta_{\mathbf{s}}(f_{\boldsymbol{\xi}_{0}})\|_{p_{2}}^{\tau}\Biggr)^{\frac{1}{\tau}}
\\
\gg \Biggl(\sum\limits_{s_{j_{0}}=2^{l_{j_{0}}}}^{2^{l_{j_{0}}+1}-1} (1 + s_{j_{0}})^{-(\xi_{j_{0}}-b_{j_{0}})\tau-\delta\tau} \prod_{j\neq j_{0}}\sum\limits_{s_{j}=2^{l_{j}}}^{2^{l_{j}+1}-1} (1 + s_{j})^{-(\xi_{j}-b_{j})\tau-t\tau} \Biggr)^{\frac{1}{\tau}}
\end{gathered}
\end{equation*}
\begin{equation}\label{eq6.17} 
\gg 2^{l_{j_{0}}(\frac{1}{\tau}-(\xi_{j_{0}}-b_{j_{0}}+\delta))}\prod_{j\neq j_{0}}2^{l_{j}(\frac{1}{\tau}-(\xi_{j}-b_{j}+t))}.
\end{equation}
Now, taking into account the inequality \eqref{eq6.17} we will have
\begin{equation}\label{eq6.18}   
   \begin{gathered}
   \sum\limits_{l_{m} =0}^{n_{m}}...\sum\limits_{l_{1} =0}^{n_{1}}
\prod_{j=1}^{m}2^{l_{j}(\frac{1}{\theta}-b_{j})\theta}
\Biggl\|\sum\limits_{s_{m}=2^{l_{m}}}^{2^{l_{m}+1}-1}...\sum\limits_{s_{1}=2^{l_{1}}}^{2^{l_{1}+1}-1}\prod_{j=1}^{m}2^{s_{j}\alpha_{j}}\delta_{\mathbf{s}}(f_{\boldsymbol{\xi}_{0}})\Biggr\|_{p, \tau}^{\theta}
\\
\gg \sum\limits_{l_{j_{0}} =0}^{n_{j_{0}}}\frac{1}{2^{l_{j_{0}}(\xi_{j_{0}}+\delta-\frac{1}{\tau}-\frac{1}{\theta})\theta}}\prod_{j\neq j_{0}}\sum\limits_{l_{j} =0}^{n_{j}} \frac{1}{2^{l_{j}(\xi_{j}+t-\frac{1}{\tau}-\frac{1}{\theta})\theta}}.
\end{gathered}
\end{equation}
Since $\xi_{j}-\frac{1}{\tau}\geqslant 0$ for $j\neq j_{0}$ and by choice $t-\frac{1}{\theta}>0$, then
\begin{equation}\label{eq6.19}  
\sum\limits_{l_{j} =0}^{\infty} \frac{1}{2^{l_{j}(\xi_{j}+t-\frac{1}{\tau}-\frac{1}{\theta})\theta}}
<\infty.
\end{equation}
On the other hand, by choosing the number $\delta$: $\xi_{j_{0}}+\delta-\frac{1}{\tau}-\frac{1}{\theta}<0$. Hence
\begin{equation}\label{eq6.20} 
\sum\limits_{l_{j_{0}} =0}^{\infty}\frac{1}{2^{l_{j_{0}}(\xi_{j_{0}}+\delta-\frac{1}{\tau}-\frac{1}{\theta})\theta}}=+\infty.
\end{equation}
Now, from \eqref{eq6.18}--\eqref{eq6.20}, according to Theorem 1.1, it follows that
the function $f_{\boldsymbol{\xi}_{0}}\notin \text{Lip}_{p, \tau, \theta}^{(\boldsymbol{\alpha}, -\mathbf{b})}(\mathbb{T}^{m})$. This contradicts assumption \eqref{eq6.15}. Statement \eqref{eq1.7} is proven for $\min\{\tau, \theta\}=\tau$.
 
 We will prove statement \eqref{eq1.8}. In Theorem 1.2 it is proved that
\begin{equation}\label{eq6.21} 
\text{Lip}_{p, \tau, \theta}^{(\boldsymbol{\alpha}, -\mathbf{b})}(\mathbb{T}^{m}) \hookrightarrow S_{p, \tau, \theta}^{\boldsymbol{\alpha}, -\mathbf{b}+\frac{1}{\max\{2, \tau, \theta\}}\mathbf{e}}B(\mathbb{T}^{m}).
\end{equation}
If $\max_{j=1,\ldots,m}\{\xi_{j}\}\leqslant \frac{1}{\max\{2, \tau, \theta\}}$, then $S_{p, \tau, \theta}^{\boldsymbol{\alpha}, -\mathbf{b}+\frac{1}{\max\{2, \tau, \theta\}}\mathbf{e}}B(\mathbb{T}^{m}) \hookrightarrow S_{p, \tau, \theta}^{\boldsymbol{\alpha}, -\mathbf{b}+\boldsymbol{\xi}}B(\mathbb{T}^{m})$. Therefore, from the inclusion \eqref{eq6.23}, it follows that $\text{Lip}_{p, \tau, \theta}^{(\boldsymbol{\alpha}, -\mathbf{b})}(\mathbb{T}^{m}) \hookrightarrow 
 S_{p, \tau, \theta}^{\boldsymbol{\alpha}, -\mathbf{b}+\boldsymbol{\xi}}B(\mathbb{T}^{m}) $ for $\max_{j=1,\ldots,m}\{\xi_{j}\}\leqslant \frac{1}{\max\{2, \tau, \theta\}}$. In \eqref{eq1.8}, sufficiency is proven.

Let us prove the necessity. Suppose that there exists a positive number $\xi_{j_{0}}>\frac{1}{\max\{2, \tau, \theta\}}$ and $\xi_{j}\leqslant\frac{1}{\max\{2, \tau, \theta\}}$ for $j\neq j_{0}$ such that
\begin{equation}\label{eq6.22}   
\text{Lip}_{p, \tau, \theta}^{(\boldsymbol{\alpha}, -\mathbf{b})}(\mathbb{T}^{m}) \hookrightarrow  S_{p, \tau, \theta}^{\boldsymbol{\alpha}, -\mathbf{b}+\boldsymbol{\xi}^{0}}B(\mathbb{T}^{m}),
\end{equation}
where  $\boldsymbol{\xi}^{0}=(\xi_{1},\ldots,\xi_{j_{0}-1}, \xi_{j_{0}}, \xi_{j_{0}+1,\ldots}, \xi_{m} )$.
 
Let $\max\{2, \tau, \theta\}=\tau$. Since  $\frac{1}{\tau}-\xi_{j_{0}}<0$, then $\frac{1}{\theta}+\frac{1}{\tau}-\xi_{j_{0}}<\frac{1}{\theta}$.  Let us choose a number $\delta$ such that $ \, \, \frac{1}{\theta}+\frac{1}{\tau}-\xi_{j_{0}}<\delta<\frac{1}{\theta}$. Now consider the function
\begin{equation*} 
   \begin{gathered}
   g_{\boldsymbol{\xi}^{0}}(\mathbf{x}) = \sum\limits_{s_{j_{0}} =1}^{\infty}2^{-s_{j_{0}}(\alpha_{j_{0}}+1-\frac{1}{p})}(s_{j_{0}} +1)^{-(\xi_{j_{0}}-b_{j_{0}})-\delta}\sum\limits_{k_{j_{0}}=2^{s_{j_{0}}-1}}^{2^{s_{j_{0}}}-1}\cos k_{j_{0}}x_{j_{0}}
\\
\times\prod_{j\neq j_{0}}\sum\limits_{s_{j} =1}^{\infty}2^{-s_{j}(\alpha_{j}+1-\frac{1}{p})}(s_{j} +1)^{-(\xi_{j}-b_{j})-t}\sum\limits_{k_{j}=2^{s_{j}-1}}^{2^{s_{j}}-1}\cos k_{j}x_{j},
\end{gathered}
\end{equation*}
 where the number   $t>\max_{j\neq j_{0}}\{\frac{1}{\theta}+ \frac{1}{\tau}-\xi_{j}\}$. 
Further, let us choose a number $1<p<p_{2}<\infty$. Then, by Theorem 3.1 in \cite{2} and according to the estimate of the Dirichlet kernel norm in Lebesgue space, we will have (see \eqref{eq6.2})
\begin{equation*} 
   \begin{gathered}
\Biggl\|\sum\limits_{s_{m}=2^{l_{m}}}^{2^{l_{m}+1}-1}...\sum\limits_{s_{1}=2^{l_{1}}}^{2^{l_{1}+1}-1}\prod_{j=1}^{m}2^{s_{j}\alpha_{j}}\delta_{\mathbf{s}}(g_{\boldsymbol{\xi}^{0}})\Biggr\|_{p, \tau}\ll \Biggl(\sum\limits_{s_{m}=2^{l_{m}}}^{2^{l_{m}+1}-1}...\sum\limits_{s_{1}=2^{l_{1}}}^{2^{l_{1}+1}-1}\prod_{j=1}^{m}2^{s_{j}(\alpha_{j}+ \frac{1}{p_{3}}-\frac{1}{p})\tau} \|\delta_{\mathbf{s}}(g_{\boldsymbol{\xi}^{0}})\|_{p_{3}}^{\tau}\Biggr)^{\frac{1}{\tau}}
\\
\ll \Biggl(\sum\limits_{s_{j_{0}}=2^{l_{j_{0}}}}^{2^{l_{j_{0}}+1}-1} (1 + s_{j_{0}})^{-(\xi_{j_{0}}-b_{j_{0}}+\delta)\tau} \prod_{j\neq j_{0}}\sum\limits_{s_{j}=2^{l_{j}}}^{2^{l_{j}+1}-1} (1 + s_{j})^{-(\xi_{j}-b_{j}+t)\tau} \Biggr)^{\frac{1}{\tau}}
\\
\ll 2^{l_{j_{0}}(\frac{1}{\tau}-(\xi_{j_{0}}-b_{j_{0}}+\delta))}\prod_{j\neq j_{0}}2^{l_{j}(\frac{1}{\tau}-(\xi_{j}-b_{j}+t))}.
\end{gathered}
\end{equation*}
Therefore 
\begin{equation}\label{eq6.23}  
   \begin{gathered}
   \sum\limits_{l_{m} =0}^{n_{m}}...\sum\limits_{l_{1} =0}^{n_{1}}
\prod_{j=1}^{m}2^{l_{j}(\frac{1}{\theta}-b_{j})\theta}
\Biggl\|\sum\limits_{s_{m}=2^{l_{m}}}^{2^{l_{m}+1}-1}...\sum\limits_{s_{1}=2^{l_{1}}}^{2^{l_{1}+1}-1}\prod_{j=1}^{m}2^{s_{j}\alpha_{j}}\delta_{\mathbf{s}}(g_{\boldsymbol{\xi}^{0}})\Biggr\|_{p, \tau}^{\theta}
\\
\ll \sum\limits_{l_{j_{0}} =0}^{n_{j_{0}}}\frac{1}{2^{l_{j_{0}}(\xi_{j_{0}}+\delta-\frac{1}{\tau}-\frac{1}{\theta})\theta}}\prod_{j\neq j_{0}}\sum\limits_{l_{j} =0}^{n_{j}} \frac{1}{2^{l_{j}(\xi_{j}+t-\frac{1}{\tau}-\frac{1}{\theta})\theta}}.
\end{gathered}
\end{equation}

Since $t>\frac{1}{\theta}+ \frac{1}{\tau}-\xi_{j}$ for $j\neq j_{0}$ and $\delta > \frac{1}{\theta}+\frac{1}{\tau}-\xi_{j_{0}}$, then
by Theorem 1.1, it follows from inequality \eqref{eq6.23} that the function $g_{\boldsymbol{\xi}^{0}}\in \text{Lip}_{p, \tau, \theta}^{(\boldsymbol{\alpha}, -\mathbf{b})}(\mathbb{T}^{m})$.
Further, according to the estimate of the norm of the Dirichlet kernel in the Lorentz space in \cite{2} we have
\begin{equation}\label{eq6.24} 
\sum\limits_{s_{m} =1}^{n_{m}}...\sum\limits_{s_{1} =1}^{n_{1}}\prod_{j=1}^{m}2^{s_{j}\alpha_{j}\theta}(1 + s_{j})^{(-b_{j}+\xi_{j})\theta}\|\delta_{\mathbf{s}}(g_{\boldsymbol{\xi}^{0}})\|_{p, \tau}^{\theta}\gg  \sum\limits_{s_{j_{0}} =1}^{n_{j_{0}}}\frac{1}{(s_{j_{0}} +1)^{\delta\theta}}\prod_{j\neq j_{0}}\sum\limits_{s_{j} =1}^{n_{j}}\frac{1}{(s_{j} +1)^{t\theta}}.
\end{equation}
Since $t\theta > 1$ and $\delta\theta<0$, it follows from inequality \eqref{eq6.26} that the function $g_{\boldsymbol{\xi}^{0}}\notin S_{p, \tau, \theta}^{\boldsymbol\alpha, -\mathbf{b}+\boldsymbol{\xi}^{0}}B(\mathbb{T}^{m})$ in the case $\max\{2, \tau, \theta\}=\tau$ and $\xi_{j_{0}}>\frac{1}{\max\{2, \tau, \theta\}}$.
Thus, the function $g_{\boldsymbol{\xi}^{0}}\in \text{Lip}_{p, \tau, \theta}^{(\boldsymbol{\alpha}, -\mathbf{b})}(\mathbb{T}^{m})$, but   $g_{\boldsymbol{\xi}^{0}}\notin S_{p, \tau, \theta}^{\boldsymbol{\alpha}, -\mathbf{b}+\boldsymbol{\xi}^{0}}B(\mathbb{T}^{m})$, if $\xi_{j_{0}}>\frac{1}{\max\{2, \tau, \theta\}}$ and $\max\{2, \tau, \theta\}=\tau$. This contradicts \eqref{eq6.22}. Statement \eqref{eq1.8} is proven $\max\{2, \tau, \theta\}=\tau$. 

 Let us prove statement \eqref{eq1.9}. Let $0<q\leqslant \min\{\tau, \theta\}$. Then, from relation \eqref{eq0.1} (see section УIntroductionФ), according to Jensen's inequality, we have 
 \begin{equation*}
 S_{p, \tau, q} ^{\boldsymbol\alpha, -\mathbf{b}+\frac{1}{\theta}\mathbf{e}}B(\mathbb{T}^{m})\hookrightarrow  S_{p, \tau, \min\{\tau,\theta\}}^{\boldsymbol\alpha, -\mathbf{b}+\frac{1}{\theta}\mathbf{e}}B(\mathbb{T}^{m}).
  \end{equation*}
   Therefore, from Theorem 1.2, it follows that 
 \begin{equation*}
 S_{p, \tau, q}^{\boldsymbol{\alpha}, -\mathbf{b}+\frac{1}{\theta}\mathbf{e}}B(\mathbb{T}^{m})\hookrightarrow \text{Lip}_{p, \tau, \theta}^{(\boldsymbol{\alpha}, -\mathbf{b})}(\mathbb{T}^{m})
 \end{equation*}
  for $q\in (0, \, \min\{\tau, \theta\})$.

 Let us prove the necessity of the condition $q\leqslant \min\{\tau, \theta\}$. Suppose that there exists a number $q> \min\{\tau, \theta\}$ such that 
\begin{equation}\label{eq6.25}  
S_{p, \tau, q}^{\boldsymbol{\alpha}, -\mathbf{b}+\frac{1}{\theta}\mathbf{e}}B(\mathbb{T}^{m})\hookrightarrow \text{Lip}_{p, \tau, \theta}^{(\boldsymbol{\alpha}, -\mathbf{b})}(\mathbb{T}^{m}).
\end{equation}
Let $\min\{\tau, \theta\}=\tau$. Then, by assumption, $\frac{1}{q}<\frac{1}{\tau}$. Let us choose a number $\delta$ such that $-b_{j_{0}}+\frac{1}{\theta}+\frac{1}{q}<\delta <-b_{j_{0}}+\frac{1}{\theta}+\frac{1}{\tau}$.
Now,  consider the function
\begin{equation*}
g_{1}(\mathbf{x}) = \sum\limits_{s_{j_{0}} =1}^{\infty}2^{-s_{j_{0}}(\alpha_{j_{0}}+1-\frac{1}{p})}(s_{j_{0}} +1)^{-\delta}\prod_{j\neq j_{0}}\cos x_{j}\sum\limits_{k_{j_{0}}=2^{s_{j_{0}}-1}}^{2^{s_{j_{0}}}-1}\cos k_{j_{0}}x_{j_{0}}.
\end{equation*}
 Then 
\begin{equation*}  
      \sum\limits_{\mathbf{s}\in \mathbb{Z}_{+}^{m}}\prod_{j=1}^{m}2^{s_{j}\alpha_{j}q}(1 + s_{j})^{(-b_{j}+\frac{1}{\theta})q}\|\delta_{\mathbf{s}}(f_{1})\|_{p, \tau}^{q}  \ll \sum\limits_{s_{j_{0}} =1}^{\infty}(s_{j_{0}} +1)^{(\delta+b_{j_{0}}-\frac{1}{\theta})q}<\infty.
\end{equation*}
Therefore, the function $g_{1}\in S_{p, \tau, q}^{\boldsymbol{\alpha}, -\mathbf{b}+\frac{1}{\theta}\mathbf{e}}B(\mathbb{T}^{m})$.
According to Lemma 6.1, the following inequality holds
\begin{equation}\label{eq6.26}  
   \Biggl\|\sum\limits_{s_{m}=2^{l_{m}}}^{2^{l_{m}+1}-1}...\sum\limits_{s_{1}=2^{l_{1}}}^{2^{l_{1}+1}-1}\prod_{j=1}^{m}2^{s_{j}\alpha_{j}}\delta_{\mathbf{s}}(g_{1})\Biggr\|_{p, \tau}\gg 2^{-l_{j_{0}}(\delta-\frac{1}{\tau})}.
\end{equation}
Now, using inequality \eqref{eq6.26}, we have
 \begin{equation}\label{eq6.27} 
   \sum\limits_{l_{m} =0}^{n_{m}}...\sum\limits_{l_{1} =0}^{n_{1}}
\prod_{j=1}^{m}2^{l_{j}(\frac{1}{\theta}-b_{j})\theta}
\Biggl\|\sum\limits_{s_{m}=2^{l_{m}}}^{2^{l_{m}+1}-1}...\sum\limits_{s_{1}=2^{l_{1}}}^{2^{l_{1}+1}-1}\prod_{j=1}^{m}2^{s_{j}\alpha_{j}}\delta_{\mathbf{s}}(g_{1})\Biggr\|_{p, \tau}^{\theta}
\\
\gg \sum\limits_{l_{j_{0}} =0}^{n_{j_{0}}}\frac{1}{2^{l_{j_{0}}(b_{j_{0}}+\delta-\frac{1}{\tau}-\frac{1}{\theta})\theta}}.
\end{equation}
Since $b_{j_{0}}+\delta-\frac{1}{\tau}-\frac{1}{\theta}<0$, it follows from \eqref{eq6.27} that 
$g_{1}\notin \text{Lip}_{p, \tau, \theta}^{(\boldsymbol{\alpha}, -\mathbf{b})}(\mathbb{T}^{m})$. This contradicts assumption \eqref{eq6.25}  when $\min\{\tau, \theta\}=\tau$.

In statement \eqref{eq1.10} the sufficiency of the condition $q\geqslant \max\{\tau, \theta\}$ for the indicated embedding follows from Theorem 1.2.
Suppose that for some $q\in (0, \, \max\{\tau, \theta\})$ the inclusion $\text{Lip}_{p, \tau, \theta}^{(\boldsymbol{\alpha}, -\mathbf{b})}(\mathbb{T}^{m}) \hookrightarrow S_{p, \tau, q} ^{\boldsymbol{\alpha}, -\mathbf{b}+\frac{1}{\theta}\mathbf{e}}B(\mathbb{T}^{m})$ holds.  
Let us choose a number $\delta$ such that $-b_{j_{0}}+\frac{1}{\theta}+\frac{1}{\tau}<\delta <-b_{j_{0}}+\frac{1}{\theta}+\frac{1}{q}$. Consider the function
\begin{equation*}
g_{2}(\mathbf{x}) = \sum\limits_{s_{j_{0}} =1}^{\infty}2^{-s_{j_{0}}(\alpha_{j_{0}}+1-\frac{1}{p})}(s_{j_{0}} +1)^{-\delta}\prod_{j\neq j_{0}}\cos x_{j}\sum\limits_{k_{j_{0}}=2^{s_{j_{0}}-1}}^{2^{s_{j_{0}}}-1}\cos k_{j_{0}}x_{j_{0}}.
\end{equation*}
Then (see Lemma 6.1 and \eqref{eq6.13}) the function $g_{2}\in \text{Lip}_{p, \tau, \theta}^{(\boldsymbol{\alpha}, -\mathbf{b})}(\mathbb{T}^{m})$, but $g_{2} \notin S_{p, \tau, q} ^{\boldsymbol{\alpha}, -\mathbf{b}+\frac{1}{\theta}\mathbf{e}}B(\mathbb{T}^{m})$. This is a contradiction.
\hfill $\Box$

\begin{flushleft}
Gabdolla Akishev,\\
 Kazakhstan Branch, \\
 Lomonosov Moscow University,\\
St. Kazhymukan, 11,\\
010010, Astana, Kazakhstan,\\
 E-mail: akishev g@mail.ru
\end{flushleft}
\end{document}